%% file: main.tex
\documentclass[article]{amsart}
\pdfoutput=1
\usepackage{montemplate}
\usepackage{standalone}

\title{An action of the Witt algebra on Khovanov--Rozansky homology}
\author{Alexis Guérin}
\address{Université Clermont Auvergne, LMBP, Campus des Cézeaux, 3 place Vasarely, TSA 60026, CS 60026, 63178 Aubière Cedex, France}
\email{\href{mailto:alexis.guerin@uca.fr}{alexis.guerin@uca.fr}}

\author{Felix Roz}
\address{Department of Mathematics Columbia University Room 509, MC 4406 2990 Broadway New York, NY 10027}
\email{\href{mailto:fr2588@columbia.edu}{fr2588@columbia.edu}}
\date{}

\begin{document}

\begin{abstract}
We construct an action of the positive part of the Witt algebra on Khovanov--Rozansky $\mathfrak{gl}_N$-link homology and show that link cobordisms induce equivariant maps between twists of the homology.
Moreover, the state spaces of simple webs are identified with standard representations of the Witt algebra on polynomials.
Some simple relations to Lee homology and genus bounds are derived from the analysis of this presentation.
\end{abstract}


\pagenumbering{arabic}

\maketitle

\setcounter{tocdepth}{1}
\tableofcontents

\section{Introduction}

\input{SF_Introduction}

\subsection{Acknowledgements}

\input{SF_Acknowledgements}

\subsection{Outline}
\input{SF_Outline}

\section{Algebraic Background}
\label{sec_AB}

\subsection{Conventions}
\input{SF_preliminairies}

\subsection{Symmetric Polynomials}

\input{SF_polynomials}

\subsection{The Witt algebra}

\input{SF_witt_algebra}

\subsection{Twisted actions}
\input{SF_twists_alternate}

\subsection{The relative homotopy category}

\input{SF_relative_homotopy_cat}

\section{The \texorpdfstring{$\Witt$}{Witt}-action on Foams}
\label{sec_WF}
\subsection{Webs and foams}

\input{SF_about_foams}


\subsection{The \texorpdfstring{$\Witt$}{Witt}-action}

\input{SF_witt_action}

\section{Twists of the \texorpdfstring{$\Witt$}{Witt}-action}
\label{sec_twist}

\subsection{Twists on webs}

\input{SF_twists_on_foams}

\subsection{Red dots}

\input{SF_red_dots_alternate}

\subsection{Useful morphisms}
\label{usual_morph_sec}

\input{SF_usual_morph}

\section{\(\Witt\)-equivariant Link Homology}
\label{sec_HI}

Going forward, single edges (\NB{\tikz{\draw (0,0) -- (0.5,0); }}), double edges (\NB{\tikz{\draw[double] (0,0) -- (0.5,0);}}), and bold edges
(\NB{\tikz[line width=.7mm]{\draw (0,0) -- (0.5,0); }}), will represent edges of thickness 1, 2, and 3 respectively in web diagrams.

\subsection{Definition of a link homology}

\input{SF_link_homology}

\subsection{The sliding lemmas}
\input{SF_Lemma}

\subsection{Reidemeister I}

\input{SF_RI}

\subsection{Reidemeister II}
\input{SF_RII}

\subsection{Reidemeister III}

\input{SF_RIII}

\subsection{Unframed Invariance}

\input{SF_unfrinv}

\section{Functoriality}
\label{sec_FU}

\input{SF_newfunctoriality}

\section{Polynomial Representations}
\label{sec_EX}

\input{SF_newcompute}

\section{\((2,m)\) Torus Links}
\label{sec_torus}
\subsection{Forktwist lemma}

\input{SF_forktwist}


\subsection{\((2,k)\) torus Links}

\input{SF_toruslinks}


\section{Future Directions}
\label{sec_future}

\input{SF_future}




\bibliographystyle{alphaurl}
\bibliography{biblio}

\end{document}

%% file: SF_Introduction.tex
Symmetries of link homology have been a consistently useful tool in the literature \cite{GOR13, GHM24, CG24, EQ23, GLW18, KhovanovSano_2025}.
In this paper, we study differential operators on  $\gl_N$-homology following the work of Qi, Robert, Sussan, and Wagner \cite{QRSW24}.
Specifically we consider an action by the Lie algebra of polynomial derivations which we call the positive half of the Witt algebra, \(\Witt\).
In other contexts, \(\Witt\) is referred to as a Cartan-type Lie algebra and denoted \(W_1\) \cite{Rudakov_1974}.
Witt actions were introduced by Khovanov and Rozansky on triply graded homology via the action of \(\Witt\) on Soergel bimodules \cite{KR16}.
A similar action was obtained in \cite{QRSW24} by lifting the action on symmetric polynomials to foams via the Robert-Wagner evaluation formula.
Using the foam action, it was shown that a copy of \(\sld\) embedded in \(\Witt\) acts naturally on \(\mathfrak{gl}_N\)-homology \cite{QRSW23}.
In this paper we show that the full Lie algebra \(\Witt\) acts naturally on \(\mathfrak{gl}_N\)-homology
and in a forthcoming paper we extend the action to exterior power colored homology.

\begin{thm}
\label{thm-isotopy}
Khovanov--Rozansky \(\mathfrak{gl}_N\)-homology is an isotopy invariant representation of \(\Witt\).
\end{thm}

\begin{rmk}
In this paper Khovanov--Rozansky homology refers to the fully equivariant theory over the ring of symmetric polynomials, denoted \(\scalars_N\).
To avoid conflicting terminology we will refer to this as \(\scalars_N\)-linear homology and reserve the term \textit{equivariant} for the \(\Witt\)-action.
\end{rmk}

The action of certain generators of \(\Witt\) have previously appeared in applications.
Based on the results in \cite{KR16}, the action of the generator $\wgen_1$ was used to define a $p$-DG structure on HOMFLY-PT homology and $\gl_{-2}$-homology, leading to categorifications of Jones and colored Jones polynomial at a root of unity \cite{QRSW21, qi_p-differential_2022}.
Similarly, the generator \(L_{-1}\) in \cite{QRSW24} agrees with the operator $\nabla$ introduced by Wang \cite{Wang24} extending the work of Shumakovitch \cite{Shumakovitch_2014}.

The Witt action should be viewed as a generalization of the grading structure on ordinary Khovanov homology.
The generator \(L_0\) acts by multiplying homogeneous elements by their quantum degree.
Therefore the \(L_0\)-module structure is the grading and \(L_0\)-equivariance is grading-preservation.
The familiar statement that the degree of a cobordism induced-map is its Euler characteristic has an analog for the \(\Witt\) action.
Grading shifts are replaced to local ``twists'' represented by red dots defined in \ref{sec_twist}.
\begin{thm}
\label{thm-functoriality}
Let \(K_1^{\lambda_1}\) and \(K_2^{\lambda_1}\) be knots decorated with red dots labelled \(\lambda_1\) and \(\lambda_2\) respectively.
If 
\(C: K_1^{\lambda_1} \to K_2^{\lambda_2}\) is a connected cobordism with
\(\lambda_2 - \lambda_1 = \chi(C)/2\), then \[\uKRW_N(C): \uKRW_N(K_1^{\lambda_1}) \to \uKRW_N(K_2^{\lambda_2})\] is \(\scalars_N\)-linear and \(\Witt\)-equivariant.
\end{thm}

In simple cases the \(\Witt\)-structure of the homology can be identified with standard representations.
\begin{prop}
Let \(U\) be the unknot and let \(U^\lambda\) be the unknot with a red dot labeled \(\lambda\).
As \(\scalars[x,y]^{S_2}\)-modules and \(\Witt\)-representations,
\begin{gather}
\uKRW_2(\varnothing) \cong \scalars[x,y]^{S_2}, \\
\uKRW_2(U) \cong (x-y)^{-1/2} \scalars[x,y],\\
\uKRW_2(U^\lambda) \cong (x-y)^{(\lambda-1)/2} \scalars[x,y].
\end{gather}
\end{prop}
In Section \ref{sec_EX} we determine the structure of these modules.
In particular, the only trivial submodule of \((x-y)^{(\lambda -1)/2}\scalars[x,y]\) is generated by \(1\) when it is present in the module.
Therefore, the inclusion is the unique up-to-scalars map \(\uKRW_2(\varnothing) \to \uKRW_2(U^\lambda)\) which is \(\scalars_2\)-linear and \(\Witt\)-equivariant.
Note that \((x-y)\) is the value associated to a handle attachment and \((x-y)^{\pm 1/2}\) is the value associated to a Reidemeister I move.

These observations allude to the results of Rasmussen and Tanaka \cite{Rasmussen_2005,Tanaka_2006}.
Further inspection of the homology groups above explains this and motivates the following conjecture:

\begin{conj}
The kernel of \(\wgen_{-1}\) is naturally isomorphic to universal Lee theory, that is the \((x^2-t)\)-deformation of Khovanov homology.
\end{conj}

%% file: SF_Acknowledgements.tex
We would like to thank Louis-Hadrien Robert, Joshua Sussan and Emmanuel Wagner for advising us on this project as well as Matt Hogancamp for helpful discussions related to twisted representations.
We thank the authors of \cite{QRSW23} and \cite{QRSW24} for allowing us to use the diagrams from their papers.

%% file: SF_Outline.tex
Sections \ref{sec_AB} and \ref{sec_WF} review the necessary algebraic and topological machinery.
In Section \ref{sec_AB} we define a framework for twisted representations, which are used to construct equivariant chain complexes,
and the relative homotopy category, which is the target of the homology functor.
These constructions are defined for a general Lie algebra \(\mathfrak{g}\), a Hopf algebra \(H\), and an associative algebra \(A\), but in the rest of the paper they will be specialized to the Witt algebra, its universal enveloping algebra, and the ring of symmetric polynomials, respectively.
Section \ref{sec_AB} also contains a brief discussion of the ring of symmetric polynomials and the Witt algebra.
Section \ref{sec_WF} reviews webs, foams, and the \(\Witt\)-action, based on \cite{QRSW24,QRSW23}. 

The remaining sections contain new results with the main goal of defining the $\Witt$-equivariant Khovanov--Rozansky $\gl_N$-homology. 
In Section~\ref{sec_twist}, we describe twists of the $\Witt$-action on foams.
Section~\ref{sec_HI} is devoted to defining the $\Witt$-module structure on Khovanov--Rozansky $\gl_N$-homology and to proving that it is isotopy invariant. 
Section~\ref{sec_FU} discusses functoriality, which follows directly from results in the previous sections and from~\cite{ETW17}.
Finally, in Section~\ref{sec_EX}, we analyze the \(\Witt\)-structure of several modules which arise in link homology specifically the state spaces of generalized theta webs.
In Section~\ref{sec_torus} we apply the results to
the homology of \((2,k)\)-torus links extending the work of the second author \cite{Roz23}.
Section~\ref{sec_future} concludes with directions for future research.

%% file: SF_preliminairies.tex
Let \(N\) be a positive integer fixed throughout the paper.
Let $\scalars$ be a unital ring in which 2 is invertible.
For all $x \in \scalars$, define $\Bar{x} := 1 - x$. 

Let \(M_i\) denote the degree $i$ subspace of a \(\Z\)-graded module \(M\).
Let $q^nM$ denote the shifted module where $(q^nM)_i = M_{i-n}$. 
Similarly if $C$ is a chain complex, let $t^nC$ denote the shifted complex with \((t^nC)_i = C_{i-n}.\)
Complexes will be cohomologically graded.
In other words, the $i$th differential goes from $C_{i-1}$ to $C_{i}$.

%% file: SF_polynomials.tex
Let \(\mathbf{a} = (a_1+\cdots+a_k) \) be an integer composition of \(n \in \N\),
let \(S_{\mathbf{a}} := S_{a_1} \times \cdots \times S_{a_k}\) be the associated Young subgroup of the symmetric group \(S_n\),
and let \( \scalars_{\mathbf{a}} := \scalars[x_1, \cdots, x_N]^{S_{\mathbf{a}}}\) be its ring of invariants.
In the special case where \(a_1=n\), \(\scalars_{\mathbf{a}}\) is the ring of symmetric polynomials.
When all \(a_i = 1\), \(\scalars_{\mathbf{a}}\) is the ring of all polynomials.
The rings \(\scalars_{\mathbf{a}}\) are graded by setting \(\deg(x_i)= 2\).

In the ring \(\scalars_{n}\), let \(e_i\), \(h_i\), and $p_i$ be the $i$th elementary, complete homogeneous, and power sum symmetric polynomials respectively.
The elements of the ring \(\scalars_N\) will be distinguished by capital letters.
For example, the variables, the elementary, complete homogeneous, and power sum symmetric polynomials will be denoted \(X_i\), \(E_i\), \(H_i\), and \(P_i\) respectively. 

For any \(n \in \Z, k \in \N\), the following quantum numbers are distinguished elements of the ring \(\Z[q,q^{-1}]\):
\[
[n] := \frac{q^n-q^{-n}}{q-q^{-1}}, \quad
[k]! := \prod_{j=1}^k[j], \quad
\croch{m}{a} := \prod_{i=1}^{a}\frac{[m+1-i]}{[i]}.
\]
Note that if $m$ is non-negative, $\croch{m}{a}=\frac{[m]!}{[a]![m-a]!}$.

%% file: SF_witt_algebra.tex
\begin{dfn}
    The \textit{Witt algebra} $\Wittt$ is the Lie algebra over $\scalars$ generated by $(\wgen_n)_{n \in \Z}$ with relations
    \begin{equation} \label{witt} [\wgen_n,\wgen_m] = (n-m) \wgen_{n+m}. \end{equation}
    for all \(n,m \in \Z\).
    Let \(\Wittt_{k}\) denote the sub Lie algebra generated by \((\wgen_n)_{n \geq k}\) for any integer \(k \geq -1\).
\end{dfn}

\begin{lem}
The map \(i: \sld \to \Wittt\) defined by
\begin{equation} 
    i : \left | \begin{array}{lll}
    e &\mapsto \wgen_{-1} \\
    h &\mapsto 2\wgen_{0} \\
    f &\mapsto -\wgen_1 
\end{array} \right. 
\end{equation}
is an injective map of Lie algebras.
\end{lem}

\begin{lem}
\label{std-witt}
The sub algebra \(\Wittt_{-1}\) is the Lie algebra of derivations on \(\scalars[x]\) under the standard representation 
\begin{equation}
    \wgen_n \mapsto - x^{n+1}\frac{\partial}{\partial x}.
\end{equation}
The operators
\begin{equation}
    \wgen_n:= -\sum_{i=1}^k x_i^{n+1}\frac{\partial}{\partial x_i}
\end{equation}
define the tensor product action of \(\Wittt_{-1}\) on \(\scalars[x_1,\ldots, x_k]\) which commutes with the permuation action of \(S_n\).
In particular this implies \(\Wittt_{-1}\) acts on the submodules \(\scalars_{\mathbf{a}}\).
\end{lem}

%% file: SF_twists_alternate.tex
\label{twist_gen}

Let \(\mathfrak{g}\) be a Lie algebra over \(\scalars\).

\begin{dfn}
A \textit{\(\mathfrak{g}\)-category},
is an \(\scalars\)-linear category \(\mathcal{C}\) such that
\(\mathrm{Hom}_\mathcal{C}(V,W)\)
is a \(\mathfrak{g}\)-representation and
\begin{equation}
    \label{enriched}
    X (f \circ g) = (X f) \circ g + f \circ (X g),
\end{equation}
for all \(X\in \mathfrak{g}\), \(f \in \mathrm{Hom}_\mathcal{C}(W,U)\) and \(g \in \mathrm{Hom}_\mathcal{C}(V,W)\).
\end{dfn}

\begin{dfn}
\label{defMC}
A \textit{Maurer--Cartan element} for an object \(V\) in a \(\mathfrak{g}\)-category \(\mathcal{C}\) is an \(\scalars\)-linear map \(\alpha_{(-)}: \mathfrak{g} \to \mathrm{End}_\mathcal{C}(V)\) such that, for all \(X,Y \in \mathfrak{g}\):
\[
\alpha_{[X,Y]} = [\alpha_X, \alpha_Y] + X\alpha_Y - Y\alpha_X,
\]
where \([-,-]\) in \( \mathrm{End}_\mathcal{C}(V) \) is the commutator.
\end{dfn}

\begin{lem}
    \label{prop_MC}
    Let \((V, \alpha)\) and \((W, \beta)\) be pairs of objects in \(\mathcal{C}\) and Maurer--Cartan elements.
    The following defines a new \(\mathfrak{g}\)-action on \( \mathrm{Hom}_\mathcal{C}(V,W) \):
    \begin{equation}
    \label{act_mc}
    X \cdot f := Xf + \beta_X \circ f - f \circ \alpha_X
    \end{equation}
    where \(X f\) denotes the original action of \(X\) on \(\mathrm{Hom}_\mathcal{C}(V, W)\).
\end{lem}



\begin{dfn}
The \(\mathfrak{g}\)-category \(\mathrm{Tw}(\mathcal{C})\)
is called the category of \textit{twisted \(\mathcal{C}\)-objects}.
Its objects are pairs of \(\mathcal{C}\)-objects and Maurer--Cartan elements.
As \(\scalars\)-modules 
\[\mathrm{Hom}_{\mathrm{Tw}(\mathcal{C})}((V,\alpha), (W, \beta)) := \mathrm{Hom}_\mathcal{C}(V,W),\]
and the \(\mathfrak{g}\)-action is defined by equation (\ref{act_mc}).
\end{dfn}

%% file: SF_relative_homotopy_cat.tex
\label{smash_sect}


Let $H$ be a Hopf algebra, let \(A\) an \(H\)-module algebra, and let \(A \# H\) be their smash product.
In the following, objects and morphisms in $A \# H$-mod will be called \textit{H-equivariant} $A$-modules and morphisms. 
There is the forgetful functor between the homotopy categories of modules over these associative algebras
\[\text{\text{For}}:\mathcal{C}(A\#H) \rightarrow \mathcal{C}(A).\]
An object in $\Ker(\text{For})$ is called \textit{relatively null-homotopic}
and is characterized by being null-homotopic after forgetting the \(H\) structure.

\begin{dfn}
    The \textit{relative homotopy category} is the Verdier quotient
    \[\RHC := \frac{\mathcal{C}(A\# H)}{\Ker(\text{For})}.\]
\end{dfn}




\begin{lem}[\cite{QRSW21}, Lemma 2.3]
    \label{lemme_RHC}
    Let $M \xrightarrow{f} N \xrightarrow{g} L$ be a short exact sequence in $\mathcal{C}(A\#H)$ that is termwise \(A\)-split exact. If $L$ is contractible, then $f$ induces an isomorphism in $\RHC$. Similarly, if $M$ is contractible then $g$ induces an isomorphism in $\RHC$.
\end{lem}

%% file: SF_about_foams.tex
In this section, we review (\(\mathfrak{gl}_N\)-)foams, following \cite{RW20a} \cite{QRSW23} and \cite{QRSW24}.

\begin{dfn}
	A \textit{web} is a finite, oriented, trivalent graph $\web = (V_\web,E_\web)$ embedded in \(\R^2\)
	endowed with a \textit{thickness function} $\ell : E(\web) \rightarrow \N$.
	The embeddeding is required to be smooth away from vertices
	and the thickness function must satisfy a flow condition: every vertex must be one of two types
	\begin{equation}
		\label{vertex1}
		\NB{\tikz[scale=0.6]{\input{sch/vertex_merge}}} \hspace{5mm} \text{or} \hspace{5mm} \NB{\tikz[scale=0.6]{\input{sch/vertex_split}}}.
	\end{equation}
	The first is called a \textit{merge}, the second a \textit{split}.
	In both cases, the edge with the largest label is called the \textit{thick} edge and the other two are \textit{thin} edges.
	Disjoint oriented circles without vertices are allowed.
\end{dfn}

\begin{figure}[ht]
  \centering
  \NB{\tikz[]{\input{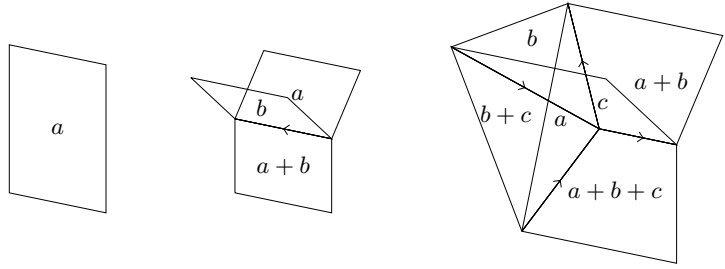}}}
  \caption{The three local models of the definition of foams. One denotes by $Y^{(a,b)}$ the one in the middle and $T^{(a,b,c)}$ the one on the right.}
  \label{deffoam}
\end{figure}

\begin{dfn}
	A \textit{foam} is a compact, finite, 2-dimensional CW complex \(F\) embedded in \(\R^2 \times [0,1]\) together with a \textit{thickness function} $\ell : F^2 \rightarrow \N$ 
	and a \textit{decoration}, $P_f \in \scalars_{\ell(f),N-\ell(f)}$, for every \(f \in F^2\).
	The embedding is required to be smooth on the interiors of the 1 and 2-cells and
	every point $p$ of $\foam$ is required to have a closed neighbourhood homeomorphic to one of the following, shown in Figure \ref{deffoam}: 
	\begin{enumerate}
		\item a disk, when $p$ belongs to a unique facet, 
		\item $Y \times[0,1]$, where $Y$ is the neighbourhood of a merge or split vertex of a web, when $p$ belongs to three facets, 
		\item the cone over the 1-skeleton of a tetrahedron with $p$ as the vertex of the cone, so that it belongs to six facets. 
	\end{enumerate}
	The 2-cells are called \textit{facets}, the 1-cells are called \textit{bindings}, and the 0-cells are called \textit{singular vertices}.
	At a binding, the thicknesses of the three facets must satisfy a flow condition analogous to the one for webs.
	Orientations for bindings are induced by those of the thin facets and by the opposite of the thick facet.
	The \textit{boundary} $\partial \foam$ of $\foam$ is the closure of the set of boundary points of facets that do not belong to a binding.
\end{dfn} 

Let \(p_i\) denote the power sum in the first \(a\) variables in \(\scalars_{a,N-a}\) and \(\widehat{p}_i\) denote the power sum in the last \(N-a\) variables in \(\scalars_{a,N-a}\). 
The following foams are abbreviations for the these common decorations on a facet of thickness \(a\):
\[ 
\NB{\tikz[scale=.8]{\input{sch/notation_deco2}}}:= \NB{\tikz[scale=.8]{\input{sch/notation_deco1}}} \hspace{5mm} \text{and} \hspace{5mm} \NB{\tikz[scale=.8]{\input{sch/notation_deco3}}}:= \NB{\tikz[scale=.8]{\input{sch/notation_antideco}}}.
\]
In particular, \(\deco_0=a\) and \(\antideco_0 = N-a\).

Partition $\foam^1$ as $\foam^1_\circ \sqcup \foam^1_-$, where $\foam^1_\circ$ is the collection of bindings diffeomorphic to circles, and $\foam^1_-$ the ones diffeomorphic to intervals.
\begin{dfn}
Let $f$ be a facet of thickness $a$,
let $b \in \foam^1_-$ have a neighbourhood homeomorphic to $Y^{(a,b)}$,
and let $v \in \foam^0$ have a neighbourhood homeomorphic to $T^{(a,b,c)}$.
The \textit{degrees} of facets, bindings, singular vertices and foams are defined by
\begin{align}
\deg_N(f) &:=   a(N-a) \chi(f),\\
\deg_N(b) &:=   ab + (a+b)(N-a-b),\\
\deg_N(v) &:=   ab+bc+ac+(a+b+c)(N-a-b-c),\\
\deg_N(\foam) &:= \sum_{f \in \foam^2} \left(\deg(P_f) - \deg_N(f)\right) + \sum_{b \in \foam^1_-} \deg_N(s) - \sum_{v \in F^0} \deg_N(v)
\end{align}
\end{dfn}

\begin{figure}[ht]
   \centering
   \begin{tikzpicture}[xscale=2.4, yscale =-3.2]
     \node (pol) at (0,0) {\NB{\tikz[font=\tiny,scale=1.5]{\input{sch/foam_polynomial}}}};
     \node[yshift= -1.4cm] at (pol) {polynomial};
     \node[yshift = -1.8cm] at (pol) {$\deg(P_1)+\deg(P_2)$};
     \node (asso) at(1.85,0) {\NB{\tikz[font=\tiny,scale=0.7]{\input{sch/foam_associativity}}}};
     \node[yshift = -1.4cm] at (asso) {associativity};
     \node[yshift= -1.8cm] at (asso) {$0$};
     \node (coasso) at (3.7,0) {\NB{\tikz[font=\tiny, xscale=-0.7, yscale=0.7]{\input{sch/foam_coassociativity}}}};
     \node[yshift = -1.4cm] at (coasso) {associativity};
     \node[yshift= -1.8cm] at (coasso) {$0$};
     \node (digcup) at (0,1) {\NB{\tikz[font=\tiny]{\input{sch/foam_digoncup_ab_arrows}}}};
     \node[yshift = -1.64cm] at (digcup) {digon-cup};
     \node[yshift= -2.04cm] at (digcup) {$-ab$};
     \node (digcap) at (2,1.1) {\NB{\tikz[font=\tiny]{\input{sch/foam_digoncap_ab_arrows}}}};
     \node[yshift = -1.4cm] at (digcap) {digon-cap};
     \node[yshift= -1.8cm] at (digcap) {$-ab$};
     \node (zip) at (4,1.1) {\NB{\tikz[font=\tiny]{\input{sch/foam_unzip_ab_arrows}}}};
     \node[yshift = -1.4cm] at (zip) {unzip};
     \node[yshift= -1.8cm] at (zip) {$ab$};
     \node (unzip) at (0,2.1) {\NB{\tikz[font=\tiny]{\input{sch/foam_zip_ab_arrows}}}};
     \node[yshift = -1.4cm] at (unzip) {zip};
     \node[yshift= -1.8cm] at (unzip) {$ab$};
     \node (cup) at (1.4,2.1) {\NB{\tikz[font=\tiny]{\input{sch/foam_cup_a_empty}}}};
     \node[yshift = -1.4cm] at (cup) {cup};
     \node[yshift= -1.8cm] at (cup) {$-a(N-a)$};
     \node (cap) at (2.6,2.1) {\NB{\tikz[font=\tiny]{\input{sch/foam_cap_a_empty}}}};
     \node[yshift = -1.4cm] at (cap) {cap};
     \node[yshift= -1.8cm] at (cap) {$-a(N-a)$};
     \node (saddle) at (4,2.1) {\NB{\tikz[font=\tiny]{\input{sch/foam_saddle_a_empty}}}};
     \node[yshift = -1.4cm] at (saddle) {saddle};
     \node[yshift= -1.8cm] at (saddle) {$a(N-a)$};
   \end{tikzpicture}
   \caption{The degree of a basic foam is given below the name of each
     of the local models.}\label{localmod}
 \end{figure}
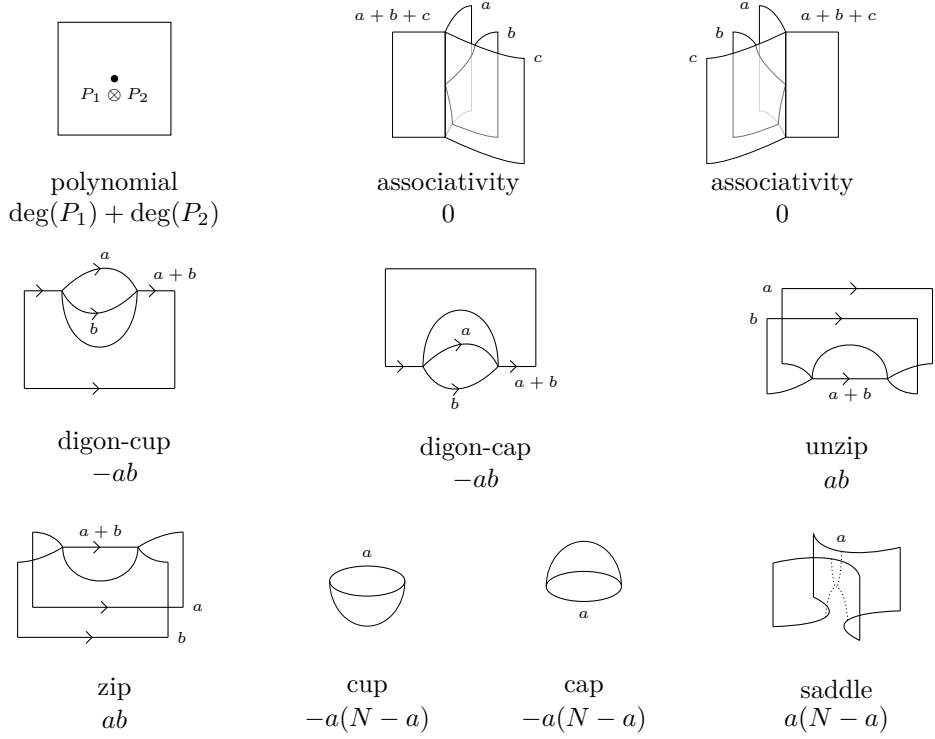

\begin{dfn}
    A foam is called \textit{basic} if it is equal to \( \Gamma \times [0,1] \) outside a cylinder \( B \times [0,1] \), and homeomorphic to one of the local models shown in Figure \ref{localmod} within the cylinder. A foam is said to be in \textit{good position} if it is a composition of basic foams and traces of isotopies.
\end{dfn}

\begin{lem}[\cite{QW24}]
Every foam is isotopic to a foam in good position.
\end{lem}

A generic section $\foam_t := \foam \cap (\R^2 \times \{t\})$ is a web.
In particular, the boundary of a foam $\foam$ is a closed web in $\R^2$.
Thus, one can define the following category.

\begin{dfn}
The category $\overline{\catf}$ has webs as objects and morphism sets
\[
\overline{\catf}(\web_0, \web_1) :=  \sfrac{\{\text{foams with } \foam_0 = -\web_0 \text{ and } \foam_1 = \web_1 \}}{\text{ambient isotopy}} ,
\]
where $-\web_0$ denotes reversed orientations.
Composition is given by stacking foams.
The identity for $\web$ is the cylinder $\web \times [0,1]$.
Decorations are multiplicative and degrees are additive under composition.
Let \(\catf\) be the graded, additive, \(\scalars_N\)-linear completion of \(\overline{\catf}\).
\end{dfn}



\begin{dfn}
Let \(\langle - \rangle_N\) denote the \(\scalars_N\)-valued Robert-Wagner \textit{foam evaluation}
of closed foams \cite{RW20a}.
For any web \(\web\), define an \(\scalars_N\)-bilinear form \(\langle \cdot ; \cdot \rangle_N\) on \(\catf(\varnothing, \web)\) by
\begin{equation}
\langle F ; G \rangle_N := \langle \overline{G} \circ F \rangle_N,
\end{equation}
where \(\overline{G}: \web \to \varnothing\) is the reflection of a foam \(G: \varnothing \to \web\) along the plane \(\R^2 \times \{1/2\}\).
The \(\scalars_N\)-module
\begin{equation}
\ssp_N(\web) := \catf(\varnothing, \web) / \mathrm{Rad}\langle \cdot ; \cdot \rangle_N
\end{equation}
is called the \textit{state space} of \(\web\).
By the universal construction of \cite{BHMV92} the state space extends to a functor \(\ssp_N: \catf \to \scalars_N\text{-mod}\).
\end{dfn}

%% file: sch/vertex_merge.tex
\begin{scope}[font=\tiny]
    \draw [->] (0,0) -- (90:1) node[pos =1, above] {$a+b$};
    \draw [-<] (0,0) -- (-30:1) node[pos =1, below] {$b$};
    \draw [-<] (0,0) -- (-150:1) node[pos =1, below] {$a$};
\end{scope}

%% file: sch/vertex_split.tex
\begin{scope}[font=\tiny]
  \begin{scope}
    \draw [-<] (0,0) -- (-90:1) node[pos =1, below] {$a+b$};
    \draw [->] (0,0) -- (30:1) node[pos =1, above] {$b$};
    \draw [->] (0,0) -- (150:1) node[pos =1, above] {$a$};
  \end{scope}
\end{scope}

%% file: sch/notation_deco2.tex
\begin{scope}
  \draw (-1,-1) rectangle (1,1);
  \coordinate (.5,.5) node[]{$\deco_i$};
\end{scope}

%% file: sch/notation_deco1.tex
\begin{scope}
  \draw (-1,-1) rectangle (1,1) coordinate [midway] (A);
  \fill (A) circle (0.5mm) node[below] {$p_i$};
\end{scope}

%% file: sch/notation_deco3.tex
\begin{scope}
  \draw (-1,-1) rectangle (1,1);
  \coordinate (.5,.5) node[]{$\antideco_i$};
\end{scope}

%% file: sch/notation_antideco.tex
\begin{scope}
    \draw (-1,-1) rectangle (1,1) coordinate [midway] (A);
  \fill (A) circle (0.5mm) node[below] {$\widehat{p}_i$};
\end{scope}

%% file: SF_witt_action.tex
Fix three elements \(s, \lambda, \mu \in \scalars\) and let
$\lambda_n := \lambda(n+1)$ and $\mu_n := \mu(n+1)$\footnote{In \cite{QRSW24}, Witt sequences were used instead. When $2$ is invertible in the ground ring, these sequences are of the form $(\lambda(n+1))_{n\in \N_{-1}}$ with $\lambda \in \scalars$.}.
Using these distinguished elements define the following operators $\Wgen_n$ on basic foams:

\begin{align}
    \label{polyn}
    \Wgen_n \left (~\NB{\tikz[scale=1.5, font=\tiny]{\input{sch/foam_polynomial}}}~\right) =  \NB{\tikz[scale=1.7, font=\tiny]{\input{sch/foam_polynomial_L}}} \hspace{1mm};
\end{align}

\begin{align}
    \Wgen_n \left (  \NB{\tikz[xscale=.9,yscale=.7,font=\tiny]{\input{sch/foam_associativity}}} \right ) = \Wgen_n \left (  \NB{\tikz[xscale=-.9, yscale=.7,font=\tiny]{\input{sch/foam_coassociativity}}} \right ) = 0 \hspace{1mm};
\end{align}

\begin{align}
    \Wgen_n \left ( ~~ \NB{\tikz[scale=1.2,font=\tiny]{\input{sch/foam_digoncup_ab_arrows}}}  \right ) = ~~ &\lambda_n ~~ \NB{\tikz[scale=1.2,font=\tiny]{\input{sch/foam_digoncup_ab_0n}}}+ \mu_n ~~ \NB{\tikz[scale=1.2,font=\tiny]{\input{sch/foam_digoncup_ab_n0}}} \\ \nonumber &+ s\sum_{k+\ell=n} \NB{\tikz[scale=1.2,font=\tiny]{\input{sch/foam_digoncup_ab_kl}}}\hspace{1mm};
\end{align}

\begin{align}
    \Wgen_n \left (~~ \NB{\tikz[scale=1.2,font=\tiny]{\input{sch/foam_digoncap_ab_arrows}}} \right ) = ~~ & -\lambda_n ~~ \NB{\tikz[scale=1.2,font=\tiny]{\input{sch/foam_digoncap_ab_0n}}} - \mu_n ~~ \NB{\tikz[scale=1.2,font=\tiny]{\input{sch/foam_digoncap_ab_n0}}} \\ \nonumber &+ \bar{s}\sum_{k+\ell=n} \NB{\tikz[scale=1.2,font=\tiny]{\input{sch/foam_digoncap_ab_kl}}}\hspace{1mm};
\end{align}

\begin{align}
    \Wgen_n \left ( ~~ \NB{\tikz[scale=1.2,font=\tiny]{\input{sch/foam_zip_ab_arrows}}}  \right ) = ~~ &\lambda_n ~~ \NB{\tikz[scale=1.2,font=\tiny]{\input{sch/foam_zip_ab_0n}}}+ \mu_n ~~ \NB{\tikz[scale=1.2,font=\tiny]{\input{sch/foam_zip_ab_n0}}} \\ \nonumber &- \bar{s}\sum_{k+\ell=n} \NB{\tikz[scale=1.2,font=\tiny]{\input{sch/foam_zip_ab_kl}}}\hspace{1mm};
\end{align}

\begin{align}
    \Wgen_n \left ( \NB{\tikz[scale=1.2,font=\tiny]{\input{sch/foam_unzip_ab_arrows}}} ~~ \right ) = ~~ & -\lambda_n  \NB{\tikz[scale=1.2,font=\tiny]{\input{sch/foam_unzip_ab_0n}}} - \mu_n  \NB{\tikz[scale=1.2,font=\tiny]{\input{sch/foam_unzip_ab_n0}}} \\ \nonumber &-s\sum_{k+\ell=n} \NB{\tikz[scale=1.2,font=\tiny]{\input{sch/foam_unzip_ab_kl}}}\hspace{1mm};
\end{align}

\begin{align}
    \Wgen_n \left ( \NB{\tikz[scale=1.5,font=\tiny]{\input{sch/foam_cup_a_empty}}} \right ) = \frac{1}{2} \sum_{k+\ell=n} \NB{\tikz[scale=1.5,font=\tiny]{\input{sch/foam_cup_a_kl}}} \hspace{1mm}; \hspace{3mm} 
    \Wgen_n \left ( \NB{\tikz[scale=1.5,font=\tiny]{\input{sch/foam_cap_a_empty}}} \right ) = \frac{1}{2} \sum_{k+\ell=n} \NB{\tikz[scale=1.5,font=\tiny]{\input{sch/foam_cap_a_kl}}}\hspace{1mm};
\end{align}

\begin{align}
    \label{saddle}
    \Wgen_n \left ( \NB{\tikz[scale=1.3,font=\tiny]{\input{sch/foam_saddle_a_empty}}} \right ) = - \frac{1}{2} \sum_{k+\ell=n}~~ \NB{\tikz[scale=1.3,font=\tiny]{\input{sch/foam_saddle_a_kl}}}\hspace{1mm};
\end{align}

These operators are extended to arbitrary foams in good position by the Leibniz rule with respect to gluing and acting by 0 on traces of isotopies.

\begin{thm}[{\cite{QRSW24}}]
\label{witt-thm}
    The operators $\Wgen_n$ 
    define an action of $\Witt$ on the $\scalars_N$-module generated by isotopy classes of foams in good position
    such that the Robert--Wagner evaluation on closed foams is an equivariant map.
\end{thm}

%% file: sch/foam_polynomial.tex
\begin{scope}
  \draw (0,0) rectangle (1,1) coordinate [midway] (A);
  \fill (A) circle (0.33mm) node[below] {$P_1 \otimes P_2$};
\end{scope}

%% file: sch/foam_polynomial_L.tex
\begin{scope}
  \draw (0,0) rectangle (1,1) coordinate [midway] (A);
  \fill (A) circle (0.29mm) node[below] {$\wgen_n\cdot (P_1 \otimes P_2)$};
\end{scope}

%% file: sch/foam_associativity.tex
\begin{scope}
  \coordinate (ABCt) at (0,2);
  \coordinate (Ot) at (1,2);
  \coordinate (At) at (1.5,2.5);
  \coordinate (Bt) at (2,2);
  \coordinate (Ct) at (2.5,1.5);
  \begin{scope}[yshift= -2cm]
    \coordinate (ABCb) at (0,2);
    \coordinate (Ob) at (1,2);
    \coordinate (Ab) at (1.5,2.5);
    \coordinate (Bb) at (2,2);
    \coordinate (Cb) at (2.5,1.5);
  \end{scope}
  \begin{scope}[fill = white, fill opacity = 0.45, draw = black,
    draw opacity =1, very thin]
    \filldraw (Ob) .. controls +(0,0) and +(-0.3, 0) .. (Ab)
    coordinate[pos=0.5] (ABb) -- (At)
    .. controls +(-0.3, 0) and +(0,0) .. (Ot) -- (Ob);
    \fill (Ob) .. controls +(0,0) and +(-0.5, 0) .. (Cb) -- (Ct)
    .. controls +(-0.5, 0) and +(0,0) .. (Ot) coordinate[pos=0.5]
    (BCt) -- (Ob) coordinate[pos=0.5] (Om);
    \filldraw (ABb) .. controls +(0,0) and +(-0.3, 0) .. (Bb)-- (Bt)
    .. controls +(-0.3, 0) and +(0,0) .. (BCt) .. controls +(0, -0.3) and
    +(0,0) .. (Om) .. controls +(0,0) and +(0, 0.3)  .. (ABb);
    \filldraw (Ob) .. controls +(0,0) and +(-0.5, 0) .. (Cb) -- (Ct)
    .. controls +(-0.5, 0) and +(0,0) .. (Ot) -- (Ob);
    \filldraw (ABCt) -- (Ot) -- (Ob) -- (ABCb) -- (ABCt);
  \end{scope}
  
  \begin{scope}
    \node[above] at (ABCt) {$a+b+c$};
    \node[right] at (At) {$a$};
    \node[right] at (Bt) {$b$};
    \node[right] at (Ct) {$c$};
  \end{scope}
\end{scope}

%% file: sch/foam_coassociativity.tex
\begin{scope}
  \coordinate (ABCt) at (0,2);
  \coordinate (Ot) at (1,2);
  \coordinate (At) at (1.5,2.5);
  \coordinate (Bt) at (2,2);
  \coordinate (Ct) at (2.5,1.5);
  \begin{scope}[yshift= -2cm]
    \coordinate (ABCb) at (0,2);
    \coordinate (Ob) at (1,2);
    \coordinate (Ab) at (1.5,2.5);
    \coordinate (Bb) at (2,2);
    \coordinate (Cb) at (2.5,1.5);
  \end{scope}
  \begin{scope}[fill = white, fill opacity = 0.45, draw = black,
    draw opacity =1, very thin]
    \filldraw (Ob) .. controls +(0,0) and +(-0.3, 0) .. (Ab)
    coordinate[pos=0.5] (ABb) -- (At)
    .. controls +(-0.3, 0) and +(0,0) .. (Ot) -- (Ob);
    \fill (Ob) .. controls +(0,0) and +(-0.5, 0) .. (Cb) -- (Ct)
    .. controls +(-0.5, 0) and +(0,0) .. (Ot) coordinate[pos=0.5]
    (BCt) -- (Ob) coordinate[pos=0.5] (Om);
    \filldraw (ABb) .. controls +(0,0) and +(-0.3, 0) .. (Bb)-- (Bt)
    .. controls +(-0.3, 0) and +(0,0) .. (BCt) .. controls +(0, -0.3) and
    +(0,0) .. (Om) .. controls +(0,0) and +(0, 0.3)  .. (ABb);
    \filldraw (Ob) .. controls +(0,0) and +(-0.5, 0) .. (Cb) -- (Ct)
    .. controls +(-0.5, 0) and +(0,0) .. (Ot) -- (Ob);
    \filldraw (ABCt) -- (Ot) -- (Ob) -- (ABCb) -- (ABCt);
  \end{scope}
  
  \begin{scope}
    \node[above] at (ABCt) {$a+b+c$};
    \node[left] at (At) {$a$};
    \node[left] at (Bt) {$b$};
    \node[left] at (Ct) {$c$};
  \end{scope}
\end{scope}

%% file: sch/foam_digoncup_ab_arrows.tex
\begin{scope}
  \begin{scope}
    \coordinate (L) at (0,0);
    \coordinate (R) at (2,0);
    \coordinate (ML) at (0.5, 0);
    \coordinate (MR) at (1.5, 0);
    \draw[->-] (L) -- (ML);
    \draw[->-] (MR) -- (R) node[above] {$a+b$};
    \draw[->-] (ML).. controls + (0.4, 0.4) and +(-0.2, 0.4) .. (MR)
    node[above, midway] {$a$};
    \draw[->-] (ML).. controls + (0.2, -0.4) and +(-0.4, -0.4) .. (MR)
    node[below, midway] {$b$};
  \end{scope}  
 \begin{scope}[yshift = -1.3cm]
    \coordinate (LB) at (0,0);
    \coordinate (RB) at (2,0);
    \draw[->-] (LB) -- (RB);
  \end{scope}  
  \draw (R) -- (RB);
  \draw (L) -- (LB);
  \draw (ML) .. controls +(0, -1) and +(0, -1) .. (MR);
\end{scope}

%% file: sch/foam_digoncup_ab_0n.tex
\begin{scope}
  \begin{scope}
    \coordinate (L) at (0,0);
    \coordinate (R) at (2,0);
    \coordinate (ML) at (0.5, 0);
    \coordinate (MR) at (1.5, 0);
    \draw[->-] (L) -- (ML);
    \draw[->-] (MR) -- (R) node[above] {$a+b$};
    \draw[->-] (ML).. controls + (0.4, 0.4) and +(-0.2, 0.4) .. (MR)
    node[above, midway] {$a$} node[pos=0.4,below] {$\deco_0$};
    \draw[->-] (ML).. controls + (0.2, -0.4) and +(-0.4, -0.4) .. (MR)
    node[below, pos=0.35] {$b$} node[pos=0.75,below] {$\deco_n$};
  \end{scope}  
 \begin{scope}[yshift = -1.3cm]
    \coordinate (LB) at (0,0);
    \coordinate (RB) at (2,0);
    \draw[->-] (LB) -- (RB);
  \end{scope}  
  \draw (R) -- (RB);
  \draw (L) -- (LB);
  \draw(ML) .. controls +(0, -1) and +(0, -1) .. (MR);
\end{scope}

%% file: sch/foam_digoncup_ab_n0.tex
\begin{scope}
  \begin{scope}
    \coordinate (L) at (0,0);
    \coordinate (R) at (2,0);
    \coordinate (ML) at (0.5, 0);
    \coordinate (MR) at (1.5, 0);
    \draw[->-] (L) -- (ML);
    \draw[->-] (MR) -- (R) node[above] {$a+b$};
    \draw[->-] (ML).. controls + (0.4, 0.4) and +(-0.2, 0.4) .. (MR)
    node[above, midway] {$a$} node[pos=0.4,below] {$\deco_n$};
    \draw[->-] (ML).. controls + (0.2, -0.4) and +(-0.4, -0.4) .. (MR)
    node[below, pos=0.35] {$b$} node[pos=0.75,below] {$\deco_0$};
  \end{scope}  
 \begin{scope}[yshift = -1.3cm]
    \coordinate (LB) at (0,0);
    \coordinate (RB) at (2,0);
    \draw[->-] (LB) -- (RB);
  \end{scope}  
  \draw (R) -- (RB);
  \draw (L) -- (LB);
  \draw (ML) .. controls +(0, -1) and +(0, -1) .. (MR);
\end{scope}

%% file: sch/foam_digoncup_ab_kl.tex
\begin{scope}
  \begin{scope}
    \coordinate (L) at (0,0);
    \coordinate (R) at (2,0);
    \coordinate (ML) at (0.5, 0);
    \coordinate (MR) at (1.5, 0);
    \draw[->-] (L) -- (ML);
    \draw[->-] (MR) -- (R) node[right] {$a+b$};
    \draw[->-] (ML).. controls + (0.4, 0.4) and +(-0.2, 0.4) .. (MR)
    node[above, midway] {$a$} node[pos=0.4,below] {$\deco_k$};
    \draw[->-] (ML).. controls + (0.2, -0.4) and +(-0.4, -0.4) .. (MR)
    node[below, pos=0.35] {$b$} node[pos=0.75,below] {$\deco_l$};
  \end{scope}  
 \begin{scope}[yshift = -1.3cm]
    \coordinate (LB) at (0,0);
    \coordinate (RB) at (2,0);
    \draw[->-] (LB) -- (RB);
  \end{scope}  
  \draw (R) -- (RB);
  \draw (L) -- (LB);
  \draw (ML) .. controls +(0, -1) and +(0, -1) .. (MR);
\end{scope}

%% file: sch/foam_digoncap_ab_arrows.tex
\begin{scope}
  \begin{scope}
    \coordinate (L) at (0,0);
    \coordinate (R) at (2,0);
    \coordinate (ML) at (0.5, 0);
    \coordinate (MR) at (1.5, 0);
    \draw[->-] (L) -- (ML);
    \draw[->-] (MR) -- (R) node[below] {$a+b$};
    \draw[->-] (ML).. controls + (0.4, 0.4) and +(-0.2, 0.4) .. (MR)
    node[above, midway] {$a$};
    \draw[->-] (ML).. controls + (0.2, -0.4) and +(-0.4, -0.4) .. (MR)
    node[below, midway] {$b$};
  \end{scope}  
 \begin{scope}[yshift = 1.3cm]
    \coordinate (LB) at (0,0);
    \coordinate (RB) at (2,0);
    \draw (LB) -- (RB);
  \end{scope}  
  \draw (R) -- (RB);
  \draw (L) -- (LB);
  \draw (ML) .. controls +(0, 1) and +(0, 1) .. (MR);
\end{scope}

%% file: sch/foam_digoncap_ab_0n.tex
 \begin{scope}
  \begin{scope}
    \coordinate (L) at (0,0);
    \coordinate (R) at (2,0);
    \coordinate (ML) at (0.5, 0);
    \coordinate (MR) at (1.5, 0);
    \draw[->-] (L) -- (ML);
    \draw[->-] (MR) -- (R) node[below] {$a+b$};
    \draw[->-] (ML).. controls + (0.4, 0.4) and +(-0.2, 0.4) .. (MR)
    node[above, pos=0.2] {$a$} node[pos=0.6,above] {$\deco_0$};
    \draw[->-] (ML).. controls + (0.2, -0.4) and +(-0.4, -0.4) .. (MR)
    node[below, pos=0.5] {$b$} node[pos=0.6, above] {$\deco_n$};
  \end{scope}  
 \begin{scope}[yshift = 1.3cm]
    \coordinate (LB) at (0,0);
    \coordinate (RB) at (2,0);
    \draw (LB) -- (RB);
  \end{scope}  
  \draw (R) -- (RB);
  \draw (L) -- (LB);
  \draw (ML) .. controls +(0, 1) and +(0, 1) .. (MR);
\end{scope}

%% file: sch/foam_digoncap_ab_n0.tex
 \begin{scope}
  \begin{scope}
    \coordinate (L) at (0,0);
    \coordinate (R) at (2,0);
    \coordinate (ML) at (0.5, 0);
    \coordinate (MR) at (1.5, 0);
    \draw[->-] (L) -- (ML);
    \draw[->-] (MR) -- (R) node[below] {$a+b$};
    \draw[->-] (ML).. controls + (0.4, 0.4) and +(-0.2, 0.4) .. (MR)
    node[above, pos=0.2] {$a$} node[pos=0.6,above] {$\deco_n$};
    \draw[->-] (ML).. controls + (0.2, -0.4) and +(-0.4, -0.4) .. (MR)
    node[below, pos=0.5] {$b$} node[pos=0.6, above] {$\deco_0$};
  \end{scope}  
 \begin{scope}[yshift = 1.3cm]
    \coordinate (LB) at (0,0);
    \coordinate (RB) at (2,0);
    \draw (LB) -- (RB);
  \end{scope}  
  \draw (R) -- (RB);
  \draw (L) -- (LB);
  \draw (ML) .. controls +(0, 1) and +(0, 1) .. (MR);
\end{scope}

%% file: sch/foam_digoncap_ab_kl.tex
 \begin{scope}
  \begin{scope}
    \coordinate (L) at (0,0);
    \coordinate (R) at (2,0);
    \coordinate (ML) at (0.5, 0);
    \coordinate (MR) at (1.5, 0);
    \draw[->-] (L) -- (ML);
    \draw[->-] (MR) -- (R) node[below] {$a+b$};
    \draw[->-] (ML).. controls + (0.4, 0.4) and +(-0.2, 0.4) .. (MR)
    node[above, pos=0.2] {$a$} node[pos=0.6,above] {$\deco_k$};
    \draw[->-] (ML).. controls + (0.2, -0.4) and +(-0.4, -0.4) .. (MR)
    node[below, pos=0.5] {$b$} node[pos=0.6, above] {$\deco_l$};
  \end{scope}  
 \begin{scope}[yshift = 1.3cm]
    \coordinate (LB) at (0,0);
    \coordinate (RB) at (2,0);
    \draw (LB) -- (RB);
  \end{scope}  
  \draw (R) -- (RB);
  \draw (L) -- (LB);
  \draw (ML) .. controls +(0, 1) and +(0, 1) .. (MR);
\end{scope}

%% file: sch/foam_zip_ab_arrows.tex
\begin{scope}
  \begin{scope}
    \coordinate (L1) at (0.2,0.4);
    \coordinate (L2) at (0,0);
    \coordinate (R1) at (2.2,0.4);
    \coordinate (R2) at (2,0);
    \coordinate (ML) at (0.6, 0.2);
    \coordinate (MR) at (1.6, 0.2);
    \draw[->-] (ML) -- (MR) node[above, midway] {$a+b$};
    \draw (MR) .. controls +(0, 0) and +(-0.3,0) .. (R1) ;
    \draw (MR) .. controls +(0, 0) and +(-0.3,0) .. (R2);
    \draw (L1) .. controls +( 0.3, 0) and +(0,0) .. (ML);
    \draw (L2) .. controls +( 0.3, 0) and +(0,0) .. (ML);
  \end{scope}  
 \begin{scope}[yshift = -1cm]
    \coordinate (L1B) at (0.2,0.4);
    \coordinate (L2B) at (0,0);
    \coordinate (R1B) at (2.2,0.4);
    \coordinate (R2B) at (2,0);
    \draw[->-] (L1B) .. controls +( 0, 0) and +(0,0) .. (R1B) node [right, pos
    = 1] {$a$};
    \draw[->-] (L2B) .. controls +( 0, 0) and +(0,0) .. (R2B) node [right, pos
    = 1] {$b$};
 \end{scope}  
  \draw (R1) -- (R1B);
  \draw (R2) -- (R2B);
  \draw (L1) -- (L1B);
  \draw (L2) -- (L2B);
  \draw (ML) .. controls +(0, -0.6) and +(0, -0.6) .. (MR);
\end{scope}

%% file: sch/foam_zip_ab_0n.tex
\begin{scope}
  \begin{scope}
    \coordinate (L1) at (0.2,0.4);
    \coordinate (L2) at (0,0);
    \coordinate (R1) at (2.2,0.4);
    \coordinate (R2) at (2,0);
    \coordinate (ML) at (0.6, 0.2);
    \coordinate (MR) at (1.6, 0.2);
    \draw[->-] (ML) -- (MR) node[above, midway] {$a+b$};
    \draw (MR) .. controls +(0, 0) and +(-0.3,0) .. (R1) ;
    \draw (MR) .. controls +(0, 0) and +(-0.3,0) .. (R2);
    \draw (L1) .. controls +( 0.3, 0) and +(0,0) .. (ML);
    \draw (L2) .. controls +( 0.3, 0) and +(0,0) .. (ML);
  \end{scope}  
 \begin{scope}[yshift = -1cm]
    \coordinate (L1B) at (0.2,0.4);
    \coordinate (L2B) at (0,0);
    \coordinate (R1B) at (2.2,0.4);
    \coordinate (R2B) at (2,0);
    \draw[->-] (L1B) .. controls +( 0, 0) and +(0,0) .. (R1B) node [right, pos = 1] {$a$} node[above,pos=0.6] {$\deco_0$};
    \draw[->-] (L2B) .. controls +( 0, 0) and +(0,0) .. (R2B) node [right, pos = 1] {$b$} node[above,pos=0.4] {$\deco_n$};
 \end{scope}  
  \draw (R1) -- (R1B);
  \draw (R2) -- (R2B);
  \draw (L1) -- (L1B);
  \draw (L2) -- (L2B);
  \draw (ML) .. controls +(0, -0.6) and +(0, -0.6) .. (MR);
\end{scope}

%% file: sch/foam_zip_ab_n0.tex
\begin{scope}
  \begin{scope}
    \coordinate (L1) at (0.2,0.4);
    \coordinate (L2) at (0,0);
    \coordinate (R1) at (2.2,0.4);
    \coordinate (R2) at (2,0);
    \coordinate (ML) at (0.6, 0.2);
    \coordinate (MR) at (1.6, 0.2);
    \draw[->-] (ML) -- (MR) node[above, midway] {$a+b$};
    \draw (MR) .. controls +(0, 0) and +(-0.3,0) .. (R1) ;
    \draw (MR) .. controls +(0, 0) and +(-0.3,0) .. (R2);
    \draw (L1) .. controls +( 0.3, 0) and +(0,0) .. (ML);
    \draw (L2) .. controls +( 0.3, 0) and +(0,0) .. (ML);
  \end{scope}  
 \begin{scope}[yshift = -1cm]
    \coordinate (L1B) at (0.2,0.4);
    \coordinate (L2B) at (0,0);
    \coordinate (R1B) at (2.2,0.4);
    \coordinate (R2B) at (2,0);
    \draw[->-] (L1B) .. controls +( 0, 0) and +(0,0) .. (R1B) node [right, pos = 1] {$a$} node[above,pos=0.6] {$\deco_n$};
    \draw[->-] (L2B) .. controls +( 0, 0) and +(0,0) .. (R2B) node [right, pos = 1] {$b$} node[above,pos=0.4] {$\deco_0$};
 \end{scope}  
  \draw (R1) -- (R1B);
  \draw (R2) -- (R2B);
  \draw (L1) -- (L1B);
  \draw (L2) -- (L2B);
  \draw (ML) .. controls +(0, -0.6) and +(0, -0.6) .. (MR);
\end{scope}

%% file: sch/foam_zip_ab_kl.tex
\begin{scope}
  \begin{scope}
    \coordinate (L1) at (0.2,0.4);
    \coordinate (L2) at (0,0);
    \coordinate (R1) at (2.2,0.4);
    \coordinate (R2) at (2,0);
    \coordinate (ML) at (0.6, 0.2);
    \coordinate (MR) at (1.6, 0.2);
    \draw[->-] (ML) -- (MR) node[above, midway] {$a+b$};
    \draw (MR) .. controls +(0, 0) and +(-0.3,0) .. (R1) ;
    \draw (MR) .. controls +(0, 0) and +(-0.3,0) .. (R2);
    \draw (L1) .. controls +( 0.3, 0) and +(0,0) .. (ML);
    \draw (L2) .. controls +( 0.3, 0) and +(0,0) .. (ML);
  \end{scope}  
 \begin{scope}[yshift = -1cm]
    \coordinate (L1B) at (0.2,0.4);
    \coordinate (L2B) at (0,0);
    \coordinate (R1B) at (2.2,0.4);
    \coordinate (R2B) at (2,0);
    \draw[->-] (L1B) .. controls +( 0, 0) and +(0,0) .. (R1B) node [right, pos = 1] {$a$} node[above,pos=0.6] {$\deco_k$};
    \draw[->-] (L2B) .. controls +( 0, 0) and +(0,0) .. (R2B) node [right, pos = 1] {$b$} node[above,pos=0.4] {$\deco_l$};
 \end{scope}  
  \draw (R1) -- (R1B);
  \draw (R2) -- (R2B);
  \draw (L1) -- (L1B);
  \draw (L2) -- (L2B);
  \draw(ML) .. controls +(0, -0.6) and +(0, -0.6) .. (MR);
\end{scope}

%% file: sch/foam_unzip_ab_arrows.tex
\begin{scope}
  \begin{scope}
    \coordinate (L1) at (0.2,0.4);
    \coordinate (L2) at (0,0);
    \coordinate (R1) at (2.2,0.4);
    \coordinate (R2) at (2,0);
    \coordinate (ML) at (0.6, 0.2);
    \coordinate (MR) at (1.6, 0.2);
    \draw[->-] (ML) -- (MR) node[below, midway] {$a+b$};
    \draw (MR) .. controls +(0, 0) and +(-0.3,0) .. (R1) ;
    \draw (MR) .. controls +(0, 0) and +(-0.3,0) .. (R2);
    \draw (L1) .. controls +( 0.3, 0) and +(0,0) .. (ML);
    \draw (L2) .. controls +( 0.3, 0) and +(0,0) .. (ML);
  \end{scope}  
 \begin{scope}[yshift = 1cm]
    \coordinate (L1B) at (0.2,0.4);
    \coordinate (L2B) at (0,0);
    \coordinate (R1B) at (2.2,0.4);
    \coordinate (R2B) at (2,0);
    \draw[->-] (L1B) .. controls +( 0, 0) and +(0,0) .. (R1B) node [left, pos
    = 0] {$a$};
    \draw[->-] (L2B) .. controls +( 0, 0) and +(0,0) .. (R2B) node [left, pos
    = 0] {$b$};
 \end{scope}  
  \draw (R1) -- (R1B);
  \draw (R2) -- (R2B);
  \draw (L1) -- (L1B);
  \draw (L2) -- (L2B);
  \draw (ML) .. controls +(0, 0.6) and +(0, 0.6) .. (MR);
\end{scope}

%% file: sch/foam_unzip_ab_0n.tex
\begin{scope}
  \begin{scope}
    \coordinate (L1) at (0.2,0.4);
    \coordinate (L2) at (0,0);
    \coordinate (R1) at (2.2,0.4);
    \coordinate (R2) at (2,0);
    \coordinate (ML) at (0.6, 0.2);
    \coordinate (MR) at (1.6, 0.2);
    \draw[->-] (ML) -- (MR) node[below, midway] {$a+b$};
    \draw (MR) .. controls +(0, 0) and +(-0.3,0) .. (R1) ;
    \draw (MR) .. controls +(0, 0) and +(-0.3,0) .. (R2);
    \draw (L1) .. controls +( 0.3, 0) and +(0,0) .. (ML);
    \draw (L2) .. controls +( 0.3, 0) and +(0,0) .. (ML);
  \end{scope}  
 \begin{scope}[yshift = 1cm]
    \coordinate (L1B) at (0.2,0.4);
    \coordinate (L2B) at (0,0);
    \coordinate (R1B) at (2.2,0.4);
    \coordinate (R2B) at (2,0);
    \draw[->-] (L1B) .. controls +( 0, 0) and +(0,0) .. (R1B) node [left, pos
    = 0] {$a$} node[pos=0.75,below] {$\deco_0$};
    \draw[->-] (L2B) .. controls +( 0, 0) and +(0,0) .. (R2B) node [left, pos
    = 0] {$b$} node[pos=0.35,below] {$\deco_n$};
 \end{scope}  
  \draw (R1) -- (R1B);
  \draw (R2) -- (R2B);
  \draw (L1) -- (L1B);
  \draw (L2) -- (L2B);
  \draw (ML) .. controls +(0, 0.6) and +(0, 0.6) .. (MR);
\end{scope}

%% file: sch/foam_unzip_ab_n0.tex
\begin{scope}
  \begin{scope}
    \coordinate (L1) at (0.2,0.4);
    \coordinate (L2) at (0,0);
    \coordinate (R1) at (2.2,0.4);
    \coordinate (R2) at (2,0);
    \coordinate (ML) at (0.6, 0.2);
    \coordinate (MR) at (1.6, 0.2);
    \draw[->-] (ML) -- (MR) node[below, midway] {$a+b$};
    \draw (MR) .. controls +(0, 0) and +(-0.3,0) .. (R1) ;
    \draw (MR) .. controls +(0, 0) and +(-0.3,0) .. (R2);
    \draw (L1) .. controls +( 0.3, 0) and +(0,0) .. (ML);
    \draw (L2) .. controls +( 0.3, 0) and +(0,0) .. (ML);
  \end{scope}  
 \begin{scope}[yshift = 1cm]
    \coordinate (L1B) at (0.2,0.4);
    \coordinate (L2B) at (0,0);
    \coordinate (R1B) at (2.2,0.4);
    \coordinate (R2B) at (2,0);
    \draw[->-] (L1B) .. controls +( 0, 0) and +(0,0) .. (R1B) node [left, pos
    = 0] {$a$} node[pos=0.75,below] {$\deco_n$};
    \draw[->-] (L2B) .. controls +( 0, 0) and +(0,0) .. (R2B) node [left, pos
    = 0] {$b$} node[pos=0.35,below] {$\deco_0$};
 \end{scope}  
  \draw (R1) -- (R1B);
  \draw (R2) -- (R2B);
  \draw (L1) -- (L1B);
  \draw (L2) -- (L2B);
  \draw (ML) .. controls +(0, 0.6) and +(0, 0.6) .. (MR);
\end{scope}

%% file: sch/foam_unzip_ab_kl.tex
\begin{scope}
  \begin{scope}
    \coordinate (L1) at (0.2,0.4);
    \coordinate (L2) at (0,0);
    \coordinate (R1) at (2.2,0.4);
    \coordinate (R2) at (2,0);
    \coordinate (ML) at (0.6, 0.2);
    \coordinate (MR) at (1.6, 0.2);
    \draw[->-] (ML) -- (MR) node[below, midway] {$a+b$};
    \draw (MR) .. controls +(0, 0) and +(-0.3,0) .. (R1) ;
    \draw (MR) .. controls +(0, 0) and +(-0.3,0) .. (R2);
    \draw (L1) .. controls +( 0.3, 0) and +(0,0) .. (ML);
    \draw (L2) .. controls +( 0.3, 0) and +(0,0) .. (ML);
  \end{scope}  
 \begin{scope}[yshift = 1cm]
    \coordinate (L1B) at (0.2,0.4);
    \coordinate (L2B) at (0,0);
    \coordinate (R1B) at (2.2,0.4);
    \coordinate (R2B) at (2,0);
    \draw[->-] (L1B) .. controls +( 0, 0) and +(0,0) .. (R1B) node [left, pos
    = 0] {$a$} node[pos=0.75,below] {$\deco_k$};
    \draw[->-] (L2B) .. controls +( 0, 0) and +(0,0) .. (R2B) node [left, pos
    = 0] {$b$} node[pos=0.35,below] {$\deco_l$};
 \end{scope}  
  \draw (R1) -- (R1B);
  \draw (R2) -- (R2B);
  \draw (L1) -- (L1B);
  \draw (L2) -- (L2B);
  \draw (ML) .. controls +(0, 0.6) and +(0, 0.6) .. (MR);
\end{scope}

%% file: sch/foam_cup_a_empty.tex
\begin{scope}
  \draw (0,0) arc (180 :0: 0.5cm and 0.2cm) 
  node[above, pos=0.5] {$a$};
  \draw[very thin] (0,0) arc (180 :0: 0.5cm and -0.6cm) node[pos=0.5, above] {};
  \draw (0,0) arc (180 :0: 0.5cm and -0.2cm);
\end{scope}

%% file: sch/foam_cup_a_kl.tex
\begin{scope}
  \draw (0,0) arc (180 :0: 0.5cm and 0.2cm) 
  node[above, pos=0.5] {$a$};
  \draw[very thin] (0,0) arc (180 :0: 0.5cm and -0.6cm);
  \draw (0,0) arc (180 :0: 0.5cm and -0.2cm) node[pos=0.5, below] {$\deco_k \antideco_l$};
\end{scope}

%% file: sch/foam_cap_a_empty.tex
\begin{scope}[-]
  \draw (0,0) arc (180 :0: 0.5cm and 0.2cm);
  \draw[very thin] (0,0) arc (180 :0: 0.5cm and 0.6cm) node[pos=0.5,
  below] {};
  \draw (0,0) arc (180 :0: 0.5cm and -0.2cm)node[ below, pos =0.5] {$a$};
\end{scope}

%% file: sch/foam_cap_a_kl.tex
\begin{scope}[-]
  \draw (0,0) arc (180 :0: 0.5cm and 0.2cm);
  \draw[very thin] (0,0) arc (180 :0: 0.5cm and 0.6cm) node[pos=0.5, below] {$\deco_k \antideco_l$};
  \draw (0,0) arc (180 :0: 0.5cm and -0.2cm)node[ below, pos =0.5] {$a$};
\end{scope}

%% file: sch/foam_saddle_a_empty.tex
\begin{scope}[scale=0.6]
\tdplotsetmaincoords{70}{25}
\begin{scope}[scale = 1.5, tdplot_main_coords]
  \tikzset{yxplane/.style={canvas is xy plane at z=#1}}
  \begin{scope}[yxplane=1]
    \coordinate (AT) at ({cos(  45)}, {sin(  45)});
    \coordinate (BT) at ({cos(135)}, {sin(135)});
    \coordinate (CT) at ({cos(225)}, {sin(225)});
    \coordinate (DT) at ({cos(315)}, {sin(315)});
    \coordinate (aT) at (  0.3, 0);
    \coordinate (bT) at ( -0.3, 0);
    \draw (AT) .. controls (aT) and (bT) .. (BT)
    coordinate[pos=0.5] (eT) node [pos =0.5, above] {$a$};
    \draw (DT) .. controls (aT) and (bT) .. (CT) coordinate[pos=0.5] (fT);
  \end{scope}

  \begin{scope}[yxplane=0]
    \coordinate (AM) at ({cos(  45)}, {sin(  45)});
    \coordinate (BM) at ({cos(135)}, {sin(135)});
    \coordinate (CM) at ({cos(225)}, {sin(225)});
    \coordinate (DM) at ({cos(315)}, {sin(315)});
    \coordinate (aM) at (0, 0.3);
    \coordinate (bM) at (0, -0.3);
    \draw (AM) .. controls (aM) and (bM) .. (DM) coordinate[pos=0.5] (eM);
    \draw (BM) .. controls (aM) and (bM) .. (CM) coordinate[pos=0.5] (fM);
\end{scope}
  \coordinate (OT) at (0,0, 0.5);
\draw (AM) -- (AT);
  \draw (BM) -- (BT);
  \draw (CM) -- (CT);
  \draw (DM) -- (DT);
  \draw[densely dotted] (eM) ..controls +(0,0,0.2) and + (0.1,0,0).. (OT);
  \draw[densely dotted] (fM) ..controls +(0,0,0.2) and + (-0.1,0,0).. (OT);
  \draw[densely dotted] (eT) ..controls +(0,0,-0.2) and + (0,0.1,0).. (OT);
  \draw[densely dotted] (fT) ..controls +(0,0,-0.2) and + (0,-0.1,0).. (OT);
\end{scope}  
\end{scope}

%% file: sch/foam_saddle_a_kl.tex
\begin{scope}[scale=0.6]
\tdplotsetmaincoords{70}{25}
\begin{scope}[scale = 1.5, tdplot_main_coords]
  \tikzset{yxplane/.style={canvas is xy plane at z=#1}}
  \begin{scope}[yxplane=1]
    \coordinate (AT) at ({cos(  45)}, {sin(  45)});
    \coordinate (BT) at ({cos(135)}, {sin(135)});
    \coordinate (CT) at ({cos(225)}, {sin(225)});
    \coordinate (DT) at ({cos(315)}, {sin(315)});
    \coordinate (aT) at (  0.3, 0);
    \coordinate (bT) at ( -0.3, 0);
    \draw (AT) .. controls (aT) and (bT) .. (BT)
    coordinate[pos=0.5] (eT) node [pos =0.5, above] {$a$};
    \draw (DT) .. controls (aT) and (bT) .. (CT) coordinate[pos=0.5] (fT);
  \end{scope}

  \begin{scope}[yxplane=0]
    \coordinate (AM) at ({cos(45)}, {sin(45)});
    \coordinate (BM) at ({cos(135)}, {sin(135)});
    \coordinate (CM) at ({cos(225)}, {sin(225)});
    \coordinate (DM) at ({cos(315)}, {sin(315)});
    \coordinate (aM) at (0, 0.3);
    \coordinate (bM) at (0, -0.3);
    \draw (AM) .. controls (aM) and (bM) .. (DM) coordinate[pos=0.5] (eM);
    \draw (BM) .. controls (aM) and (bM) .. (CM) coordinate[pos=0.5] (fM);
\end{scope}
  \coordinate (OT) at (0,0, 0.5);
    \draw (AM) -- (AT) node[pos=.5,left] {$\antideco_\ell$};
    \draw (BM) -- (BT) ;
    \draw (CM) -- (CT) node[pos=.5,right] {$\deco_k$};
    \draw (DM) -- (DT);
  \draw[densely dotted] (eM) ..controls +(0,0,0.2) and + (0.1,0,0).. (OT);
  \draw[densely dotted] (fM) ..controls +(0,0,0.2) and + (-0.1,0,0).. (OT);
  \draw[densely dotted] (eT) ..controls +(0,0,-0.2) and + (0,0.1,0).. (OT);
  \draw[densely dotted] (fT) ..controls +(0,0,-0.2) and + (0,-0.1,0).. (OT);
\end{scope}  
\end{scope}

%% file: SF_twists_on_foams.tex
In the language of Section \ref{twist_gen},
Theorem \ref{witt-thm} implies that
\(\catf\) is a \(\Witt\)-category.

\begin{lem}
A foam $\foam \in \catf(\web, \web^\prime)$,
induces a \(\Witt\)-equivariant map \(\ssp_N(\foam): \ssp_N(\web) \to \ssp_N(\web^\prime)\) if and only if \(\Witt\) acts trivially on \(F\).
\end{lem}

\begin{proof}
For any \(X \in \Witt\) and \(v \in \ssp_N(\web)\),
\(\foam (X \cdot v) = X \cdot (\foam v) - (X \cdot \foam) v.\) 
Thus, a map $\ssp_N(\foam)$ intertwines the actions of $\Witt$ if and only if for all \(X \in \Witt\), \(X \cdot F = 0\).
\end{proof}

\begin{cor}
    \label{lemm_eq}
    For any \(\foam \in \Hom_{\text{Tw}(\catf)}((\web,\tau),(\web^\prime,\tau^\prime))\), the induced map 
    \(\ssp_N(\foam):\ssp_N(\web,\tau) \to \ssp_N(\web^\prime, \tau^\prime)\) is
    $\Witt$-equivariant if and only if, for every \(\wgen_n \in \Witt\)
    \begin{equation}
        \label{eq_action}
        \Wgen_n \foam + \tau_{\wgen_n} \circ \foam - \foam \circ \tau^\prime_{\wgen_n} = 0.
    \end{equation}
\end{cor}


\begin{dfn}
For each edge $e$ in $\web$ with thickness \(a\),
let \(D_e := \scalars_{a, N-a}\)
be the algebra of polynomial decorations on that edge and
let
\[D_{\web} := \bigotimes_{e \in E(\web)} D_e\]
be the algebra of all decorations on the web.
Let \({d: D_{\web} \to \catf(\web,\web)}\).
be the \(\Witt\)-equivariant map defined by
sending a decoration
to the 
cylinder 
decorated by it.
\end{dfn}

Since the algebra of decorations is commutative,
the condition for a Maurer--Cartan element \(\alpha\) which lands in \(d(D_\web) \subset \catf(\web,\web)\) reduces to
\[
(n-m)\alpha_{\wgen_{n+m}} = {\wgen_n}\alpha_{\wgen_m} - {\wgen_m}\alpha_{\wgen_n}.
\]
In \cite{QRSW23}, this condition was called \textit{\(\Witt\)-flatness} for the map \[ \function{\alpha}{\Witt}{D_\web}{\wgen_n}{\alpha_{\wgen_n}}.\]

%% file: SF_red_dots_alternate.tex
The twist maps above are local in the sense that they involve decorations on a single edge.
Therefore they can be conveniently tracked by a diagrammatic formalism.
Let \(p_k^{(e)}\) and \(\widehat{p_k}^{(e)}\) denote the decorations \(p_k\) and \(\widehat{p_k}\) in \(R_{a,N-a} = D_e \subset D_\web\).

\begin{dfn}
\label{def_dots}
A \textit{red-dotted} web is a web \(\web\) endowed with a finite collection \(\cold\) of \textit{red dots} of the following three types on the interiors of the edges of \(\web\):
\begin{itemize}
	\item A hollow red dot $\rd$, associated to the twist
	\[\tau_1: L_n \mapsto -\sum_{k=0}^np_k^{(e)}p_{n-k}^{(e)}.\]
	\item A solid red dot $\srd$, associated to the twist
	\[\tau_2: L_n \mapsto -\sum_{k=0}^np^{(e)}_k\widehat{p}^{(e)}_{n-k}.\]
	\item A triangular red dot $\trd$, associated to the twist
	\[\tau_3: L_n \mapsto -(n+1) p^{(e)}_n.\]
\end{itemize}
\end{dfn}
The linearity of twists implies that multiple red dots can lie on the same edge with \(\scalars\)-valued labels representing linear combinations.

For three elements \(\omega_1, \omega_2, \lambda \in \scalars\), the twisted action can be described by local pictures:
\begin{align}
    \wgen_n \cdot \, \NB{\tikz[font=\small, scale=1.5]{\input{sch/def_dots_circle_hollow}}} \, &=
    \wgen_n \NB{\tikz[font=\small, scale=1.5]{\input{sch/def_dots_blank}}}
    -\omega_1 \sum_{k=0}^n \NB{\tikz[font=\small, scale=1.5]{\input{sch/def_dots_hollow_blank}}}
    \\ 
    \wgen_n \cdot \, \NB{\tikz[font=\small, scale=1.5]{\input{sch/def_dots_circle_solid}}} \, &=
    \wgen_n \NB{\tikz[font=\small, scale=1.5]{\input{sch/def_dots_blank}}}
    -\omega_2 \sum_{k=0}^n \NB{\tikz[font=\small, scale=1.5]{\input{sch/def_dots_solid_blank}}}
    \\
    \wgen_n \cdot \, \NB{\tikz[font=\small, scale=1.5]{\input{sch/def_dots_trian_hollow}}} \, &=
    \wgen_n \NB{\tikz[font=\small, scale=1.5]{\input{sch/def_dots_blank}}}
    -\lambda(n+1) ~ \NB{\tikz[font=\small, scale=1.5]{\input{sch/def_dots_trian_hollow_blank}}} 
\end{align} 
Following Proposition \ref{prop_MC}, the twists on the target web are subtracted, while the twists on the source web are added. 

\begin{lem}
    \label{flatness}
    The maps in Definition (\ref{def_dots}) define Maurer-Cartan elements.
\end{lem}

\begin{proof}
Since each twist lands in \(D_\web\), the Maurer-Cartan condition reduces to flatness.
    Assume without loss of generality that $n\leq m$. Then
    \begin{align*}
        \wgen_n \cdot \tau_1(\wgen_m) - \wgen_m \cdot \tau_1 (\wgen_n) &= \sum_{i=0}^n (m-n) p^{(e)}_i p^{(e)}_{n+m-i} +  \sum_{i=n+1}^m (m-i) p^{(e)}_ip^{(e)}_{n+m-i} \\
	&\quad-  \sum_{i=n}^{m+n} (n-i) p^{(e)}_ip^{(e)}_{n+m-i} -  \sum_{i=m}^{n+m} (i-m) p^{(e)}_ip^{(e)}_{n+m-i}\\
	&= \sum_{i=0}^n (m-n) p^{(e)}_i p^{(e)}_{n+m-i} +  \sum_{i=m}^{m+n} (m-n) p^{(e)}_ip^{(e)}_{n+m-i} \\
	&\quad+  \sum_{i=n}^{m-1} (i-n) p^{(e)}_ip^{(e)}_{n+m-i} +  \sum_{i=n+1}^{m} (m-i) p^{(e)}_ip^{(e)}_{n+m-i} \\
        &= (n-m) \tau_1(\wgen_{m+n}). 
    \end{align*}
    Thus $\tau_1$ is flat. The computation for $\tau_2$ is similar and is omitted.
    The computation showing that $\tau_3$ is flat is as follows.
    \begin{align*}
        \wgen_n \cdot \tau_3(\wgen_m) - \wgen_m \cdot \tau_3 (\wgen_n) & = (m(m+1) - n(n+1)) p^{(e)}_{n+m} \\
        & = (n-m) \tau_3(\wgen_{m+n}).
    \end{align*}
\end{proof}

\begin{rmk}
    The triangular red dot labelled by $\frac{\lambda}{2} \in \scalars$ agrees with the hollow green dots labelled by $\lambda$ for \(\sld\) defined in \cite{QRSW23}. 
\end{rmk}

\begin{exa}
The collection of twists which make a fixed foam equivariant is usually not unique.
For example, the pair below of twists on the unzip foam induce \(\Witt\)-equivariant morphisms.
       \[\NB{\tikz[font=\tiny, scale=.5]{\input{sch/ex_sliding_1}}} \rightarrow \NB{\tikz[font=\tiny,scale=.5]{\input{sch/ex_sliding_2}}} \hspace{5mm} \text{and} \hspace{5mm}\NB{\tikz[font=\tiny,scale=.5]{\input{sch/ex_sliding_3}}} \rightarrow \NB{\tikz[font=\tiny,scale=.5]{\input{sch/ex_sliding_4}}} \ . \]
This illustrates that red dots can migrate between the source to the target along a shared facet while picking up a sign.
The situation is analogous to the fact that grading shifts are relative and not absolute.
Since twists are local, more variation is possible.
\end{exa}

\begin{lem}
    \label{lemm_rd_1}
    The Leibniz rule implies the following relations between twists: 
    \begin{gather}
        \NB{\tikz[font=\small]{\input{sch/relation_add_dots1}}} = \NB{\tikz[font=\small]{\input{sch/relation_add_dots2}}}, \\
        \NB{\tikz[font=\small]{\input{sch/relation_add_dots3}}} = \NB{\tikz[font=\small]{\input{sch/relation_add_dots4}}}, \\
        \NB{\tikz[font=\small]{\input{sch/relation_add_dots5}}} = \NB{\tikz[font=\small]{\input{sch/relation_add_dots6}}}.
    \end{gather}
\end{lem}

\begin{dfn} The following red dots abbreviate twists around trivalent vertices.
\begin{gather}
	\NB{\tikz[font=\small, scale=0.5]{\input{sch/notation_reddot2}}} :=
	\NB{\tikz[font=\small, scale=0.5]{\input{sch/notation_reddot1}}}, \hspace{3mm}
	\NB{\tikz[font=\small, yscale=-1, scale=0.5]{\input{sch/notation_reddot4}}} :=
	\NB{\tikz[font=\small, yscale=-1, scale=0.5]{\input{sch/notation_reddot3}}}, \\
	\NB{\tikz[font=\small, scale=0.5]{\input{sch/notation_reddot2_solid}}} :=
	\NB{\tikz[font=\small, scale=0.5]{\input{sch/notation_reddot1_solid}}}, \hspace{3mm}
	\NB{\tikz[font=\small, yscale=-1, scale=0.5]{\input{sch/notation_reddot4_solid}}} :=
	\NB{\tikz[font=\small, yscale=-1, scale=0.5]{\input{sch/notation_reddot3_solid}}}.
\end{gather}
\end{dfn}

\begin{lem}
\label{lemm_rd_3}
Dot migration rules \cite[(11)]{RW20a} imply the following relation between twists:
\begin{gather}
	\label{lem_lever}
	\NB{\tikz[font=\small, scale=0.5]{\input{sch/notation_reddot2_solid}}} = \NB{\tikz[font=\small, scale=0.5]{\input{sch/notation_reddot2}}},
	\hspace{3mm}
	\NB{\tikz[font=\small, yscale=-1, scale=0.5]{\input{sch/notation_reddot4_solid}}} = \NB{\tikz[font=\small, yscale=-1, scale=0.5]{\input{sch/notation_reddot4}}},\\
	\NB{\tikz[font=\small, scale=.7]{\input{sch/reidem_lemm_21}}} = \NB{\tikz[font=\small, scale=0.7]{\input{sch/reidem_lemm_22}}} \hspace{4mm}\text{and} \hspace{6mm} \NB{\tikz[font=\small, scale=0.7]{\input{sch/relation_rd_merge}}} = \NB{\tikz[font=\small, scale=0.7]{\input{sch/relation_rd_split}}}.
\end{gather}
\end{lem}

\begin{lem}
\label{circ_tri}
On a facet of thickness one, the dot twist reduces to the triangle twist:
\begin{equation}
	\NB{\tikz[font=\small]{\input{sch/lemm_circ_tri_1}}} = \NB{\tikz[font=\small]{\input{sch/lemm_circ_tri_2}}}.
\end{equation}
\end{lem}

%% file: sch/def_dots_circle_hollow.tex
\begin{scope}
  \draw[thin] (0,0) rectangle (1,0.7);
  \draw[thick] (0,0.7) -- (1,0.7)  coordinate [pos=0.5] (A);
  \filldraw[draw= \crd, fill = white] (A) circle (.7mm)
  node[below, yshift=-0.5mm] {$\omega_1$};
\end{scope}

%% file: sch/def_dots_blank.tex
\begin{scope}
  \draw[thin] (0,0) rectangle (1,0.7);
  \draw[thick] (0,0.7) -- (1,0.7)  coordinate [pos=0.5] (A);
\end{scope}

%% file: sch/def_dots_hollow_blank.tex
\begin{scope}
  \draw[thin] (0,0) rectangle (1,0.7);
  \draw[thick] (0,0.7) -- (1,0.7)  coordinate [pos=0.5] (A);
  \node (B) at (.55,.3) {$\deco_k \deco_{n-k}$};
\end{scope}

%% file: sch/def_dots_circle_solid.tex
\begin{scope}
  \draw[thin] (0,0) rectangle (1,0.7);
  \draw[thick] (0,0.7) -- (1,0.7)  coordinate [pos=0.5] (A);
  \filldraw[draw= \crd, fill = \csrd] (A) circle (.7mm)
  node[below, yshift=-0.5mm] {$\omega_2$};
\end{scope}

%% file: sch/def_dots_solid_blank.tex
\begin{scope}
  \draw[thin] (0,0) rectangle (1,0.7);
  \draw[thick] (0,0.7) -- (1,0.7)  coordinate [pos=0.5] (A);
  \node (B) at (.55,.3) {$\deco_k \hat{\deco}_{n-k}$};
\end{scope}

%% file: sch/def_dots_trian_hollow.tex
\begin{scope}
  \draw[thin] (0,0) rectangle (1,0.7);
  \draw[thick] (0,0.7) -- (1,0.7)  coordinate [pos=0.5] (A); 
  \node (C) at (0.5,0.7) [draw=\crd, fill=white,regular polygon, regular polygon sides=3,inner sep=1.5pt]{};
  \node (D) at (0.5,0.7) [below,yshift=-0.5mm]{$\lambda$};   
\end{scope}

%% file: sch/def_dots_trian_hollow_blank.tex
\begin{scope}
  \draw[thin] (0,0) rectangle (1,0.7);
  \draw[thick] (0,0.7) -- (1,0.7)  coordinate [pos=0.5] (A); 
  \node (B) at (.55,.3) {$\deco_n$};
\end{scope}

%% file: sch/ex_sliding_2.tex
\begin{scope}
    \coordinate (bl) at (-0.5, -1);
    \coordinate (br) at ( 0.5, -1);
    \coordinate (bm) at (  0,-0.3);
    \coordinate (tl) at (-.5,  1);
    \coordinate (tr) at ( 0.5,  1);
    \coordinate (tm) at (  0, 0.3);
    \draw[->, color=blue] (br) -- (tr);
    \draw[->] (bl) -- (tl);
    \draw (-.5, 1) arc (-180:0:-0.5) -- (-1.5,-1) arc (-180:0:0.5);
    \draw[color = blue] (.5, 1) arc (180:0:0.5)  -- (1.5,-1) arc (0:-180:0.5);
\end{scope}

%% file: sch/ex_sliding_4.tex
\begin{scope}
    \coordinate (bl) at (-0.5, -1);
    \coordinate (br) at ( 0.5, -1);
    \coordinate (bm) at (  0,-0.3);
    \coordinate (tl) at (-.5,  1);
    \coordinate (tr) at ( 0.5,  1);
    \coordinate (tm) at (  0, 0.3);
    \draw[->,color=blue] (br) -- (tr);
    \draw[->] (bl) -- (tl);
    \draw (-.5, 1) arc (-180:0:-0.5) -- (-1.5,-1) arc (-180:0:0.5);
    \draw[color=blue] (.5, 1) arc (180:0:0.5) -- (1.5,-1) arc (0:-180:0.5) node[pos=.5, draw=\crd, fill= white, circle , inner sep=1.5pt]{} node[pos=.5, above]{$\frac{s}{2}$};
\end{scope}

%% file: sch/relation_add_dots1.tex
\begin{scope}
  \coordinate (m) at (0,0);
  \coordinate (t) at (2,0);
  \draw[->] (m) -- (t) node[pos = 1, right] {$a$} coordinate[pos=0.3] (ga) coordinate[pos=0.7] (gb);
  \filldraw[draw= \crd, fill = white] (ga) circle (1mm) node[yshift= 3mm] {$\omega_1$};
  \filldraw[draw= \crd, fill = white] (gb) circle (1mm) node[yshift=3mm] {$\omega_2$};
\end{scope}

%% file: sch/relation_add_dots2.tex
\begin{scope}
  \coordinate (m) at (0,0);
  \coordinate (t) at (2,0);
  \draw[->] (m) -- (t) node[pos = 1, right] {$a$} coordinate[pos=0.5] (ga) coordinate[pos=0.7] (gb);
  \filldraw[draw= \crd, fill = white] (ga) circle (1mm) node[yshift= 3mm] {$\omega_1+\omega_2$};
\end{scope}

%% file: sch/relation_add_dots3.tex
\begin{scope}
  \coordinate (m) at (0,0);
  \coordinate (t) at (2,0);
  \draw[->] (m) -- (t) node[pos = 1, right] {$a$} coordinate[pos=0.3] (ga) coordinate[pos=0.7] (gb);
  \filldraw[draw= \crd, fill = \csrd] (ga) circle (1mm) node[yshift= 3mm] {$\omega_1$};
  \filldraw[draw= \crd, fill = \csrd] (gb) circle (1mm) node[yshift=3mm] {$\omega_2$};
\end{scope}

%% file: sch/relation_add_dots4.tex
\begin{scope}
  \coordinate (m) at (0,0);
  \coordinate (t) at (2,0);
  \draw[->] (m) -- (t) node[pos = 1, right] {$a$} coordinate[pos=0.5] (ga) coordinate[pos=0.7] (gb);
  \filldraw[draw= \crd, fill = \csrd] (ga) circle (1mm) node[yshift= 3mm] {$\omega_1+\omega_2$};
\end{scope}

%% file: sch/relation_add_dots6.tex
\begin{scope}
  \coordinate (m) at (0,0);
  \coordinate (t) at (2,0);
  \draw[->] (m) -- (t) node[pos = 1, right] {$a$} coordinate[pos=0.5] (ga) node[pos= 0.5 , draw=\crd, fill=white,regular polygon, regular polygon sides=3,inner sep=1.2pt] {} node[pos=0.5, above, yshift=0.5mm]{$\omega_1 + \omega_2$};
\end{scope}

%% file: sch/notation_reddot2.tex
\begin{scope}
  \coordinate (a) at (1,0); 
  \coordinate (b) at (-1,0);
  \coordinate (c) at (0,-2);
  \coordinate (d) at (0,-1);
  \draw[->] (d) -- (a) node[pos=1, above] {$b$};
  \draw[->] (d) -- (b) node[pos=1, above] {$a$};
  \draw[>-] (c) -- (d) node[pos=0, below] {$a+b$} node[pos=1, draw=\crd, fill= white, circle , inner sep=1.5pt]{} node[pos=1, right]{$\omega$};
\end{scope}

%% file: sch/notation_reddot4.tex
\begin{scope}
  \coordinate (a) at (1,0); 
  \coordinate (b) at (-1,0);
  \coordinate (c) at (0,-2);
  \coordinate (d) at (0,-1);
  \draw[-<] (d) -- (a) node[pos=1, below] {$b$};
  \draw[-<] (d) -- (b) node[pos=1, below] {$a$};
  \draw[<-] (c) -- (d) node[pos=0, above] {$a+b$} node[pos=1, draw=\crd, fill= white, circle , inner sep=1.5pt]{} node[pos=1, right]{$\omega$};
\end{scope}

%% file: sch/notation_reddot2_solid.tex
\begin{scope}
  \coordinate (a) at (1,0); 
  \coordinate (b) at (-1,0);
  \coordinate (c) at (0,-2);
  \coordinate (d) at (0,-1);
  \draw[->] (d) -- (a) node[pos=1, above] {$b$};
  \draw[->] (d) -- (b) node[pos=1, above] {$a$};
  \draw[>-] (c) -- (d) node[pos=0, below] {$a+b$} node[pos=1, draw=\crd, fill= \csrd, circle , inner sep=1.5pt]{} node[pos=1, right]{$\omega$};
\end{scope}

%% file: sch/notation_reddot4_solid.tex
\begin{scope}
  \coordinate (a) at (1,0); 
  \coordinate (b) at (-1,0);
  \coordinate (c) at (0,-2);
  \coordinate (d) at (0,-1);
  \draw[-<] (d) -- (a) node[pos=1, below] {$b$};
  \draw[-<] (d) -- (b) node[pos=1, below] {$a$};
  \draw[<-] (c) -- (d) node[pos=0, above] {$a+b$} node[pos=1, draw=\crd, fill= \csrd, circle , inner sep=1.5pt]{} node[pos=1, right]{$\omega$};
\end{scope}

%% file: sch/lemm_circ_tri_1.tex
\begin{scope}
  \coordinate (m) at (0,0);
  \coordinate (t) at (2,0);
  \draw[->] (m) -- (t) node[pos = 1, right,font=\tiny] {$1$} coordinate[pos=0.5] (ga) node[pos= 0.5 , draw=\crd, fill=white,regular polygon, regular polygon sides=3,inner sep=1.2pt] {} node[pos=0.5, above, yshift=0.5mm]{$\omega$};
\end{scope}

%% file: sch/lemm_circ_tri_2.tex
\begin{scope}
  \coordinate (m) at (0,0);
  \coordinate (t) at (2,0);
  \draw[->] (m) -- (t) node[pos = 1, right,font=\tiny] {$1$} coordinate[pos=0.5] (ga) coordinate[pos=0.7] (gb);
  \filldraw[draw= \crd, fill = white] (ga) circle (1mm) node[yshift= 3mm] {$\omega$};
\end{scope}

%% file: SF_usual_morph.tex
\begin{lem}
    \label{rd_asso}
    The two orientations for the following foam:  
     \[\NB{\tikz[font=\tiny, scale=0.7]{\input{sch/foam_associativity}}} \]
    induce $\Witt$-equivariant isomorphisms:
    \[\mapasso : \NB{\tikz[font=\tiny, scale=0.7]{\input{sch/web_associativity1}}} \rightarrow \NB{\tikz[font=\tiny, scale=0.7,xscale=-1]{\input{sch/web_associativity2}}} \hspace{5.5mm} \text{and} \hspace{5.5mm} \mapasso : \NB{\tikz[font=\tiny, scale=0.7]{\input{sch/web_associativity3}}} \rightarrow \NB{\tikz[font=\tiny, scale=0.7]{\input{sch/web_associativity4}}}. \]
\end{lem}

\begin{lem}
    \label{rd_digon}
    The foams:
    \[\NB{\tikz[font=\tiny, scale=1]{\input{sch/foam_digoncap_empty}}} \hspace{5.5mm} \text{and} \hspace{5.5mm} \NB{\tikz[font=\tiny, scale=1]{\input{sch/foam_digoncup_empty}}}\]
    induce $\Witt$-equivariant morphisms:
    \[\mapdcap : \NB{\tikz[font=\tiny, scale=0.7]{\input{sch/web_digoncap_beg}}} \rightarrow  \NB{\tikz[font=\tiny, scale=0.7]{\input{sch/web_l_ab}}} \hspace{5.5mm} \text{and} \hspace{5.5mm} \mapdcup : \NB{\tikz[font=\tiny, scale=0.7]{\input{sch/web_l_ab}}} \rightarrow \NB{\tikz[font=\tiny, scale=0.7]{\input{sch/web_digoncup_end}}}. \]

\end{lem}

\begin{proof}
    \begin{align*}
        \wgen_n & \left ( \hspace{2mm} \NB{\tikz[font=\tiny, scale=1]{\input{sch/foam_digoncup_empty}}} \right ) = \lambda_n ~~~ \NB{\tikz[font=\tiny, scale=1]{\input{sch/foam_digoncup_0n}}} + \mu_n ~~~ \NB{\tikz[font=\tiny, scale=1]{\input{sch/foam_digoncup_n0}}} \\ & + s \sum_{k+l=n} ~~~ \NB{\tikz[font=\tiny, scale=1]{\input{sch/foam_digoncup_kl}}} 
         - b \lambda_n ~~~ \NB{\tikz[font=\tiny, scale=1]{\input{sch/foam_digoncup_N}}} - a\mu_n ~~~ \NB{\tikz[font=\tiny, scale=1]{\input{sch/foam_digoncup_n}}} \\
        & -\frac{s}{2} \sum_{k+l=n} \left ( \hspace{2mm} \NB{\tikz[font=\tiny, scale=1]{\input{sch/foam_digoncup_abkl}}} - ~~~\NB{\tikz[font=\tiny, scale=1]{\input{sch/foam_digoncup_akl}}} - ~~~ \NB{\tikz[font=\tiny, scale=1]{\input{sch/foam_digoncup_bkl}}} \right ). 
    \end{align*}
    Let $\deco^a$ denote a decoration on the facet of thickness $a$.
    The computation between the parenthesis can be summarized by
    \[
    \sum_{k+l=n} \deco_k^{a+b}\deco_l^{a+b} - \deco_k^{a}\deco_l^a -\deco_k^b\deco_l^b = 2 \sum_{k+l=n} \deco_k^a\deco_l^b,
    \]
    using the identity \(\deco_i^{a+b} = \deco_i^a + \deco_i^b.\)
    Thus, \[\wgen_n \left ( \hspace{2mm} \NB{\tikz[font=\tiny, scale=1]{\input{sch/foam_digoncup_empty}}} \right ) = 0.\]
\end{proof}
The proofs of the following lemmas are analogous.

\begin{lem}
    \label{rd_cp}
    The foams
    \[\NB{\tikz[font=\tiny, scale=1]{\input{sch/foam_cap_a_empty}}} \hspace{5.5mm} \text{and} \hspace{5.5mm} \NB{\tikz[font=\tiny, scale=1]{\input{sch/foam_cup_a_empty}}}\]
    induce $\Witt$-equivariant morphisms
    \[ \mapcap: \NB{\tikz[font=\tiny, scale=0.7]{\input{sch/web_circle_cap}}} \rightarrow \varnothing \hspace{5.5mm} \text{and} \hspace{5.5mm} \mapcup: \varnothing \rightarrow \NB{\tikz[font=\tiny, scale=0.7]{\input{sch/web_circle_cup}}}.\]
\end{lem}

\begin{lem}
    \label{rd_saddle}
    The foam 
    \[\NB{\tikz[font=\tiny, scale=1]{\input{sch/foam_saddle_a_empty}}} \]
    induces an $\Witt$-equivariant morphism
    \[\mapsaddle: \NB{\tikz[font=\tiny, scale=.8]{\input{sch/web_saddle_a_beg}}} \rightarrow \NB{\tikz[font=\tiny, scale=.8]{\input{sch/web_saddle_a_end}}}.\]
\end{lem}

\begin{lem}
    \label{rd_zip}
    The foams
    \[\NB{\tikz[font=\tiny, scale=1]{\input{sch/foam_zip_a_empty}}} \hspace{5.5mm} \text{and} \hspace{5.5mm} \NB{\tikz[font=\tiny, scale=1]{\input{sch/foam_unzip_a_empty}}}\]
    induce $\Witt$-equivariant morphisms
    \[ \mapzip: \NB{\tikz[font=\tiny, scale=0.7]{\input{sch/web_two_arrows_ab}}} \rightarrow \NB{\tikz[font=\tiny, scale=0.7]{\input{sch/web_zip_end}}} \hspace{5.5mm} \text{and} \hspace{5.5mm} \mapunzip: \NB{\tikz[font=\tiny, scale=0.7]{\input{sch/web_unzip_beg}}} \rightarrow \NB{\tikz[font=\tiny, scale=0.7]{\input{sch/web_two_arrows_ab}}}.\]
\end{lem}


%% file: sch/web_associativity1.tex
\begin{scope}
  \coordinate (b) at ( 0,-0.50);
  \coordinate (m) at (0,0);
  \coordinate (tr) at (+1,  1);
  \coordinate (tm) at (0,  1);
  \coordinate (tl) at (-1,  1);
  \draw[>-, >= to] (b) -- (m) node[pos=0, below] {$a+b+c$};
  \draw[-to] (m) .. controls +(0,0) and + (0, -0.5) .. (tl) node[pos =
  1, above] {$a$} coordinate[pos= 0.5] (ml) coordinate[pos = 0.25] (lml);
  \draw[-to] (m) .. controls +(0,0) and + (0, -0.5) .. (tr) node[pos =
  1, above] {$c$};
  \draw[-to] (ml)  .. controls +(0,0) and + (0, -0.5) .. (tm) node[pos =
  1, above] {$b$};
\end{scope}

%% file: sch/web_associativity2.tex
\begin{scope}
  \coordinate (b) at ( 0,-0.50);
  \coordinate (m) at (0,0);
  \coordinate (tr) at (+1,  1);
  \coordinate (tm) at (0,  1);
  \coordinate (tl) at (-1,  1);
  \draw[>-, >= to] (b) -- (m) node[pos=0, below] {$a+b+c$};
  \draw[-to] (m) .. controls +(0,0) and + (0, -0.5) .. (tl) node[pos =
  1, above] {$c$} coordinate[pos= 0.5] (ml) coordinate[pos = 0.25] (lml);
  \draw[-to] (m) .. controls +(0,0) and + (0, -0.5) .. (tr) node[pos =
  1, above] {$a$};
  \draw[-to] (ml)  .. controls +(0,0) and + (0, -0.5) .. (tm) node[pos =
  1, above] {$b$};
\end{scope}

%% file: sch/web_associativity3.tex
\begin{scope}
  \coordinate (t) at ( 0,0.50);
  \coordinate (m) at (0,0);
  \coordinate (br) at (+1,  -1);
  \coordinate (bm) at ( 0,  -1);
  \coordinate (bl) at (-1,  -1);
  \draw[-to] (m) -- (t) node[pos=1, above] {$a+b+c$};
  \draw[>-, >=to] (bl) .. controls +(0,0.5) and + (0, 0) .. (m) node[pos =
  0, below] {$a$} coordinate[pos= 0.5] (ml) coordinate[pos = 0.75] (lml);
  \draw[>-, >=to] (br) .. controls +(0,0.5) and + (0, -0) .. (m) node[pos =
  0, below] {$c$};
  \draw[>-, >=to] (bm)  .. controls +(0,0.5) and + (0, 0) .. (ml) node[pos =
  0, below] {$b$};
\end{scope}

%% file: sch/foam_digoncap_empty.tex
\begin{scope}[yscale=-1, font=\tiny]
  \begin{scope}
    \coordinate (Lt) at (0,2);
    \coordinate (Rt) at (2,2);
    \coordinate (lt) at (0.5,2);
    \coordinate (rt) at (1.5,2);
    \coordinate (Lb) at (0,1);
    \coordinate (Rb) at (2,1);
  \end{scope}
  \draw (lt) .. controls +(0.3, 0.3) and +(-.3, 0.3) ..  (rt) node
  [pos=0.5, yshift= -0.1cm] {$a$};
  \draw[white!50!black] (lt) .. controls +(0.3, -.3) and +(-.3, -.3) ..  (rt) node
  [pos=0.5, yshift= 0.1cm, black] {$b$};
  \draw (lt) .. controls +(0 , -.9) and +(0, -.9) ..  (rt);
  \draw (rt) -- (Rt) node[pos=1, below] {$a+b$} -- (Rb) -- (Lb) -- (Lt) -- (lt);
\end{scope}

%% file: sch/foam_digoncup_empty.tex
\begin{scope}
  \begin{scope}
    \coordinate (Lt) at (0,2);
    \coordinate (Rt) at (2,2);
    \coordinate (lt) at (0.5,2);
    \coordinate (rt) at (1.5,2);
    \coordinate (Lb) at (0,1);
    \coordinate (Rb) at (2,1);
  \end{scope}
  \draw (lt) .. controls +(0.3, 0.3) and +(-.3, 0.3) ..  (rt) node
  [pos=0.5, yshift= 0.1cm] {$a$};
  \draw (lt) .. controls +(0.3, -.3) and +(-.3, -.3) ..  (rt) node
  [pos=0.5, yshift= -0.1cm] {$b$};
  \draw (lt) .. controls +(0 , -.9) and +(0, -.9) ..  (rt);
  \draw (rt) -- (Rt) node[pos=1, above] {$a+b$} -- (Rb) -- (Lb) -- (Lt) -- (lt);
\end{scope}

%% file: sch/web_digoncap_beg.tex
\begin{scope}
  \coordinate (bm) at ( 0, -1);
  \coordinate (cm) at ( 0, -0.6);
  \coordinate (sm) at ( 0,  0.6);
  \coordinate (tm) at ( 0, 1);
  \draw[->] (sm) -- (tm) node[above, pos =1] {$a+b$};
  \draw[->-] (cm) .. controls + ( 0.6, 0.6) and + ( 0.6, -0.6) .. (sm)
  node [pos = 0.5, right] {$b$} node[pos= 0.8 , draw=\crd, fill=white,regular polygon, regular polygon sides=3,inner sep=0.9pt] {} node[pos=0.8, right]{$a\mu$}; ;
  \draw[->-] (cm) .. controls + (-0.6, 0.6) and + (-0.6, -0.6) .. (sm)
  node [pos = 0.5,  left] {$a$} node[pos= 0.8 , draw=\crd, fill=white,regular polygon, regular polygon sides=3,inner sep=0.9pt] {} node[pos=0.8, left]{$b\lambda$};
  \draw[>-] (bm) -- (cm) node[below, pos =0] {$a+b$} node[pos=1, draw=\crd, fill= white, circle , inner sep=1.5pt]{} node[pos=1, right]{$-\bar{s}$} ;
\end{scope}

%% file: sch/web_l_ab.tex
\begin{scope}
  \coordinate (bm) at ( 0, -1);
  \coordinate (tm) at ( 0, 1);
  \draw[>->] (bm) -- (tm) node[above, pos =1] {$a+b$} node[below, pos =0] {$a+b$};
\end{scope}

%% file: sch/web_digoncup_end.tex
\begin{scope}
  \coordinate (bm) at ( 0, -1);
  \coordinate (cm) at ( 0, -0.6);
  \coordinate (sm) at ( 0,  0.6);
  \coordinate (tm) at ( 0, 1);
  \draw[->] (sm) -- (tm) node[above, pos =1] {$a+b$};
  \draw[->-] (cm) .. controls + ( 0.6, 0.6) and + ( 0.6, -0.6) .. (sm)
  node [pos = 0.5, right] {$b$} node[pos= 0.8 , draw=\crd, fill=white,regular polygon, regular polygon sides=3,inner sep=0.9pt] {} node[pos=0.8, right]{$a\mu$}; ;
  \draw[->-] (cm) .. controls + (-0.6, 0.6) and + (-0.6, -0.6) .. (sm)
  node [pos = 0.5,  left] {$a$} node[pos= 0.8 , draw=\crd, fill=white,regular polygon, regular polygon sides=3,inner sep=0.9pt] {} node[pos=0.8, left]{$b\lambda$};
  \draw[>-] (bm) -- (cm) node[below, pos =0] {$a+b$} node[pos=1, draw=\crd, fill= white, circle , inner sep=1.5pt]{} node[pos=1, right]{$s$} ;
\end{scope}

%% file: sch/foam_digoncup_0n.tex
\begin{scope}
  \begin{scope}
    \coordinate (Lt) at (0,2);
    \coordinate (Rt) at (2,2);
    \coordinate (lt) at (0.5,2);
    \coordinate (rt) at (1.5,2);
    \coordinate (Lb) at (0,1);
    \coordinate (Rb) at (2,1);
  \end{scope}
  \draw (lt) .. controls +(0.3, 0.3) and +(-.3, 0.3) ..  (rt) node[above, pos=0.3] {$a$} node[below, pos =0.6] {$\deco_n$};
  \draw (lt) .. controls +(0.3, -.3) and +(-.3, -.3) ..  (rt)node[below, pos =0.25] {$b$} node[below, pos =0.65] {$\deco_0$};
  \draw (lt) .. controls +(0 , -.9) and +(0, -.9) ..  (rt);
  \draw (rt) -- (Rt) node[pos=1, above] {$a+b$} -- (Rb) -- (Lb) -- (Lt) -- (lt);
\end{scope}

%% file: sch/foam_digoncup_n0.tex
\begin{scope}
  \begin{scope}
    \coordinate (Lt) at (0,2);
    \coordinate (Rt) at (2,2);
    \coordinate (lt) at (0.5,2);
    \coordinate (rt) at (1.5,2);
    \coordinate (Lb) at (0,1);
    \coordinate (Rb) at (2,1);
  \end{scope}
  \draw (lt) .. controls +(0.3, 0.3) and +(-.3, 0.3) ..  (rt) node[above, pos=0.3] {$a$} node[below, pos =0.6] {$\deco_0$};
  \draw (lt) .. controls +(0.3, -.3) and +(-.3, -.3) ..  (rt)node[below, pos =0.25] {$b$} node[below, pos =0.65] {$\deco_n$};
  \draw (lt) .. controls +(0 , -.9) and +(0, -.9) ..  (rt);
  \draw (rt) -- (Rt) node[pos=1, above] {$a+b$} -- (Rb) -- (Lb) -- (Lt) -- (lt);
\end{scope}

%% file: sch/foam_digoncup_kl.tex
\begin{scope}
  \begin{scope}
    \coordinate (Lt) at (0,2);
    \coordinate (Rt) at (2,2);
    \coordinate (lt) at (0.5,2);
    \coordinate (rt) at (1.5,2);
    \coordinate (Lb) at (0,1);
    \coordinate (Rb) at (2,1);
  \end{scope}
  \draw (lt) .. controls +(0.3, 0.3) and +(-.3, 0.3) ..  (rt) node[above, pos=0.3] {$a$} node[below, pos =0.6] {$\deco_k$};
  \draw (lt) .. controls +(0.3, -.3) and +(-.3, -.3) ..  (rt)node[below, pos =0.25] {$b$} node[below, pos =0.65] {$\deco_l$};
  \draw (lt) .. controls +(0 , -.9) and +(0, -.9) ..  (rt);
  \draw (rt) -- (Rt) node[pos=1, above] {$a+b$} -- (Rb) -- (Lb) -- (Lt) -- (lt);
\end{scope}

%% file: sch/foam_digoncup_N.tex
\begin{scope}
  \begin{scope}
    \coordinate (Lt) at (0,2);
    \coordinate (Rt) at (2,2);
    \coordinate (lt) at (0.5,2);
    \coordinate (rt) at (1.5,2);
    \coordinate (Lb) at (0,1);
    \coordinate (Rb) at (2,1);
  \end{scope}
  \draw (lt) .. controls +(0.3, 0.3) and +(-.3, 0.3) ..  (rt) node[above, pos=0.3] {$a$} node[below, pos =0.6] {$\deco_n$};
  \draw (lt) .. controls +(0.3, -.3) and +(-.3, -.3) ..  (rt)node[below, pos =0.25] {$b$} ;
  \draw (lt) .. controls +(0 , -.9) and +(0, -.9) ..  (rt);
  \draw (rt) -- (Rt) node[pos=1, above] {$a+b$} -- (Rb) -- (Lb) -- (Lt) -- (lt);
\end{scope}

%% file: sch/foam_digoncup_n.tex
\begin{scope}
  \begin{scope}
    \coordinate (Lt) at (0,2);
    \coordinate (Rt) at (2,2);
    \coordinate (lt) at (0.5,2);
    \coordinate (rt) at (1.5,2);
    \coordinate (Lb) at (0,1);
    \coordinate (Rb) at (2,1);
  \end{scope}
  \draw (lt) .. controls +(0.3, 0.3) and +(-.3, 0.3) ..  (rt) node[above, pos=0.3] {$a$};
  \draw (lt) .. controls +(0.3, -.3) and +(-.3, -.3) ..  (rt)node[below, pos =0.25] {$b$} node[below, pos =0.65] {$\deco_n$};
  \draw (lt) .. controls +(0 , -.9) and +(0, -.9) ..  (rt);
  \draw (rt) -- (Rt) node[pos=1, above] {$a+b$} -- (Rb) -- (Lb) -- (Lt) -- (lt);
\end{scope}

%% file: sch/foam_digoncup_abkl.tex
\begin{scope}
  \begin{scope}
    \coordinate (Lt) at (0,2);
    \coordinate (Rt) at (2,2);
    \coordinate (lt) at (0.5,2);
    \coordinate (rt) at (1.5,2);
    \coordinate (Lb) at (0,1);
    \coordinate (Rb) at (2,1);
  \end{scope}
  \draw (lt) .. controls +(0.3, 0.3) and +(-.3, 0.3) ..  (rt) node[above, pos=0.3] {$a$};
  \draw (lt) .. controls +(0.3, -.3) and +(-.3, -.3) ..  (rt)node[below, pos =0.25] {$b$};
  \draw (lt) .. controls +(0 , -.9) and +(0, -.9) ..  (rt);
  \draw (rt) -- (Rt) node[pos=1, above] {$a+b$} node[below, pos =0.5] {$\deco_l$} -- (Rb) -- (Lb) -- (Lt) -- (lt) node[below, pos =0.5] {$\deco_k$};
\end{scope}

%% file: sch/foam_digoncup_akl.tex
\begin{scope}
  \begin{scope}
    \coordinate (Lt) at (0,2);
    \coordinate (Rt) at (2,2);
    \coordinate (lt) at (0.5,2);
    \coordinate (rt) at (1.5,2);
    \coordinate (Lb) at (0,1);
    \coordinate (Rb) at (2,1);
  \end{scope}
  \draw (lt) .. controls +(0.3, 0.3) and +(-.3, 0.3) ..  (rt) node[above, pos=0.3] {$a$} node[below, pos =0.5] {$\deco_k \deco_l$};
  \draw (lt) .. controls +(0.3, -.3) and +(-.3, -.3) ..  (rt)node[below, pos =0.25] {$b$};
  \draw (lt) .. controls +(0 , -.9) and +(0, -.9) ..  (rt);
  \draw (rt) -- (Rt) node[pos=1, above] {$a+b$} -- (Rb) -- (Lb) -- (Lt) -- (lt);
\end{scope}

%% file: sch/foam_digoncup_bkl.tex
\begin{scope}
  \begin{scope}
    \coordinate (Lt) at (0,2);
    \coordinate (Rt) at (2,2);
    \coordinate (lt) at (0.5,2);
    \coordinate (rt) at (1.5,2);
    \coordinate (Lb) at (0,1);
    \coordinate (Rb) at (2,1);
  \end{scope}
  \draw (lt) .. controls +(0.3, 0.3) and +(-.3, 0.3) ..  (rt) node[above, pos=0.3] {$a$};
  \draw (lt) .. controls +(0.3, -.3) and +(-.3, -.3) ..  (rt)node[above, pos =0.5] {$b$} node[below, pos =0.5,yshift=.25mm] {$\deco_k\deco_l$};
  \draw (lt) .. controls +(0 , -.9) and +(0, -.9) ..  (rt);
  \draw (rt) -- (Rt) node[pos=1, above] {$a+b$} -- (Rb) -- (Lb) -- (Lt) -- (lt);
\end{scope}

%% file: sch/web_circle_cap.tex
\begin{scope}[scale= 1]
  \draw[->] (0,0) arc (0:360: 0.8) node[pos=1, right] {$a$} coordinate[pos=0.45]
  (ga); 
  \filldraw[draw= \crd, fill = \csrd] (ga) circle (1mm)
  node[left, black] {$-\frac{1}{2}$};
\end{scope}

%% file: sch/web_circle_cup.tex
\begin{scope}[scale= 1]
  \draw[->] (0,0) arc (0:360: 0.8) node[pos=1, right] {$a$} coordinate[pos=0.45]
  (ga); 
  \filldraw[draw= \crd, fill = \csrd] (ga) circle (1mm)
  node[left, black] {$\frac{1}{2}$};
\end{scope}

%% file: sch/web_saddle_a_beg.tex
\begin{scope}[scale= 1]
  \draw [->] (45:0.8) arc (-45:-135:0.8) node [pos=1, above] {$a$};
  \draw [<-] (-45:0.8) arc (45: 135:0.8) node [pos=0, below] {$a$};
\end{scope}

%% file: sch/web_saddle_a_end.tex
\begin{scope}[scale= 1]
    \draw [->] (45:0.8) arc (135:225:0.8) coordinate[pos=0.5] (ga) node [pos=1, below] {$a$};
    \draw [<-] (135:0.8) arc (-135:-225:-0.8) node [pos=0, above] {$a$};
    \filldraw[draw= \crd, fill = \csrd] (ga) circle (1mm) node[right, black] {$-\frac{1}{2}$};
\end{scope}

%% file: sch/foam_zip_a_empty.tex
\begin{scope}
  \begin{scope}
    \coordinate (L1) at (0.2,0.4);
    \coordinate (L2) at (0,0);
    \coordinate (R1) at (2.2,0.4);
    \coordinate (R2) at (2,0);
    \coordinate (ML) at (0.6, 0.2);
    \coordinate (MR) at (1.6, 0.2);
    \draw[->-] (ML) -- (MR) node[above, midway] {$a+b$};
    \draw (MR) .. controls +(0, 0) and +(-0.3,0) .. (R1) ;
    \draw (MR) .. controls +(0, 0) and +(-0.3,0) .. (R2);
    \draw (L1) .. controls +( 0.3, 0) and +(0,0) .. (ML);
    \draw (L2) .. controls +( 0.3, 0) and +(0,0) .. (ML);
  \end{scope}  
 \begin{scope}[yshift = -1cm]
    \coordinate (L1B) at (0.2,0.4);
    \coordinate (L2B) at (0,0);
    \coordinate (R1B) at (2.2,0.4);
    \coordinate (R2B) at (2,0);
    \draw[->-] (L1B) .. controls +( 0, 0) and +(0,0) .. (R1B) node [right, pos
    = 1] {$a$};
    \draw[->-] (L2B) .. controls +( 0, 0) and +(0,0) .. (R2B) node [right, pos
    = 1] {$b$};
 \end{scope}  
  \draw (R1) -- (R1B);
  \draw (R2) -- (R2B);
  \draw (L1) -- (L1B);
  \draw (L2) -- (L2B);
  \draw (ML) .. controls +(0, -0.6) and +(0, -0.6) .. (MR);
\end{scope}

%% file: sch/foam_unzip_a_empty.tex
\begin{scope}
  \begin{scope}
    \coordinate (L1) at (0.2,0.4);
    \coordinate (L2) at (0,0);
    \coordinate (R1) at (2.2,0.4);
    \coordinate (R2) at (2,0);
    \coordinate (ML) at (0.6, 0.2);
    \coordinate (MR) at (1.6, 0.2);
    \draw[->-] (ML) -- (MR) node[below, midway] {$a+b$};
    \draw (MR) .. controls +(0, 0) and +(-0.3,0) .. (R1) ;
    \draw (MR) .. controls +(0, 0) and +(-0.3,0) .. (R2);
    \draw (L1) .. controls +( 0.3, 0) and +(0,0) .. (ML);
    \draw (L2) .. controls +( 0.3, 0) and +(0,0) .. (ML);
  \end{scope}  
 \begin{scope}[yshift = 1cm]
    \coordinate (L1B) at (0.2,0.4);
    \coordinate (L2B) at (0,0);
    \coordinate (R1B) at (2.2,0.4);
    \coordinate (R2B) at (2,0);
    \draw[->-] (L1B) .. controls +( 0, 0) and +(0,0) .. (R1B) node [left, pos
    = 0] {$a$};
    \draw[->-] (L2B) .. controls +( 0, 0) and +(0,0) .. (R2B) node [left, pos
    = 0] {$b$};
 \end{scope}  
  \draw (R1) -- (R1B);
  \draw (R2) -- (R2B);
  \draw (L1) -- (L1B);
  \draw (L2) -- (L2B);
  \draw (ML) .. controls +(0, 0.6) and +(0, 0.6) .. (MR);
\end{scope}

%% file: sch/web_two_arrows_ab.tex
\begin{scope}
    \coordinate (bl) at (-0.5, -1);
    \coordinate (br) at ( 0.5, -1);
    \coordinate (tl) at (-0.5,  1);
    \coordinate (tr) at ( 0.5,  1);
    \draw[->] (bl) -- (tl) node[pos=0,below] {$a$};
    \draw[->] (br) -- (tr) node[pos=0,below] {$b$};
\end{scope}

%% file: SF_link_homology.tex
Let \(D\) be a diagram of a link \(L\) in \(\R^3\).

\begin{dfn}
The complex \(\CKRW(L)\) is defined using the cube of resolutions construction associated to the braiding complexes
\begin{gather}
  \NB{\tikz[xscale = 0.6]{\input{sch/positive_crossing}}} :=
\NB{\tikz[xscale = 3.5, yscale = 3]
{\node (i0) at (0.2, 0){$\NB{\tikz[font= \tiny,scale=0.6]{\input{sch/reidem_lemm_14}}}$};
\node (i1) at (-1, 0) {$q\NB{\tikz[font= \tiny,scale=0.6]{\input{sch/def_homo_cross_pos}}} $};
\draw[->] (i1) -- (i0) coordinate[pos=0.45] (a);
\node[above] at (a) {\NB{\tikz[font=\tiny, scale=.5]{\input{sch/foam_unzip_empty}}}};
  }},\\
   \NB{\tikz[xscale =- 0.6]{\input{sch/positive_crossing}}}:=
    \NB{ \tikz[xscale = 3.5, yscale = 3]{
    \node (i0) at (-1, 0) {$ \, \NB{\tikz[font= \tiny, scale=0.6]{\input{sch/reidem_lemm_14}}}$};
    \node (i1) at (0.2, 0) {$q^{-1}\NB{\tikz[font= \tiny, scale=0.6]{\input{sch/def_homo_cross_neg}}}$ };
    \draw[->] (i0) -- (i1) coordinate[pos=0.6] (b);
    \node[above] at (b) {\NB{\tikz[font=\tiny, scale=.5]{\input{sch/foam_zip_empty}}}}; }} .
\end{gather}
twisted with red dots to ensure that the differentials are \(\Witt\)-equivariant.
The underlined terms are in cohomological degree 0.
The cohomology \( \KRW_N(L,\scalars):=H^*(\CKRW(L,R))\) is equipped with the action of \(\Witt\).
The dependence on the parameters \(\lambda, \mu,\) and \(s\) is suppressed.
\end{dfn}

\begin{thm}
The chain complex $\CKRW(-)$ is invariant under framed isotopy in the relative homotopy category.
\end{thm}
The rest of the section is dedicated to proving the following theorem by showing invariance under the Reidemeister moves.

%% file: sch/positive_crossing.tex
\begin{scope}[font=\tiny]
  \draw[->] (0.5, -0.5) ..controls +(0,0.3) and +(0,-0.3) .. (-0.5, 0.5);
  \fill[white] (0,0) circle (2mm);
  \draw[->] (-0.5, -0.5) ..controls +(0,0.3) and +(0,-0.3) .. (0.5,0.5);
\end{scope}

%% file: sch/reidem_lemm_14.tex
\begin{scope}
    \coordinate (bl) at (-0.5, -1);
    \coordinate (br) at ( 0.5, -1);
    \coordinate (tl) at (-0.5,  1);
    \coordinate (tr) at ( 0.5,  1);
    \draw[->] (bl) -- (tl);
    \draw[->] (br) -- (tr);
\end{scope}

%% file: sch/foam_unzip_empty.tex
\begin{scope}
  \begin{scope}
    \coordinate (L1) at (0.2,0.4);
    \coordinate (L2) at (0,0);
    \coordinate (R1) at (2.2,0.4);
    \coordinate (R2) at (2,0);
    \coordinate (ML) at (0.6, 0.2);
    \coordinate (MR) at (1.6, 0.2);
    \draw[double] (ML) -- (MR) node[below, midway] {}; 
    \draw (MR) .. controls +(0, 0) and +(-0.3,0) .. (R1) ;
    \draw (MR) .. controls +(0, 0) and +(-0.3,0) .. (R2);
    \draw (L1) .. controls +( 0.3, 0) and +(0,0) .. (ML);
    \draw (L2) .. controls +( 0.3, 0) and +(0,0) .. (ML);
  \end{scope}  
 \begin{scope}[yshift = 1cm]
    \coordinate (L1B) at (0.2,0.4);
    \coordinate (L2B) at (0,0);
    \coordinate (R1B) at (2.2,0.4);
    \coordinate (R2B) at (2,0);
    \draw[->-] (L1B) .. controls +( 0, 0) and +(0,0) .. (R1B) node [left, pos
    = 0] {};
    \draw[->-] (L2B) .. controls +( 0, 0) and +(0,0) .. (R2B) node [left, pos
    = 0] {};
 \end{scope}  
  \draw (R1) -- (R1B);
  \draw (R2) -- (R2B);
  \draw (L1) -- (L1B);
  \draw (L2) -- (L2B);
  \draw (ML) .. controls +(0, 0.6) and +(0, 0.6) .. (MR);
\end{scope}

%% file: sch/foam_zip_empty.tex
\begin{scope}
  \begin{scope}
    \coordinate (L1) at (0.2,0.4);
    \coordinate (L2) at (0,0);
    \coordinate (R1) at (2.2,0.4);
    \coordinate (R2) at (2,0);
    \coordinate (ML) at (0.6, 0.2);
    \coordinate (MR) at (1.6, 0.2);
    \draw[double] (ML) -- (MR) node[above, midway] {};
    \draw (MR) .. controls +(0, 0) and +(-0.3,0) .. (R1) ;
    \draw (MR) .. controls +(0, 0) and +(-0.3,0) .. (R2);
    \draw (L1) .. controls +( 0.3, 0) and +(0,0) .. (ML);
    \draw (L2) .. controls +( 0.3, 0) and +(0,0) .. (ML);
  \end{scope}  
 \begin{scope}[yshift = -1cm]
    \coordinate (L1B) at (0.2,0.4);
    \coordinate (L2B) at (0,0);
    \coordinate (R1B) at (2.2,0.4);
    \coordinate (R2B) at (2,0);
    \draw[->-] (L1B) .. controls +( 0, 0) and +(0,0) .. (R1B) node [right, pos
    = 1] {};
    \draw[->-] (L2B) .. controls +( 0, 0) and +(0,0) .. (R2B) node [right, pos
    = 1] {};
 \end{scope}  
  \draw (R1) -- (R1B);
  \draw (R2) -- (R2B);
  \draw (L1) -- (L1B);
  \draw (L2) -- (L2B);
  \draw (ML) .. controls +(0, -0.6) and +(0, -0.6) .. (MR);
\end{scope}

%% file: SF_Lemma.tex
\label{sliding_sect}

We begin with a technical lemma which is essential in the remainder of the section.
The proof below is a clarification and adaptation of
similar lemmas that appear in \cite{KR16} and \cite{QRSW23}.

\begin{lem}
    \label{lemm_cross}
     For any $\omega \in \scalars$, there is an isomorphism in the relative homotopy category between the complexes
\begin{align}
{}_xC &:= \NB{\tikz[font=\tiny]{\input{sch/reidem_lemm_9}}} = 
\NB{
  \begin{tikzpicture}[xscale =4, yscale=2]
    \node (BT) at (0,3) {$\NB{\tikz[font=\tiny,scale=0.6]{\input{sch/reidem_lemm_11}}}$};
    \node (AT) at (1.1,3) {$\NB{\tikz[font=\tiny,scale=0.6]{\input{sch/reidem_lemm_13}}}$};
    \draw[-to] (BT) -- (AT) node[pos=0.5, above, scale = 0.7] {$\mapunzip$}; 
    \end{tikzpicture}}\\
C^y &:=\NB{\tikz[font=\tiny]{\input{sch/reidem_lemm_10}}} = 
\NB{
  \begin{tikzpicture}[xscale =4, yscale=2]
    \node (BT) at (0,3) {$\NB{\tikz[font=\tiny,scale=0.6]{\input{sch/reidem_lemm_15}}}$};
    \node (AT) at (1.1,3) {$\NB{\tikz[font=\tiny,scale=0.6]{\input{sch/reidem_lemm_16}}}$};
    \draw[-to] (BT) -- (AT) node[pos=0.5, above, scale = 0.7] {$\mapunzip$}; 
    \end{tikzpicture}}.
\end{align}
\end{lem}

\begin{proof}
Let the left and right facets of \(\mapunzip\) be called the \(x\) and \(y\) facets respectively.
The map \(\mapzip\) (see Lemma (\ref{rd_zip})) is a chain homotopy between the decoration maps \(x\) and \(y\), but it cannot be twisted into equivariance.
Our solution is to transfer twists from the \(x\) facet to the \(y\) facet via an intermediate formal variable \(z\).
To facilitate this, extend the base ring from \(\scalars_N\) to \(\scalars_N[z] := \scalars_{N,1}\),
and the state space \(S\) of any web to \(S[z] := S \otimes_{\scalars_N} \scalars_{N,1}\).
The complexes \({}_xC[z]\) and \(C^y[z]\) are defined as termwise extensions of \({}_xC\) and \(C^y\).

Let
$S[z]^{\gamma}$ denote the representation twisted by
$\wgen_n \mapsto -\gamma(n+1)z^n$,
analogous to the red triangle,
and
let $S^{\star_x^z}[z]$ denote the one twisted by
$\wgen_n \mapsto -h_n(x,z)$,
analogous to the solid red dot.

The variable \(z\) can be imagined as belonging to a cone point that every facet is connected to.
This idea motivates the following algebraic construction.
The map \(C^{\star_x^z}[z]^{\omega} \xrightarrow{x-z} C[z]^{\omega}\) given by the difference of the \(x\) decoration map and multiplication by \(z\) is equivariant and its cokernel is equivariantly isomorphic to the original complex \({}_xC\).
That the \(\star\) twist makes the map \(x-z\) equivariant while the \(\omega\) twist on \(z\) corresponds to the \(\omega\) twist in the cokernel.
Similarly, \(C^y \cong \Coker({C^y}^{\star_y^z}[z]^{\omega} \xrightarrow{y-z} {C^y}[z]^{\omega})\).
In the homotopy category, cokernels are isomorphic to cones and the extra room allows for additional twists on \(\Cone(y-z)\) that make the isomorphism with \(\Cone(x-z)\) corresponding to the chain homotopy \(\mapzip\) equivariant.
This will prove that ${}_xC \cong C^y$.
\begin{equation}
\NB{
  \begin{tikzpicture}[xscale =4, yscale=2]
    \node (CT) at (1.65,3) {\(\Coker(x-z)\)};
    \node (Ct) at (1.65,1.9) {\(\Cone(x-z)\)};
    \node (Cb) at (1.65,0.8) {\(\Cone(\Id)\)};
    \node (zeroT) at (1.65,3.5) {$0$};
    \node (zerob) at (1.65,0.3) {$0$};
    \draw[densely dotted] (CT) -- (1.34, 3);
    \draw[densely dotted] (Ct) -- (1.34, 1.9);
    \draw[densely dotted] (Cb) -- (1.34, 0.8);
    \draw[densely dotted, rounded corners] (-0.22, 3.4) rectangle (1.34, 2.59);
    \draw[densely dotted, rounded corners] (-1.45, 2.3) rectangle (1.34, 1.49);
    \draw[densely dotted, rounded corners] (-1.45, 1.2) rectangle (1.34,0.35);
    \draw [-to] (CT) -- (zeroT);
    \draw [-to] (Ct) -- (CT);
    \draw [-to] (Cb) -- (Ct);
    \draw [-to] (zerob) -- (Cb);
    \node (BT) at (0,3) {$\NB{\tikz[font=\tiny,scale=0.6]{\input{sch/reidem_lemm_11}}}$};
    \node (AT) at (1.1,3) {$\NB{\tikz[font=\tiny,scale=0.6]{\input{sch/reidem_lemm_13}}}$};
    \node (Zt) at (-1.1,1.9) {$\NB{\tikz[font=\tiny,scale=0.6]{\input{sch/reidem_lemm_12}}}^{\star_x^z}[z]^{\omega}$};
    \node (Bt) at (0,1.9) {$\NB{\tikz[font=\tiny,scale=0.6]{\input{sch/reidem_lemm_12}}}[z]^{\omega}\oplus \NB{\tikz[font=\tiny,scale=0.6]{\input{sch/reidem_lemm_14}}}^{\star_x^z}[z]^{\omega}$};
    \node (At) at (1.1,1.9) {$\NB{\tikz[font=\tiny,scale=0.6]{\input{sch/reidem_lemm_14}}}[z]^{\omega}$};
    \node (Zb) at (-1.1,0.8) {$\NB{\tikz[font=\tiny,scale=0.6]{\input{sch/reidem_lemm_12}}}^{\star_x^z}[z]^{\omega}$};
    \node (Bb) at (0,0.8) {$\NB{\tikz[font=\tiny,scale=0.6]{\input{sch/reidem_lemm_12}}}^{\star_x^z}[z]^{\omega}\oplus \NB{\tikz[font=\tiny,scale=0.6]{\input{sch/reidem_lemm_14}}}^{\star_x^z}[z]^{\omega}$};
    \node (Ab) at (1.1,0.8) {$\NB{\tikz[font=\tiny,scale=0.6]{\input{sch/reidem_lemm_14}}}^{\star_x^z}[z]^{\omega}$};
    \draw[-to] (BT) -- (AT) node[pos=0.5, above, scale = 0.7] {$\mapunzip$};
    \draw[-to] (Bt) -- (At) node[pos=0.5, above, scale = 0.7] {$\left(\mapunzip \,\, x -z\right)$};
    \draw[-to] (Bb) -- (Ab) node[pos=0.5, above, scale = 0.7] {$\left(\mapunzip \,\, 1 \right)$};
    \draw[-to] (Zt) -- (Bt) node[pos=0.5, above, scale = 0.7]
    {$\left( \begin{array}{c}z-x \\ \mapunzip\end{array}\right)$};
    \draw[-to] (Zb) -- (Bb) node[pos=0.5, above, scale = 0.7]
    {$\left( \begin{array}{c} 1 \\ -\mapunzip\end{array}\right)$};
    \draw[-to] (Zb) -- (Zt) node[pos=0.5, left, scale=0.7] {$-\Id$};
    \draw[-to] (Ab) -- (At) node[pos=0.5, left, scale=0.7] {$x-z$};
    \draw[-to] (Bb) -- (Bt) node[pos=0.5, left, scale=0.6]
    {$\begin{pmatrix}  x-z & 0 \\ 0 & 1 \end{pmatrix}$};
    \draw[-to] (Bt) -- (BT) node[pos=0.5, left, scale=0.7]
    {$\left(z\mapsto x \,\,\, 0 \right)$};
    \draw[-to] (At) -- (AT) node[pos=0.5, left, scale=0.7]
    {$z\mapsto x$};
    \end{tikzpicture}}
\ ,
\end{equation}

\begin{equation}
\NB{
  \begin{tikzpicture}[xscale =4, yscale=2]
    \node (CT) at (1.65,3)   {\(\Coker(y-z)\)};
    \node (Ct) at (1.65,1.9) {\(\Cone^\prime(y-z)\)};
    \node (Cb) at (1.65,0.8) {\(\Cone^\prime(\Id)\)};
    \node (zeroT) at (1.65,3.5) {$0$};
    \node (zerob) at (1.65,0.3) {$0$};
    \draw[densely dotted] (CT) -- (1.34, 3);
    \draw[densely dotted] (Ct) -- (1.34, 1.9);
    \draw[densely dotted] (Cb) -- (1.34, 0.8);
    \draw[densely dotted, rounded corners] (-0.22, 3.4) rectangle (1.34, 2.59);
    \draw[densely dotted, rounded corners] (-1.45, 2.3) rectangle (1.34, 1.49);
    \draw[densely dotted, rounded corners] (-1.45, 1.2) rectangle (1.34,0.35);
    \draw [-to] (CT) -- (zeroT);
    \draw [-to] (Ct) -- (CT);
    \draw [-to] (Cb) -- (Ct);
    \draw [-to] (zerob) -- (Cb);
    \node (BT) at (0,3) {$\NB{\tikz[font=\tiny,scale=0.6]{\input{sch/reidem_lemm_15}}}$};
    \node (AT) at (1.1,3) {$\NB{\tikz[font=\tiny,scale=0.6]{\input{sch/reidem_lemm_16}}}$};
    \node (Zt) at (-1.1,1.9) {$\NB{\tikz[font=\tiny,scale=0.6]{\input{sch/reidem_lemm_12}}}^{\star_x^z}[z]^{\omega}$};
    \node (Bt) at (0,1.9) {$\NB{\tikz[font=\tiny,scale=0.6]{\input{sch/reidem_lemm_12}}}[z]^{\omega}\oplus \NB{\tikz[font=\tiny,scale=0.6]{\input{sch/reidem_lemm_14}}}^{\star_x^z}[z]^{\omega}$};
    \node (At) at (1.1,1.9) {$\NB{\tikz[font=\tiny,scale=0.6]{\input{sch/reidem_lemm_14}}}[z]^{\omega}$};
    \node (Zb) at (-1.1,0.8) {$\NB{\tikz[font=\tiny,scale=0.6]{\input{sch/reidem_lemm_12}}}^{\star_x^z}[z]^{\omega}$};
    \node (Bb) at (0,0.8) {$\NB{\tikz[font=\tiny,scale=0.6]{\input{sch/reidem_lemm_12}}}^{\star_y^z}[z]^{\omega}\oplus \NB{\tikz[font=\tiny,scale=0.6]{\input{sch/reidem_lemm_14}}}^{\star_x^z}[z]^{\omega}$};
    \node (Ab) at (1.1,0.8) {$\NB{\tikz[font=\tiny,scale=0.6]{\input{sch/reidem_lemm_14}}}^{\star_y^z}[z]^{\omega}$};
    \draw[-to] (BT) -- (AT) node[pos=0.5, above, scale = 0.7] {$\mapunzip$};
    \draw[-to] (Bt) -- (At) node[pos=0.5, above, scale = 0.7] {$\left(\mapunzip \,\, y -z\right)$};
    \draw[-to] (Bb) -- (Ab) node[pos=0.5, above, scale = 0.7] {$\left(\mapunzip \,\, 1 \right)$};
    \draw[-to] (Zt) -- (Bt) node[pos=0.5, above, scale = 0.7]
    {$\left( \begin{array}{c}z-y \\ \mapunzip\end{array}\right)$};
    \draw[-to] (Zb) -- (Bb) node[pos=0.5, above, scale = 0.7]
    {$\left( \begin{array}{c} 1 \\ -\mapunzip\end{array}\right)$};
    \draw[-to] (Zb) -- (Zt) node[pos=0.5, left, scale=0.7] {$-\Id$};
    \draw[-to] (Ab) -- (At) node[pos=0.5, left, scale=0.7] {$y-z$};
    \draw[-to] (Bb) -- (Bt) node[pos=0.5, left, scale=0.6]
    {$\begin{pmatrix}  y-z & 0 \\ 0 & 1 \end{pmatrix}$};
    \draw[-to] (Bt) -- (BT) node[pos=0.5, left, scale=0.7]
    {$\left(z\mapsto y \,\, 0 \right)$};
    \draw[-to] (At) -- (AT) node[pos=0.5, left, scale=0.7]
    {$z\mapsto y$};
    \draw[-to, \crd] ($(Bt)+(0.04,0.15)$) .. controls +(-0.1,0.1) and +(0.1, 0.1) .. +(-0.2, 0) node[scale=0.7,\crd, pos=0.8, above] {(1)};
    \draw[-to, \crd] ($(Bb)+(0.08, -0.22)$) .. controls +(-0.1,-0.1) and +(0.1, -0.1) .. +(-0.25, 0) node[scale=0.7,\crd, pos=0.2, below] {(2)};
    \end{tikzpicture}}
\ .
\end{equation}

The red arrows indicate additional twists on the direct sum objects defined by
\begin{align}
&\textcolor{\crd}{(1)} : 
L_n \mapsto
(h_n(x,z) - h_n(x,y))
\left(
\begin{array}{cc}
    0 & 0 \\
    0 & \mapzip
\end{array}
\right)\\
&\textcolor{\crd}{(2)} :
L_n \mapsto
\frac{h_n(y,z) - h_n(x,z)}{x-y}
\left(
\begin{array}{cc}
    0 & 0 \\
    0 & \mapzip
\end{array}
\right)
\end{align}

The additional twists on \(\Cone^\prime(y-z)\) make the following isomorphism to ${\Cone(x-y)}$ equivariant:
\begin{equation}
\NB{
  \begin{tikzpicture}[xscale =4, yscale=2]
    \node (Zb) at (-1.1,0.8) {$\NB{\tikz[font=\tiny,scale=0.6]{\input{sch/reidem_lemm_12}}}^{\star_x^z}[z]^{\omega}$};
    \node (Bb) at (0,0.8) {$\NB{\tikz[font=\tiny,scale=0.6]{\input{sch/reidem_lemm_12}}}[z]^{\omega}\oplus \NB{\tikz[font=\tiny,scale=0.6]{\input{sch/reidem_lemm_14}}}^{\star_x^z}[z]^{\omega}$};
    \node (Ab) at (1.1,0.8) {$\NB{\tikz[font=\tiny,scale=0.6]{\input{sch/reidem_lemm_14}}}[z]^{\omega}$};
    \node (Zt) at (-1.1,1.9) {$\NB{\tikz[font=\tiny,scale=0.6]{\input{sch/reidem_lemm_12}}}^{\star_x^z}[z]^{\omega}$};
    \node (Bt) at (0,1.9) {$\NB{\tikz[font=\tiny,scale=0.6]{\input{sch/reidem_lemm_12}}}[z]^{\omega}\oplus \NB{\tikz[font=\tiny,scale=0.6]{\input{sch/reidem_lemm_14}}}^{\star_x^z}[z]^{\omega}$};
    \node (At) at (1.1,1.9) {$\NB{\tikz[font=\tiny,scale=0.6]{\input{sch/reidem_lemm_14}}}[z]^{\omega}$};
    \draw[-to] (Bt) -- (At) node[pos=0.5, above, scale = 0.7] {$\left(\mapunzip \,\, x-z\right)$};
    \draw[-to] (Bb) -- (Ab) node[pos=0.5, above, scale = 0.7] {$\left(\mapunzip \,\, y-z \right)$};
    \draw[-to] (Zt) -- (Bt) node[pos=0.5, above, scale = 0.7]
    {$\left( \begin{array}{c}z-x \\ \mapunzip\end{array}\right)$};
    \draw[-to] (Zb) -- (Bb) node[pos=0.5, above, scale = 0.7]
    {$\left( \begin{array}{c} z-y \\ \mapunzip\end{array}\right)$};
    \draw[-to] (Zb) -- (Zt) node[pos=0.5, left, scale=0.7] {$\Id$};
    \draw[-to] (Ab) -- (At) node[pos=0.5, left, scale=0.7] {$\Id$};
    \draw[-to] (Bb) -- (Bt) node[pos=0.65, left, scale=0.6]
    {$\begin{pmatrix}  1 & \mapzip \\ 0 & 1 \end{pmatrix}$};
     \draw[-to, \crd] ($(Bb)+(0.04,-0.17)$) .. controls
     +(-0.1,-0.1) and +(0.1, -0.1) .. +(-0.2, 0)
     node[scale=0.7,\crd, pos=0.2, below] {$(1)$};
    \end{tikzpicture}}
\ ,
\end{equation}
\end{proof}

\begin{rmk}
The twisted action of $\Witt$ on an element \(\foam = \foam_1 + \foam_2\) in the state space of the direct sum
\begin{equation}
    M=
    \NB{\begin{tikzpicture}[xscale =4, yscale=2]
        \node (Bb) at (0,0.8) {$\NB{\tikz[font=\tiny,scale=0.6]{\input{sch/reidem_lemm_12}}}^{\star_y^z}[z]^{\omega}\oplus \NB{\tikz[font=\tiny,scale=0.6]{\input{sch/reidem_lemm_14}}}^{\star_x^z}[z]^{\omega}$};
        \draw[-to, \crd] ($(Bb)+(0.08, -0.22)$) .. controls +(-0.1,-0.1) and +(0.1, -0.1) .. +(-0.25, 0) node[scale=0.7,\crd, pos=0.2, below] {(2)};
    \end{tikzpicture}}
\end{equation}
is given explicitly by
\begin{align*}
    \wgen_n \cdot_{M} \foam =  \wgen_n \cdot \foam - & (\lambda(n+1) x^n+ \mu(n+1) y^n) \foam_1 - (n+1)sx^n \foam_1 - h(z,y)\foam_1 \\
    &-h_n(z,x) \foam_2 + \frac{h_n(y,z)-h_n(x,z)}{x-y} \mapzip \circ \foam_2. 
\end{align*}
\end{rmk}

Let \(\omega, \epsilon, \lambda,\) and \(\mu\) be elements of \(\scalars\).
The following lemmas are consequences of dot migration.

\begin{lem}
    \label{lemm_cross2}
    There is an isomorphism in the relative homotopy category: 
    \begin{align}
     \NB{\tikz[font=\tiny]{\input{sch/reidem_lemm_19}}} \cong \NB{\tikz[font=\tiny]{\input{sch/reidem_lemm_20}}}.
    \end{align}
\end{lem}

\begin{cor}
    \label{fund_lemm}
    There are isomorphisms in the relative homotopy category: 
    \begin{align} \label{lemreid1}
        \NB{\tikz[font=\tiny]{\input{sch/reidem_lemm_1}}} \cong \NB{\tikz[font=\tiny]{\input{sch/reidem_lemm_2}}} \qquad \text{and} \qquad \NB{\tikz[font=\tiny]{\input{sch/reidem_lemm_3}}} \cong \NB{\tikz[font=\tiny]{\input{sch/reidem_lemm_4}}}, 
    \end{align}    
\end{cor}

\begin{proof}

The first isomorphism of (\ref{lemreid1}) is deduced from Lemmas \ref{lemm_cross}, \ref{lemm_cross2}, and the properties of red dots:
\begin{align}
    \NB{\tikz[font=\tiny]{\input{sch/reidem_lemm_1}}} \cong \NB{\tikz[font=\tiny]{\input{sch/reidem_lemm_23}}} \cong \NB{\tikz[font=\tiny]{\input{sch/reidem_lemm_24}}} \cong 
    \NB{\tikz[font=\tiny]{\input{sch/reidem_lemm_2}}}
\end{align}

To prove the second isomorphism in (\ref{lemreid1}) apply Proposition \ref{ReidemeisterII}, which does not rely on it.
\begin{align}
    \NB{\tikz[font=\tiny]{\input{sch/reidem_lemm_3}}} \overset{\text{RII}}{\cong} \NB{\tikz[scale=.8,font=\tiny]{\input{sch/reidem_lemm_25}}} \cong \NB{\tikz[scale=.8,font=\tiny]{\input{sch/reidem_lemm_26}}} \overset{\text{RII}}{\cong}
    \NB{\tikz[font=\tiny]{\input{sch/reidem_lemm_4}}}
\end{align}
\end{proof}

\begin{cor}
    There is an isomorphism in the relative homotopy category:
    \[\NB{\tikz[font=\tiny]{\input{sch/reidem_lemm_17}}} \cong \NB{\tikz[font=\tiny]{\input{sch/reidem_lemm_18}}}.\]
\end{cor}

\begin{lem}
    \label{migr_srd}
     Let $\omega$ and $\epsilon$ be elements of $\scalars$, there are isomorphisms in the relative homotopy category:
        \begin{align} \label{lemreid2}
        \NB{\tikz[font=\tiny]{\input{sch/reidem_lemm_5}}} \cong \NB{\tikz[font=\tiny]{\input{sch/reidem_lemm_6}}} \qquad \text{and} \qquad \NB{\tikz[font=\tiny]{\input{sch/reidem_lemm_7}}} \cong \NB{\tikz[font=\tiny]{\input{sch/reidem_lemm_8}}}.
    \end{align}
\end{lem}

\begin{proof}
It suffices to prove the relation
\begin{equation}
\label{lemm_solid}
\NB{\tikz[font=\tiny]{\input{sch/reidem_lemm_27}}} \cong \NB{\tikz[font=\tiny]{\input{sch/reidem_lemm_28}}}.
\end{equation}
The rest of the proof follows as in Corollary \ref{fund_lemm}.

To avoid dealing with the $N-1$ variables in this twist, introduce a new twist $L_n \mapsto -\sum_{k=0}^n p^{(e)}_k P_{n-k}$ recorded diagrammatically by a black dot.
Since $P_{n-k} = p^{(e)}_{n-k} + \hat{p}_ {n-k}^{(e)}$,
\[\NB{\tikz[font=\small]{\input{sch/relation_hollow_black_dots_2}}} =\NB{\tikz[font=\small]{\input{sch/relation_hollow_black_dots_1}}}.\]
Lemma \ref{lemm_cross} implies that
(\ref{lemm_solid}) follows from
\begin{equation}
\label{lemm_black}
\NB{\tikz[font=\tiny]{\input{sch/reidem_lemm_29}}} \cong \NB{\tikz[font=\tiny]{\input{sch/reidem_lemm_30}}}.
\end{equation}
The proof of (\ref{lemm_black}) is identical to the proof of Lemma \ref{lemm_cross} with
\(S[z]^{\gamma}\) redefined using the twist
\(\wgen_n \mapsto -\gamma\left (\sum_{j=0}^n z^jP_{n-j} \right ) z^k.\)
\end{proof}

%% file: sch/reidem_lemm_9.tex
\begin{scope}[font=\tiny]
  \draw[->] (0.5, -0.5) ..controls +(0,0.3) and +(0,-0.3) .. (-0.5, 0.5) node[pos=1, above] {};
  \fill[white] (0,0) circle (2mm);
  \draw[->] (-0.5, -0.5) ..controls +(0,0.3) and +(0,-0.3) .. (0.5, 0.5) node[pos=0.2, draw=\crd, fill= white, circle , inner sep=2pt]{} node[pos=0.2, left]{$\omega$};
\end{scope}

%% file: sch/reidem_lemm_13.tex
\begin{scope}
    \coordinate (bl) at (-0.5, -1);
    \coordinate (br) at ( 0.5, -1);
    \coordinate (tl) at (-0.5,  1);
    \coordinate (tr) at ( 0.5,  1);
    \draw[->] (bl) -- (tl) node[pos=0.2, draw=\crd, fill= white, circle , inner sep=1.5pt]{} node[pos=0.2, left]{$\omega$};
    \draw[->] (br) -- (tr);
\end{scope}

%% file: sch/reidem_lemm_10.tex
\begin{scope}[font=\tiny]
  \draw[->] (0.5, -0.5) ..controls +(0,0.3) and +(0,-0.3) .. (-0.5,0.5) node[pos=1, above] {};
  \fill[white] (0,0) circle (2mm);
  \draw[->] (-0.5, -0.5) ..controls +(0,0.3) and +(0,-0.3) .. (0.5, 0.5) node[pos=0.75, draw=\crd, fill= white, circle , inner sep=2pt]{} node[pos=0.75, right]{$\omega$};
\end{scope}

%% file: sch/reidem_lemm_15.tex
\begin{scope}
  \coordinate (bl) at (-0.5, -1);
  \coordinate (br) at ( 0.5, -1);
  \coordinate (bm) at (  0,-0.3);
  \coordinate (tl) at (-0.5,  1);
  \coordinate (tr) at ( 0.5,  1);
  \coordinate (tm) at (  0, 0.3);
  \draw[>-]  (bl) .. controls +( 0, 0.5) and +(0,0) .. (bm);
  \draw[>-]  (br) .. controls +( 0, 0.5) and +(0,0) .. (bm);
  \draw[<-]  (tl) .. controls +( 0, -0.5) and +(0,0) .. (tm) coordinate[pos = 0.25] node[pos= 0.35 , draw=\crd, fill=white,regular polygon, regular polygon sides=3,inner sep=0.9pt] {} node [pos=0.3, left]{$\lambda$};
  \draw[<-]  (tr) .. controls +( 0, -0.5) and +(0,0) .. (tm)
  coordinate[pos = 0.25] node[pos= 0.25 , draw=\crd, fill=white,regular polygon, regular polygon sides=3,inner sep=0.9pt] {} node[pos=0.2, right]{$\mu$} 
  node[pos=0.6, draw=\crd, fill= white, circle , inner sep=1.5pt]{} node[pos=0.6, right]{$\omega$};
  \draw [double] (bm) -- (tm) node[pos=0, draw=\crd, fill= white, circle , inner sep=1.5pt]{} node[pos=0, right]{$s$};
\end{scope}

%% file: sch/reidem_lemm_16.tex
\begin{scope}
    \coordinate (bl) at (-0.5, -1);
    \coordinate (br) at ( 0.5, -1);
    \coordinate (tl) at (-0.5,  1);
    \coordinate (tr) at ( 0.5,  1);
    \draw[->] (bl) -- (tl) ;
    \draw[->] (br) -- (tr) node[pos=0.7, draw=\crd, fill= white, circle , inner sep=1.5pt]{} node[pos=0.7, left]{$\omega$};
\end{scope}

%% file: sch/reidem_lemm_12.tex
\begin{scope}
  \coordinate (bl) at (-0.5, -1);
  \coordinate (br) at ( 0.5, -1);
  \coordinate (bm) at (  0,-0.3);
  \coordinate (tl) at (-0.5,  1);
  \coordinate (tr) at ( 0.5,  1);
  \coordinate (tm) at (  0, 0.3);
  \draw[>-]  (bl) .. controls +( 0, 0.5) and +(0,0) .. (bm);
  \draw[>-]  (br) .. controls +( 0, 0.5) and +(0,0) .. (bm);
  \draw[<-]  (tl) .. controls +( 0, -0.5) and +(0,0) .. (tm) coordinate[pos = 0.25] node[pos= 0.35 , draw=\crd, fill=white,regular polygon, regular polygon sides=3,inner sep=0.9pt] {} node [pos=0.3, left]{$\lambda$};
  \draw[<-]  (tr) .. controls +( 0, -0.5) and +(0,0) .. (tm)
  coordinate[pos = 0.25] node[pos= 0.35 , draw=\crd, fill=white,regular polygon, regular polygon sides=3,inner sep=0.9pt] {} node[pos=0.3, right]{$\mu$};
  \draw [double] (bm) -- (tm) node[pos=0, draw=\crd, fill= white, circle , inner sep=1.5pt]{} node[pos=0, right]{$s$};
\end{scope}

%% file: sch/reidem_lemm_19.tex
\begin{scope}[font=\tiny]
  \draw[->] (0.5, -0.5) ..controls +(0,0.3) and +(0,-0.3) .. (-0.5, 0.5) node[pos=1, above] {} coordinate[pos =0.2] (t1) ;
  \fill[white] (0,0) circle (2mm);
  \draw[->] (-0.5, -0.5) ..controls +(0,0.3) and +(0,-0.3) .. (0.5, 0.5) node[pos=1, above] {} coordinate[pos =0.2] (t2);
  \filldraw[draw= \crd, fill = white] (t2) circle (1mm)
  node[left] {$\omega$};
  \filldraw[draw= \crd, fill = white] (t1) circle (1mm)
  node[right] {$\omega$};
\end{scope}

%% file: sch/reidem_lemm_20.tex
\begin{scope}[font=\tiny]
  \draw[->] (0.5, -0.5) ..controls +(0,0.3) and +(0,-0.3) .. (-0.5,0.5) node[pos=1, above] {} coordinate[pos =0.75] (t1);
  \fill[white] (0,0) circle (2mm);
  \draw[->] (-0.5, -0.5) ..controls +(0,0.3) and +(0,-0.3) .. (0.5, 0.5) node[pos=1, above] {} coordinate[pos =0.75] (t2);
  \filldraw[draw= \crd, fill = white] (t2) circle (1mm)
  node[right] {$\omega$};
  \filldraw[draw= \crd, fill = white] (t1) circle (1mm)
  node[left] {$\omega$};
\end{scope}

%% file: sch/reidem_lemm_1.tex
\begin{scope}[font=\tiny]
  \draw[->] (0.5, -0.5) ..controls +(0,0.3) and +(0,-0.3) .. (-0.5, 0.5) node[pos=1, above] {} coordinate[pos =0.2] (t1) ;
  \fill[white] (0,0) circle (2mm);
  \draw[->] (-0.5, -0.5) ..controls +(0,0.3) and +(0,-0.3) .. (0.5, 0.5) node[pos=1, above] {} coordinate[pos =0.2] (t2);
  \filldraw[draw= \crd, fill = white] (t2) circle (1mm)
  node[left] {$\omega$};
  \filldraw[draw= \crd, fill = white] (t1) circle (1mm)
  node[right] {$\epsilon$};
\end{scope}

%% file: sch/reidem_lemm_2.tex
\begin{scope}[font=\tiny]
  \draw[->] (0.5, -0.5) ..controls +(0,0.3) and +(0,-0.3) .. (-0.5,0.5) node[pos=1, above] {} coordinate[pos =0.75] (t1);
  \fill[white] (0,0) circle (2mm);
  \draw[->] (-0.5, -0.5) ..controls +(0,0.3) and +(0,-0.3) .. (0.5, 0.5) node[pos=1, above] {} coordinate[pos =0.75] (t2);
  \filldraw[draw= \crd, fill = white] (t2) circle (1mm)
  node[right] {$\omega$};
  \filldraw[draw= \crd, fill = white] (t1) circle (1mm)
  node[left] {$\epsilon$};
\end{scope}

%% file: sch/reidem_lemm_3.tex
\begin{scope}[font=\tiny]
  \draw[->] (-0.5, -0.5) ..controls +(0,0.3) and +(0,-0.3) .. (0.5,0.5) node[pos=1, above] {} coordinate[pos =0.2] (t1);
  \fill[white] (0,0) circle (2mm);
  \draw[->] (0.5,-0.5) ..controls +(0,0.3) and +(0,-0.3) .. (-0.5, 0.5) node[pos=1, above] {} coordinate[pos =0.2] (t2);
  \filldraw[draw= \crd, fill = white] (t2) circle (1mm)
  node[right] {$\epsilon$};
  \filldraw[draw= \crd, fill = white] (t1) circle (1mm)
  node[left] {$\omega$};
\end{scope}

%% file: sch/reidem_lemm_4.tex
\begin{scope}[font=\tiny]
  \draw[->] (-0.5, -0.5) ..controls +(0,0.3) and +(0,-0.3) .. (0.5,0.5) node[pos=1, above] {} coordinate[pos =0.75] (t1);
  \fill[white] (0,0) circle (2mm);
  \draw[->] (0.5,-0.5) ..controls +(0,0.3) and +(0,-0.3) .. (-0.5, 0.5) node[pos=1, above] {} coordinate[pos =0.75] (t2);
  \filldraw[draw= \crd, fill = white] (t2) circle (1mm)
  node[left] {$\epsilon$};
  \filldraw[draw= \crd, fill = white] (t1) circle (1mm)
  node[right] {$\omega$};
\end{scope}

%% file: sch/reidem_lemm_23.tex
\begin{scope}[font=\tiny]
  \draw[->] (0.5, -0.5) ..controls +(0,0.3) and +(0,-0.3) .. (-0.5, 0.5) node[pos=1, above] {} coordinate[pos =0.2] (t1) ;
  \fill[white] (0,0) circle (2mm);
  \draw[->] (-0.5, -0.5) ..controls +(0,0.3) and +(0,-0.3) .. (0.5, 0.5) node[pos=1, above] {} coordinate[pos =0.15] (t2) coordinate[pos=0.35] (t3);
  \filldraw[draw= \crd, fill = white] (t1) circle (1mm)
  node[right] {$\epsilon$};
  \filldraw[draw= \crd, fill = white] (t2) circle (1mm)
  node[left] {$\epsilon$};
  \filldraw[draw= \crd, fill = white] (t3) circle (1mm)
  node[left] {$\omega-\epsilon$};
\end{scope}

%% file: sch/reidem_lemm_24.tex
\begin{scope}[font=\tiny]
  \draw[->] (0.5, -0.5) ..controls +(0,0.3) and +(0,-0.3) .. (-0.5, 0.5) node[pos=1, above] {} coordinate[pos =0.2] (t1) ;
  \fill[white] (0,0) circle (2mm);
  \draw[->] (-0.5, -0.5) ..controls +(0,0.3) and +(0,-0.3) .. (0.5, 0.5) node[pos=1, above] {} coordinate[pos =0.2] (t2) coordinate[pos=0.7] (t3);
  \filldraw[draw= \crd, fill = white] (t1) circle (1mm)
  node[right] {$\epsilon$};
  \filldraw[draw= \crd, fill = white] (t2) circle (1mm)
  node[left] {$\epsilon$};
  \filldraw[draw= \crd, fill = white] (t3) circle (1mm)
  node[right] {$\omega-\epsilon$};
\end{scope}

%% file: sch/reidem_lemm_5.tex
\begin{scope}[font=\tiny]
  \draw[->] (0.5, -0.5) ..controls +(0,0.3) and +(0,-0.3) .. (-0.5, 0.5) node[pos=1, above] {} coordinate[pos =0.2] (t1);
  \fill[white] (0,0) circle (2mm);
  \draw[->] (-0.5, -0.5) ..controls +(0,0.3) and +(0,-0.3) .. (0.5, 0.5) node[pos=1, above] {} coordinate[pos =0.2] (t2);
  \filldraw[draw= \crd, fill = \csrd] (t2) circle (1mm)
  node[left] {$\omega$};
  \filldraw[draw= \crd, fill = \csrd] (t1) circle (1mm)
  node[right] {$\epsilon$};
\end{scope}

%% file: sch/reidem_lemm_6.tex
\begin{scope}[font=\tiny]
  \draw[->] (0.5, -0.5) ..controls +(0,0.3) and +(0,-0.3) .. (-0.5,0.5) node[pos=1, above] {} coordinate[pos =0.75] (t1);
  \fill[white] (0,0) circle (2mm);
  \draw[->] (-0.5, -0.5) ..controls +(0,0.3) and +(0,-0.3) .. (0.5, 0.5) node[pos=1, above] {} coordinate[pos =0.75] (t2);
  \filldraw[draw= \crd, fill = \csrd] (t2) circle (1mm)
  node[right] {$\omega$};
  \filldraw[draw= \crd, fill = \csrd] (t1) circle (1mm)
  node[left] {$\epsilon$};
\end{scope}

%% file: sch/reidem_lemm_7.tex
\begin{scope}[font=\tiny]
  \draw[->] (-0.5, -0.5) ..controls +(0,0.3) and +(0,-0.3) .. (0.5,0.5) node[pos=1, above] {} coordinate[pos =0.2] (t1);
  \fill[white] (0,0) circle (2mm);
  \draw[->] (0.5,-0.5) ..controls +(0,0.3) and +(0,-0.3) .. (-0.5, 0.5) node[pos=1, above] {} coordinate[pos =0.2] (t2);
  \filldraw[draw=\crd, fill = \csrd] (t2) circle (1mm)
  node[right] {$\epsilon$};
  \filldraw[draw= \crd, fill = \csrd] (t1) circle (1mm)
  node[left] {$\omega$};
\end{scope}

%% file: sch/reidem_lemm_8.tex
\begin{scope}[font=\tiny]
  \draw[->] (-0.5, -0.5) ..controls +(0,0.3) and +(0,-0.3) .. (0.5,0.5) node[pos=1, above] {} coordinate[pos =0.75] (t1);
  \fill[white] (0,0) circle (2mm);
  \draw[->] (0.5,-0.5) ..controls +(0,0.3) and +(0,-0.3) .. (-0.5, 0.5) node[pos=1, above] {} coordinate[pos =0.75] (t2);
  \filldraw[draw= \crd, fill = \csrd] (t2) circle (1mm)
  node[left] {$\epsilon$};
  \filldraw[draw= \crd, fill = \csrd] (t1) circle (1mm)
  node[right] {$\omega$};
\end{scope}

%% file: sch/reidem_lemm_27.tex
\begin{scope}[font=\tiny]
  \draw[->] (0.5, -0.5) ..controls +(0,0.3) and +(0,-0.3) .. (-0.5, 0.5) node[pos=1, above] {} coordinate[pos =0.2] (t1);
  \fill[white] (0,0) circle (2mm);
  \draw[->] (-0.5, -0.5) ..controls +(0,0.3) and +(0,-0.3) .. (0.5, 0.5) node[pos=1, above] {} coordinate[pos =0.2] (t2);
  \filldraw[draw= \crd, fill = \csrd] (t2) circle (1mm)
  node[left] {$\omega$};
\end{scope}

%% file: sch/reidem_lemm_28.tex
\begin{scope}[font=\tiny]
  \draw[->] (0.5, -0.5) ..controls +(0,0.3) and +(0,-0.3) .. (-0.5,0.5) node[pos=1, above] {} coordinate[pos =0.75] (t1);
  \fill[white] (0,0) circle (2mm);
  \draw[->] (-0.5, -0.5) ..controls +(0,0.3) and +(0,-0.3) .. (0.5, 0.5) node[pos=1, above] {} coordinate[pos =0.75] (t2);
  \filldraw[draw= \crd, fill = \csrd] (t2) circle (1mm)
  node[right] {$\omega$};
\end{scope}

%% file: sch/relation_hollow_black_dots_2.tex
\begin{scope}
  \coordinate (m) at (  0,  0);
  \coordinate (t) at (  2, 0);
  \draw[->] (m) -- (t) node[pos = 1, right] {$a$} coordinate[pos=0.5] (ga);
  \filldraw[draw= \crd, fill = \csrd] (ga) circle (1mm) node[yshift= 3mm] {$\omega$};
\end{scope}

%% file: sch/relation_hollow_black_dots_1.tex
\begin{scope}
  \coordinate (m) at (0,0);
  \coordinate (t) at (2,0);
  \draw[->] (m) -- (t) node[pos = 1, right] {$a$} coordinate[pos=0.2] (ga) coordinate[pos=0.8] (gb) coordinate[pos=0.5,yshift=6mm] (gc) coordinate[pos=0.5,yshift=-6mm] (gd);
  \filldraw[draw= \crd, fill = white] (ga) circle (1mm) node[yshift= 3.5mm] {$-\omega$};
  \filldraw[draw= black, fill = black] (gb) circle (1mm) node[yshift=3.5mm] {$\omega$};
\end{scope}

%% file: sch/reidem_lemm_29.tex
\begin{scope}[font=\tiny]
  \draw[->] (0.5, -0.5) ..controls +(0,0.3) and +(0,-0.3) .. (-0.5, 0.5) node[pos=1, above] {} coordinate[pos =0.2] (t1);
  \fill[white] (0,0) circle (2mm);
  \draw[->] (-0.5, -0.5) ..controls +(0,0.3) and +(0,-0.3) .. (0.5, 0.5) node[pos=1, above] {} coordinate[pos =0.2] (t2);
  \filldraw[draw=black, fill = black] (t2) circle (1mm)
  node[left] {$\omega$};
\end{scope}

%% file: sch/reidem_lemm_30.tex
\begin{scope}[font=\tiny]
  \draw[->] (0.5, -0.5) ..controls +(0,0.3) and +(0,-0.3) .. (-0.5,0.5) node[pos=1, above] {} coordinate[pos =0.75] (t1);
  \fill[white] (0,0) circle (2mm);
  \draw[->] (-0.5, -0.5) ..controls +(0,0.3) and +(0,-0.3) .. (0.5, 0.5) node[pos=1, above] {} coordinate[pos =0.75] (t2);
  \filldraw[draw=black, fill = black] (t2) circle (1mm)
  node[right] {$\omega$};
\end{scope}

%% file: SF_RI.tex
\begin{prop}
    There are isomorphisms in the relative homotopy category:
    \begin{equation} 
    \label{reid1}
    \NB{\tikz[font=\tiny, scale=0.6, yscale=-1]{\input{sch/reidem1_3}}} \cong t^{-1}q^{N-1}\NB{\tikz[font=\tiny, scale=0.6]{\input{sch/reidem1_4}}} \qquad \text{and} \qquad   \NB{\tikz[font=\tiny, scale=0.6]{\input{sch/reidem1_1}}} \cong t q^{1-N} \NB{\tikz[font=\tiny, scale=0.6]{\input{sch/reidem1_2}}}.
    \end{equation}
\end{prop}

\begin{proof}
We only present the proof for the first isomorphism since the second uses the same ideas.
Thanks to the Lemma \ref{migr_srd}, the existence of the following isomorphism is enough: 
\begin{equation}
	\label{iso_rei1}
	\NB{\tikz[font=\tiny, scale=0.6, yscale=-1]{\input{sch/reidem1_5}}} \cong t^{-1}q^{N-1}\NB{\tikz[font=\tiny, scale=0.6]{\input{sch/reidem1_6}}}.
\end{equation}
The first complex is defined to be:
\[
\NB{\tikz[font=\tiny, scale=0.6,yscale=-1]{\input{sch/reidem1_5}}} = \NB{ \tikz[xscale = 3, yscale = 3]{
\node (AT) at (0,0) {$t^{-1}\NB{\tikz[font=\tiny, scale=0.6]{\input{sch/reidem1_7}}}$};
\node (CT) at (1.5,0) {$q^{-1}\NB{\tikz[font=\tiny, scale=0.6]{\input{sch/reidem1_8}}}$};
\draw[->] (AT) -- (CT) node[pos=0.5, above]{$\mapzip$}; }}.
\]
These two complexes embed in a short exact sequence of $\Witt$-equivariant complexes that splits after forgetting the $\Witt$-action: 

\[
\NB{\tikz[xscale = 3, yscale = 2.8]{
	\node (i0) at (0, 0) {$t^{-1}\NB{\tikz[font= \tiny, scale=.6]{\input{sch/reidem1_7}}}$};
	\node (i1) at (1.5, 0) { $q^{-1}\NB{\tikz[font= \tiny,scale=0.6]{\input{sch/reidem1_8}}}$ };
	\node (i2) at (0, 1) {$t^{-1}q^{N-1}\NB{\tikz[font= \tiny, scale=0.6]{\input{sch/reidem1_6}}}$};
	\node (i3) at (1.5, -1.6) {$\oplus [N-1] \NB{\tikz[font= \tiny, scale=0.6]{\input{sch/reidem1_6}}}$};
	\node (i4) at (0, -1.6) {$t^{-1}\oplus [N-1] \NB{\tikz[font= \tiny, scale=0.6]{\input{sch/reidem1_6}}}$};
\draw[->] (i0) -- (i1)  node[pos=0.5, above] {$\mapzip$};
\draw[->] (i0) -- (i2)  node[pos=0.5, left] {\NB{\tikz[font= \tiny, scale=0.9,-]{\input{sch/foam_cap_empty}}}};
	\draw[->] (i4) -- (i0)  node[pos=0.5, left] {$\begin{pmatrix} 
	\NB{\tikz[font= \tiny, scale=0.9,-]{\input{sch/foam_cup_empty}}} \\
	\vdots \\
	\NB{\tikz[font= \tiny, scale=0.9,-]{\input{sch/foam_cup_decoN-2}}}
	\end{pmatrix}$};
\draw[->] (i3) -- (i1) node[pos=0.5, right] {$\begin{pmatrix} 
	\mapzip \circ \NB{\tikz[font= \tiny, scale=0.9,-]{\input{sch/foam_cup_empty}}} \\
	\vdots \\
	\mapzip \circ \NB{\tikz[font= \tiny, scale=0.9,-]{\input{sch/foam_cup_decoN-2}}} 
	\end{pmatrix}$};
	\draw[->] (i4) -- (i3)  node[pos=0.5, above] {$\Id$};
}}.
\]
\vspace{0.5cm}

The splitting of the bottom complex does not respect the \(\Witt\)-structure, but
by Corollary \ref{lemme_RHC} there is an isomorphism between the top two rows in the relative homotopy category, since it is contractible.
\end{proof}

%% file: sch/reidem1_3.tex
\begin{scope}[xshift = 2.5cm]
  \draw  (0, 1) -- (0, 0.5) .. controls +(0,-0.5) and +(0, 0.5)
  .. (1, -0.5) arc (180:360:0.5) -- (2,0);  
    \fill[white] (0.5, 0) circle (1mm);
  \draw[<-] (0, -1) -- (0, -0.5) .. controls +(0,0.5) and +(0, -0.5)
  .. (1, 0.5) arc (180:0:0.5) -- (2,0);
\end{scope}

%% file: sch/reidem1_4.tex
\begin{scope}
  \draw [->] (0, -1) -- +(0,2);
       \coordinate (B) at (0, 0);
        \filldraw[draw= \crd, fill = \csrd] (B) circle (1mm)
  node[left]{$\frac{1}{2}$};
\end{scope}

%% file: sch/reidem1_1.tex
\begin{scope}[xshift = 2.5cm]
  \draw [<-] (0, 1) -- (0, 0.5) .. controls +(0,-0.5) and +(0, 0.5) .. (1, -0.5) arc (180:360:0.5) -- (2,0);  
    \fill[white] (0.5, 0) circle (1mm);
  \draw (0, -1) -- (0, -0.5) .. controls +(0,0.5) and +(0, -0.5) .. (1, 0.5) arc (180:0:0.5) -- (2,0);
\end{scope}

%% file: sch/reidem1_2.tex
\begin{scope}
  \draw [->] (0, -1) -- +(0,2);
       \coordinate (B) at (0, 0);
        \filldraw[draw= \crd, fill = \csrd] (B) circle (1mm)
  node[left]{$-\frac{1}{2}$};
\end{scope}

%% file: sch/reidem1_6.tex
\begin{scope}
  \draw [->] (0, -1) -- +(0,2);
   \coordinate (O) at (1, 0);
       \filldraw[draw= \crd, fill = \csrd] (O);
     \coordinate (A) at (0, -.6);
        \filldraw[draw= \crd, fill = white] (A);
\end{scope}

%% file: sch/reidem1_7.tex
\begin{scope}
    \coordinate (bl) at (-0.5, -1);
    \coordinate (br) at ( 0.5, -1);
    \coordinate (tl) at (-0.5,  1);
    \coordinate (tr) at ( 0.5,  1);
    \coordinate (Oo) at (2.5,0);
    \coordinate (ml) at (-0.5,  -.8);
    \coordinate (Ml) at (-0.5,  .8);
    \coordinate (mr) at (0.5,  -.6);
    \coordinate (Mr) at (0.5,  .6);
    \draw[->] (bl) -- (tl);  
    \draw[->] (br) -- (tr);
    \draw (.5, 1) arc (180:0:0.5) -- (1.5,-1) arc (180:0:-0.5);
    \coordinate (A) at (1.5, -.8);
    \filldraw[draw= \crd, fill = \csrd] (A) circle (1mm) node[right] {$-\frac{1}{2}$};
    \coordinate (B) at (1.5, .8);
\end{scope}

%% file: sch/reidem1_8.tex
\begin{scope}
  \coordinate (bl) at (-0.5, -1);
  \coordinate (br) at ( 0.5, -1);
  \coordinate (bm) at (  0,-0.3);
  \coordinate (tl) at (-0.5,  1);
  \coordinate (tr) at ( 0.5,  1);
  \coordinate (tm) at (  0, 0.3);
  \draw[>-]  (bl) .. controls +( 0, 0.5) and +(0,0) .. (bm);
  \draw[>-]  (br) .. controls +( 0, 0.5) and +(0,0) .. (bm);
  \draw[<-]  (tl) .. controls +( 0, -0.5) and +(0,0) .. (tm) coordinate[pos = 0.25] node[pos= 0.35 , draw=\crd, fill=white,regular polygon, regular polygon sides=3,inner sep=0.9pt] {} node [pos=0.3, left]{$\lambda$};
  \draw[<-]  (tr) .. controls +( 0, -0.5) and +(0,0) .. (tm)
  coordinate[pos = 0.25] node[pos= 0.35 , draw=\crd, fill=white,regular polygon, regular polygon sides=3,inner sep=0.9pt] {} node[pos=0.3, right]{$\mu$};
  \draw [double] (bm) -- (tm) node[pos=0, draw=\crd, fill= white, circle , inner sep=1.5pt]{} node[pos=0, right]{$-\bar{s}$};
  \draw (.5, 1) arc (180:0:0.5) -- (1.5,-.85) arc (180:0:-0.5);
    \coordinate (A) at (1.5, -.8);
    \coordinate (B) at (1.5, .8);
    \filldraw[draw= \crd, fill = \csrd] (A) circle (1mm)
    node[right] {$-\frac{1}{2}$};
    \coordinate (Oo) at (2.7,0);
\end{scope}

%% file: sch/foam_cap_empty.tex
\begin{scope}[-]
  \draw (0,0) arc (180 :0: 0.5cm and 0.2cm);
  \draw[very thin] (0,0) arc (180 :0: 0.5cm and 0.6cm) node[pos=0.5,
  below] {};
  \draw (0,0) arc (180 :0: 0.5cm and -0.2cm)node[ below, pos =0.5] {};
\end{scope}

%% file: sch/foam_cup_empty.tex
\begin{scope}
  \draw (0,0) arc (180 :0: 0.5cm and 0.2cm) 
  node[above, pos=0.5] {};
  \draw[very thin] (0,0) arc (180 :0: 0.5cm and -0.6cm) node[pos=0.5, above] {};
  \draw (0,0) arc (180 :0: 0.5cm and -0.2cm);
\end{scope}

%% file: sch/foam_cup_decoN-2.tex
\begin{scope}
  \draw (0,0) arc (180 :0: 0.5cm and 0.2cm) 
  node[above, pos=0.5] {};
  \draw[very thin] (0,0) arc (180 :0: 0.5cm and -0.6cm) node[pos=0.5, above, scale=0.8] {$\deco_{N-2}$};
  \draw (0,0) arc (180 :0: 0.5cm and -0.2cm);
\end{scope}

%% file: SF_RII.tex
The following result can be proved by straightforward computations using Section~\ref{usual_morph_sec} and foam relations.
\begin{lem}
    \label{lemm_rei2}
    There is a short exact sequence of state spaces that splits when the $\Witt$-action is forgotten:
\begin{equation}
    \NB{\tikz[xscale = 3, yscale = 3]{\node (i0) at (0, 0) {\NB{\tikz[font= \tiny,scale=0.6]{\input{sch/reidem2a_3}}} };
    \node (i1) at (-1,0) { \NB{\tikz[font= \tiny, scale=0.6]{\input{sch/reidem2a_2}}} };
    \node (i3) at (-2, 0) { \NB{\tikz[font= \tiny,scale=0.6]{\input{sch/reidem2a_1}}} };
    \draw[->] (i3) -- (i1)  node[pos=0.5, above] {$\mapdcup$};
    \draw[->] (i1) -- (i0)  node[pos=0.5, above] {$\mapdcap$};
    \draw[-to] (i0) .. controls +(-.25, -0.55) and +(+.25, -0.55) .. (i1)node[pos=0.5, below] {$\mapdcup^x - {}^x\mapdcup$};
   \draw[-to] (i1) .. controls +(-.25, -0.55) and +(+.25, -0.55) .. (i3)node[pos=0.5, below] {$\mapdcap^x - {}^x\mapdcap$};
      }}
\end{equation}
where the splitting maps are provided by:
\[
\mapdcap^x-{}^x\mapdcap =
\NB{\tikz[scale=.9]{\input{sch/foam_digoncap_10}}} - \NB{\tikz[scale=.9]{\input{sch/foam_digoncap_01}}}
\hspace{1mm}; \hspace{3mm}
\mapdcup^x-{}^x\mapdcup =
\NB{\tikz[scale=.9]{\input{sch/foam_digoncup_10}}} - \NB{\tikz[scale=.9]{\input{sch/foam_digoncup_01}}}.
\]
\end{lem}

\begin{prop}
    \label{ReidemeisterII}
    There are isomorphisms in the relative homotopy category:
    \begin{equation}
        \NB{\tikz[font=\tiny,scale=.7]{\input{sch/reidem2_1}}} \cong \NB{\tikz[font=\tiny,scale=.7]{\input{sch/reidem2_2}}} \cong \NB{\tikz[font=\tiny, xscale=-1,scale=.7]{\input{sch/reidem2_1}}}.
    \end{equation}
\end{prop}

\begin{proof}
The constructions of both isomorphisms are the same so we only prove the first. By definition  
\begin{equation}
    \NB{\tikz[scale=.7]{\input{sch/reidem2_1}}} \hspace{5mm} = \hspace{5mm}
   \NB{\tikz[xscale = 4, yscale = 2]{
    \node (i0) at (0, 0) {$q^{-1}\NB{\tikz[font=\tiny,scale=0.6]{\input{sch/reidem2a_4}}}$}; 
    \node (i1) at (-.75, .6) {\NB{\tikz[font= \tiny,scale=0.6]{\input{sch/reidem2_2}}}};
    \node (i2) at (-.75, -.6) {\NB{\tikz[font= \tiny, scale=0.6]{\input{sch/reidem2a_6}}}};
    \node (i3) at (-1.5, 0) {$q\NB{\tikz[font=\tiny, scale=0.6]{\input{sch/reidem2a_5}}}$}; 
    \draw[->] (i3) -- (i1) node[pos=0.5, above]{$\mapunzip$};
    \draw[->] (i3) -- (i2) node[pos=0.5, above]{$-\mapzip$};
    \draw[->] (i1) -- (i0) node[pos=0.5, above]{$\mapzip$};
    \draw[->] (i2) -- (i0) node[pos=0.5, above]{$\mapunzip$};
    }}.
\end{equation}

Thanks to Lemma \ref{lemm_rei2} this complex fits into in a short exact sequence of complexes that splits when forgetting the $\Witt$-action. 

    \begin{equation}
    \NB{\tikz[xscale = 3.5, yscale = 2.25]{
    \node (i0) at (0, 0) {$q^{-1}\NB{\tikz[font=\tiny,scale=0.6]{\input{sch/reidem2a_4}}}$}; 
    \node (i1) at (-1, .5) {\NB{\tikz[font= \tiny,scale=0.6]{\input{sch/reidem2_2}}}};
    \node (i2) at (-1, -.5) {\NB{\tikz[font= \tiny, scale=0.6]{\input{sch/reidem2a_6}}}};
    \node (i3) at (-2, 0) {$q\NB{\tikz[font=\tiny, scale=0.6]{\input{sch/reidem2a_5}}}$}; 
    \draw[->] (i3) -- (i1) node[pos=0.5, above]{$\mapunzip$};
    \draw[->] (i3) -- (i2) node[pos=0.5, above]{$-\mapzip$};
    \draw[->] (i1) -- (i0) node[pos=0.5, above]{$\mapzip$};
    \draw[->] (i2) -- (i0) node[pos=0.5, above]{$\mapunzip$};
    \node (i5) at (-2.2, -.9) { };
    \node (i6) at (0.21, .65) { };
    
    \node (i4) at (-1, 1.5) {$q^{-1} \NB{\tikz[font= \tiny,scale=0.6]{\input{sch/reidem2a_4}}}$ };
    \node (i7) at (-1, 2.5) {\NB{\tikz[font = \tiny,scale=0.6]{\input{sch/reidem2_2}}}};
    \node (i8) at (-2,2) {$0$};
    \node (i9) at (0,2) {$q^{-1}\NB{\tikz[font = \tiny,scale=0.6]{\input{sch/reidem2a_4}}}$};
    \draw[->] (i8) -- (i4);
    \draw[->] (i8) -- (i7);
    \draw[->] (i4) -- (i9) node[pos=0.5, above]{$\Id$};
    \draw[->] (i7) -- (i9) node[pos=0.5, above]{$\mapzip$};
    \node (i10) at (-2.2, 1.35) { };
    \node (i11) at (0.21, 2.65) { };

    \node (i12) at (-1, -1.5) {$0$};
    \node (i13) at (-1, -2.5) {\NB{\tikz[font = \tiny,scale=0.6]{\input{sch/reidem2a_5}}}};
    \node (i14) at (-2,-2) {$q\NB{\tikz[font = \tiny,scale=0.6]{\input{sch/reidem2a_5}}}$};
    \node (i15) at (0,-2) {$0$};
    \draw[->] (i14) -- (i12);
    \draw[->] (i14) -- (i13) node[pos=0.5, above]{$-\Id$};
    \draw[->] (i12) -- (i15);
    \draw[->] (i13) -- (i15);
    \node (i16) at (-2.2, -2.65) { };
    \node (i17) at (0.21, -1.53) { };

    \draw[->] (i2) -- (i12);
    \draw[->] (i4) -- (i1);
    \node[draw, densely dotted ,fit=(i5) (i6) ,inner sep=1ex,rounded corners] {};
    \node[draw, densely dotted ,fit=(i10) (i11) ,inner sep=1ex,rounded corners] {};
    \node[draw, densely dotted ,fit=(i16) (i17) ,inner sep=1ex,rounded corners] {};
    }} .
    \end{equation}

The last complex is contractible, thus by using Corollary \ref{lemme_RHC}, the first and the second complexes are isomorphic in the relative homotopy category. Moreover, in the first complex, the terms
\[ \NB{\tikz[font = \tiny,scale=0.6]{\input{sch/reidem2a_4}}} \ ,\]
can be contracted resulting in the desired isomorphism.\qedhere 
\end{proof}

Note that the proof of the previous proposition does not require any result of Section \ref{sliding_sect}.

\begin{prop}
    There is an isomorphism in the relative homotopy category: 
     \begin{equation}
         \NB{\tikz[font=\tiny]{\input{sch/reidem2b_1}}} \cong \NB{\tikz[font=\tiny]{\input{sch/reidem2b_2}}}.
     \end{equation}
\end{prop}

\begin{proof}
Additivity and sliding of dots in Lemma \ref{migr_srd} implies the isomorphism
\[\NB{\tikz[font=\tiny]{\input{sch/reidem2b_2}}} \cong \NB{\tikz[font=\tiny]{\input{sch/reidem2b_2p}}}.\]
By definition, the second complex is given by:
\begin{align}
    \NB{\tikz[xscale = 4, yscale = 2.5]{\node (i0) at (0, 0) {$q^{-1}\NB{\tikz[font=\tiny,scale=0.6]{\input{sch/reidem2b_6}}}$}; 
    \node (i1) at (-.95, .6) {\NB{\tikz[font= \tiny,scale=0.6]{\input{sch/reidem2b_4}}}};
    \node (i2) at (-.95, -.6) {\NB{\tikz[font= \tiny, scale=0.6]{\input{sch/reidem2b_5}}}};
    \node (i3) at (-1.9, 0) {$q\NB{\tikz[font=\tiny, scale=0.6]{\input{sch/reidem2b_3}}}$}; 
    \draw[->] (i3) -- (i1) node[pos=0.5, above]{$\mapunzip$};
    \draw[->] (i3) -- (i2) node[pos=0.5, above]{$-\mapzip$};
    \draw[->] (i1) -- (i0) node[pos=0.5, above]{$\mapzip$};
    \draw[->] (i2) -- (i0) node[pos=0.5, above]{$\mapunzip$}; }}.
\end{align}

The following map is a homotopy equivalence when forgetting the \(\Witt\)-action, so it induces an isomorphism in the relative homotopy category:
\begin{equation}
    \label{isomR2b}
    \NB{\tikz[xscale = 4, yscale = 2.5]{\node (i0) at (0, 0) {$q^{-1}\NB{\tikz[font=\tiny,scale=0.6]{\input{sch/reidem2b_6}}}$}; 
    \node (i1) at (-.95, .6) {\NB{\tikz[font= \tiny,scale=0.6]{\input{sch/reidem2b_4}}}};
    \node (i2) at (-.95, -.6) {\NB{\tikz[font= \tiny, scale=0.6]{\input{sch/reidem2b_5}}}};
    \node (i3) at (-1.9, 0) {$q\NB{\tikz[font=\tiny, scale=0.6]{\input{sch/reidem2b_3}}}$}; 
    \node (i4) at (-.95, 1.5) { \NB{\tikz[font= \tiny,scale=0.6]{\input{sch/reidem2b_1}}} };
    \draw[->] (i3) -- (i1) node[pos=0.5, above]{$\mapunzip$};
    \draw[->] (i3) -- (i2) node[pos=0.5, above]{$-\mapzip$};
    \draw[->] (i1) -- (i0) node[pos=0.5, above]{$\mapzip$};
    \draw[->] (i2) -- (i0) node[pos=0.5, above]{$\mapunzip$};
    \node (i5) at (-2.2, -1.05) { };
    \node (i6) at (0.35, .75) { };
   \draw[->] (-.95, 1.2) -- (i1) node[pos=0.5, left, font=\small] {$\begin{pmatrix} \mapcup \circ \mapsaddle \\ 
   \mapsaddle \circ \mapzip \mapzip \circ \mapcup \mapcup
    \end{pmatrix}$};
    \node[draw, densely dotted ,fit=(i4) (i4) ,inner sep=1ex, rounded corners] {};
    \node[draw, densely dotted ,fit=(i5) (i6) ,inner sep=1ex,rounded corners] {};
    }}    .
\end{equation}

Let us check that the morphisms are indeed $\Witt$-equivariant. 
The map \(\mapsaddle \circ \mapzip \mapzip \circ \mapcup \mapcup\) is represented by the movie
\[
\mymoviefour[scale=0.75]{\NB{\tikz[scale=.85]{\input{sch/reidem2b_1}}}}{\NB{\tikz[scale=.8]{\input{sch/movie_r2b_1}}}}{\NB{\tikz[scale=.95]{\input{sch/movie_r2b_2}}}}{\NB{\tikz[scale=0.55]{\input{sch/movie_r2b_3}}}}
\]
Lemmas (\ref{rd_cp}), (\ref{rd_saddle}), and (\ref{rd_zip}) show that the morphism $\mapsaddle \circ \mapzip \mapzip \circ \mapcup \mapcup$ is equivariant with the following twist 
\[\mapsaddle\circ\mapzip\mapzip \circ \mapcup \mapcup \, \, : \, \, \NB{\tikz[font=\tiny, scale=0.9]{\input{sch/reidem2b_1}}} \longrightarrow  \NB{\tikz[font=\tiny, scale=0.6]{\input{sch/reidem2b_7}}}=
\raisebox{2.5mm}{\NB{\tikz[font=\tiny, scale=0.6]{\input{sch/reidem2b_12}}}}
. \]
The identity is due to Lemma \ref{lemm_rd_3}.
In this case, the relation takes the form
\[\NB{\tikz[font=\tiny, scale=0.9,yscale=-1]{\input{sch/reidem2b_8}}} = \NB{\tikz[font=\tiny, scale=0.9,yscale=-1]{\input{sch/reidem2b_9}}}.\]

Thus, the map $\mapsaddle \circ \mapzip \mapzip \circ \mapcup \mapcup$ in (\ref{isomR2b}) is indeed $\Witt$-equivariant. The easier computation for $\mapcup \circ \mapsaddle$ is similar to the previous one and is thus omitted.
\end{proof}

%% file: sch/reidem2a_3.tex
\begin{scope}
  \coordinate (bl) at (-0.5, -1);
  \coordinate (br) at ( 0.5, -1);
  \coordinate (bm) at (  0,-0.3);
  \coordinate (tl) at (-0.5,  1);
  \coordinate (tr) at ( 0.5,  1);
  \coordinate (tm) at (  0, 0.3);
  \draw[>-]  (bl) .. controls +( 0, 0.5) and +(0,0) .. (bm);
  \draw[>-]  (br) .. controls +( 0, 0.5) and +(0,0) .. (bm);
  \draw[<-]  (tl) .. controls +( 0, -0.5) and +(0,0) .. (tm);
  \draw[<-]  (tr) .. controls +( 0, -0.5) and +(0,0) .. (tm);
  \draw [double] (bm) -- (tm) node[pos=1, draw=\crd, fill= white, circle , inner sep=1.5pt]{} node[pos=1, right, yshift=-.5mm]{$1$};
\end{scope}

%% file: sch/reidem2a_2.tex
\begin{scope}
    \coordinate (bl) at (-0.5, -1);
    \coordinate (br) at ( 0.5, -1);
    \coordinate (bm) at (  0,-0.3);
    \coordinate (tl) at (-0.5,  1);
    \coordinate (tr) at ( 0.5,  1);
    \coordinate (tm) at (  0, 0.3);
    \coordinate (brr) at ( 1.5, -1);
    \coordinate (trr) at ( 1.5, 1);
    \draw[>-]  (bl) .. controls +( 0, 0.5) and +(0,0) .. (bm) ;
    \draw[>-]  (br) .. controls +( 0, 0.5) and +(0,0) .. (bm);
    \draw[<-]  (tl) .. controls +( 0, -0.5) and +(0,0) .. (tm);
    \draw[<-]  (tr) .. controls +( 0, -0.5) and +(0,0) .. (tm);
    \draw [double] (bm) -- (tm) node[pos=1, draw=\crd, fill= white, circle , inner sep=1.5pt]{} node[pos=1, right, yshift=-.5mm]{$s$}; 
    \coordinate (blU) at (-0.5, 1);
    \coordinate (brU) at ( 0.5, 1);
    \coordinate (bmU) at (  0,1.7);
    \coordinate (tlU) at (-0.5,  3);
    \coordinate (trU) at ( 0.5,  3);
    \coordinate (tmU) at (  0, 2.3);
    \coordinate (brrU) at ( 1.5, 1);
    \coordinate (trrU) at ( 1.5, 3);
    \draw[]  (blU) .. controls +( 0, 0.5) and +(0,0) .. (bmU) node[pos= 0.35 , draw=\crd, fill=white,regular polygon, regular polygon sides=3,inner sep=0.9pt] {} node[pos=0.3, left]{$\lambda$};
    \draw[]  (brU) .. controls +( 0, 0.5) and +(0,0) .. (bmU) node[pos= 0.35 , draw=\crd, fill=white,regular polygon, regular polygon sides=3,inner sep=0.9pt] {} node[pos=0.3, right]{$\mu$};
    \draw[<-]  (tlU) .. controls +( 0, -0.5) and +(0,0) .. (tmU) coordinate[pos = 0.25] (gaU) ;
    \draw[<-]  (trU) .. controls +( 0, -0.5) and +(0,0) .. (tmU) coordinate[pos = 0.25] (gbU) ;  
    \draw [double] (bmU) -- (tmU);
\end{scope}

%% file: sch/reidem2a_1.tex
\begin{scope}
  \coordinate (bl) at (-0.5, -1);
  \coordinate (br) at ( 0.5, -1);
  \coordinate (bm) at (  0,-0.3);
  \coordinate (tl) at (-0.5,  1);
  \coordinate (tr) at ( 0.5,  1);
  \coordinate (tm) at (  0, 0.3);
  \draw[>-]  (bl) .. controls +( 0, 0.5) and +(0,0) .. (bm);
  \draw[>-]  (br) .. controls +( 0, 0.5) and +(0,0) .. (bm);
  \draw[<-]  (tl) .. controls +( 0, -0.5) and +(0,0) .. (tm);
  \draw[<-]  (tr) .. controls +( 0, -0.5) and +(0,0) .. (tm);
  \draw [double] (bm) -- (tm);
\end{scope}

%% file: sch/foam_digoncap_10.tex
\begin{scope}[font =\tiny]
  \begin{scope}
    \coordinate (L) at (0,0);
    \coordinate (R) at (2,0);
    \coordinate (ML) at (0.5, 0);
    \coordinate (MR) at (1.5, 0);
    \draw[->-] (ML).. controls + (0.4, 0.4) and +(-0.2, 0.4) .. (MR)
    node[above, pos =0.3] {}; 
    \draw[->-] (ML).. controls + (0.2, -0.4) and +(-0.4, -0.4) .. (MR)
    node[above, pos =0.7, yshift = -0.4mm] {$\deco_1$}; 
  \end{scope}  
 \begin{scope}[yshift = 1.3cm]
    \coordinate (LB) at (0,0);
    \coordinate (RB) at (2,0);
  \end{scope}  
  \draw[double] (MR)--(R) -- (RB) -- (LB) -- (L)-- (ML);
  \draw (ML) .. controls +(0, 1) and +(0, 1) .. (MR);
\end{scope}

%% file: sch/foam_digoncap_01.tex
\begin{scope}[font =\tiny]
  \begin{scope}
    \coordinate (L) at (0,0);
    \coordinate (R) at (2,0);
    \coordinate (ML) at (0.5, 0);
    \coordinate (MR) at (1.5, 0);
    \draw[->-] (ML).. controls + (0.4, 0.4) and +(-0.2, 0.4) .. (MR)
    node[above, pos =0.3] {$\deco_1$};    
    \draw[->-] (ML).. controls + (0.2, -0.4) and +(-0.4, -0.4) .. (MR)
    node[above, pos =0.7, yshift = -0.5mm] {}; 
  \end{scope}  
 \begin{scope}[yshift = 1.3cm]
    \coordinate (LB) at (0,0);
    \coordinate (RB) at (2,0);
  \end{scope}  
  \draw[double] (MR)--(R) -- (RB) -- (LB) -- (L)-- (ML);
  \draw (ML) .. controls +(0, 1) and +(0, 1) .. (MR);
\end{scope}

%% file: sch/foam_digoncup_10.tex
\begin{scope}[font=\tiny]
  \begin{scope}
    \coordinate (L) at (0,0);
    \coordinate (R) at (2,0);
    \coordinate (ML) at (0.5, 0);
    \coordinate (MR) at (1.5, 0);
    \draw[->-] (ML).. controls + (0.4, 0.4) and +(-0.2, 0.4) .. (MR) node[below, pos =0.3] {};
    \draw[->-] (ML).. controls +(0.2, -0.4) and +(-0.4, -0.4) .. (MR) node[below, pos =0.73] {$\deco_1$};
  \end{scope}  
 \begin{scope}[yshift = -1.3cm]
    \coordinate (LB) at (0,0);
    \coordinate (RB) at (2,0);
  \end{scope}  
  \draw[double] (MR) -- (R) -- (RB) -- (LB) -- (L) -- (ML);
  \draw (ML) .. controls +(0, -1) and +(0, -1) .. (MR);
\end{scope}

%% file: sch/foam_digoncup_01.tex
\begin{scope}[font=\tiny]
  \begin{scope}
    \coordinate (L) at (0,0);
    \coordinate (R) at (2,0);
    \coordinate (ML) at (0.5, 0);
    \coordinate (MR) at (1.5, 0);
    \draw[->-] (ML).. controls + (0.4, 0.4) and +(-0.2, 0.4) .. (MR) node[below, pos =0.3] {$\deco_1$};
    \draw[->-] (ML).. controls +(0.2, -0.4) and +(-0.4, -0.4) .. (MR) node[below, pos =0.75] {};
  \end{scope}  
 \begin{scope}[yshift = -1.3cm]
    \coordinate (LB) at (0,0);
    \coordinate (RB) at (2,0);

  \end{scope}  
  \draw[double] (MR) -- (R) -- (RB) -- (LB) -- (L) -- (ML);
  \draw (ML) .. controls +(0, -1) and +(0, -1) .. (MR);
\end{scope}

%% file: sch/reidem2_2.tex
\begin{scope}
    \coordinate (bl) at (-0.5, -1);
    \coordinate (br) at ( 0.5, -1);
    \coordinate (tl) at (-0.5,  1);
    \coordinate (tr) at ( 0.5,  1);
    \draw[->] (bl) -- (tl);
    \draw[->] (br) -- (tr);
\end{scope}

%% file: sch/reidem2b_1.tex
\begin{scope}[xshift = 2.5cm]
  \draw [<-](0, 0) .. controls +(1,.5) 
  .. (2, 0);
    \draw [->](0, 1) .. controls +(1,-.5) 
  .. (2, 1);
\end{scope}

%% file: sch/reidem2b_2p.tex
\begin{scope}[xshift = 2.5cm]
  \draw [<-] (0, -0.5) .. controls +(.5,0) and +(-.5, 0)
  .. (1, 0.5)  ; 
    \fill[white] (0.5, 0) circle (1mm);
    \draw(0, 0.5) .. controls +(.5,0) and +(-.5, 0)
  .. (1, -0.5) ;
  \draw[] (1, 0.5) .. controls +(.5,0) and +(-.5, 0)
  .. (2, -0.5)  node[pos=0.25, draw=\crd, fill= \csrd, circle , inner sep=1.5pt]{} node[pos=0.25, left, yshift=-2.5mm]{$\frac{1}{2}$}
  node[pos=0.7, draw=\crd, fill= \csrd, circle , inner sep=1.5pt]{} node[pos=0.7, above, xshift=2mm, yshift=-1mm]{$-\frac{1}{2}$};
    \fill[white] (1.5, 0) circle (1mm);  
      \draw[<-] (2, 0.5)  .. controls +(-.5,0) and +(.5, 0)
  .. (1, -0.5); 
\end{scope}

%% file: sch/reidem2b_4.tex
\begin{scope}
    \coordinate (bl) at (0, -1);
    \coordinate (br) at (2, -1);
    \coordinate (tl) at (0,  1);
    \coordinate (tr) at (2,  1);
    \coordinate (ml) at (-0.5,  -.8);
    \coordinate (Ml) at (-0.5,  .8);
    \coordinate (mr) at (0.5,  -.6);
    \coordinate (Mr) at (0.5,  .6);
    \draw[<-<] (bl) -- (tl);
    \draw[>->] (br) -- (tr) node[pos=0.5, draw=\crd, fill= \csrd, circle , inner sep=1.5pt]{} node[pos=0.5, right]{$-\frac{1}{2}$};
    \draw [-<] (.5, 0) arc (180:0:0.5) node[pos=0.5, draw=\crd, fill= \csrd, circle , inner sep=1.5pt]{} node[pos=0.5, above]{$\frac{1}{2}$};
    \draw (.5, 0) arc (-180:0:0.5) ;
\end{scope}

%% file: sch/movie_r2b_1.tex
\begin{scope}
    \coordinate (a) at (0.6, -0.3);
    \coordinate (b) at (0.6, 1.3);
    \draw [<-](0, 0) .. controls +(1,.5).. (2, 0);
    \draw [->](0, 1) .. controls +(1,-.5).. (2, 1);
    \draw [-<] (a) arc (180:0:0.4);
    \draw [-] (a) arc (-180:0:0.4);
    \draw [-<] (b) arc (180:0:0.4);
    \draw [-] (b) arc (-180:0:0.4);
\end{scope}

%% file: sch/movie_r2b_3.tex
\begin{scope}
    \coordinate (bl) at (-0.5, -1.5);
    \coordinate (br) at ( 0.5, -1);
    \coordinate (tr) at ( 0.5,  1);
    \coordinate (tl) at (-0.5,  1.5);
    \coordinate (bm) at (  0,-0.4);
    \coordinate (tm) at (  0, 0.4);
    \coordinate (bl2) at (1.5, -1);
    \coordinate (br2) at ( 2.5, -1.5);
    \coordinate (bm2) at (  2,-0.4);
    \coordinate (tl2) at (1.5,  1);
    \coordinate (tr2) at ( 2.5,  1.5);
    \coordinate (tm2) at ( 2, 0.4);
    \draw[<-]  (bl) .. controls +( 0, 0.5) and +(0,0) .. (bm);
    \draw[<-]  (br) .. controls +( 0, 0.5) and +(0,0) .. (bm);
    \draw[>-]  (tl) .. controls +( 0, -0.5) and +(0,0) .. (tm) coordinate[pos = 0.25] (ga) ;
    \draw[>-]  (tr) .. controls +( 0, -0.5) and +(0,0) .. (tm) coordinate[pos = 0.25] (gb) ;
    \draw [double] (bm) -- (tm);
    \draw (.5, 1) arc (180:0:0.5);
    \draw (1.5,-1) arc (180:0:-0.5);
    \draw[>-]  (bl2) .. controls +( 0, 0.5) and +(0,0) .. (bm2);
    \draw[>-]  (br2) .. controls +( 0, 0.5) and +(0,0) .. (bm2);
    \draw[<-]  (tl2) .. controls +( 0, -0.5) and +(0,0) .. (tm2);
    \draw [double] (bm2) -- (tm2);
    \draw[<-]  (tr2) .. controls +( 0, -0.5) and +(0,0) .. (tm2); 
\end{scope}

%% file: sch/reidem2b_8.tex
\begin{scope}
  \coordinate (m) at (  0,  0);
  \coordinate (t) at (  0, 0.75);
  \coordinate (br) at (+.5,  -.8);
  \coordinate (bl) at (-.5,  -.8);
  \draw[<-] (bl) .. controls +(0,0.5) and + (0, 0) .. (m) node[pos=0.3, draw=\crd, fill= \csrd, circle , inner sep=1.5pt]{} node[pos=0.3, left]{$\frac{1}{2}$};
  \draw[<-] (br) .. controls +(0,0.5) and + (0, 0) .. (m) node[pos=0.3]{} ;
  \draw[double] (m) -- (t) node[pos=0.7, draw=\crd, fill= \csrd, circle , inner sep=1.5pt]{} node[pos=0.7, right]{$-\frac{1}{2}$} 
  node[pos=0, draw=\crd, fill= white, circle , inner sep=1.5pt]{} node[pos=0, right]{$-\bar{s}$};
\end{scope}

%% file: SF_RIII.tex
\begin{prop}
    There are isomorphisms in the relative homotopy category: 
    \begin{equation}
        \label{reideIII}
        \NB{\tikz[scale=0.5]{\input{sch/reidem3_1}}} \cong \NB{\tikz[ xscale=-.5,yscale=.5]{\input{sch/reidem3_1_1}}}.
    \end{equation} 
\end{prop}

\begin{proof}
    In the following proof, standard differentials and $q$-shifts are omitted for readability. 
    The complex on the left-hand side in (\ref{reideIII}) is given by:

\begin{equation}
    \providecommand{\myxscale}{0.55}
    \providecommand{\myyscale}{0.45}
    \NB{\tikz[scale=0.5]{\input{sch/reidem3_1}}} =
    \NB{\tikz[xscale = 2.8, yscale = 2.5]{
       \node (i0) at (0, 0) { \NB{\tikz[font= \tiny, xscale=\myxscale, yscale=\myyscale]{\input{sch/reidem3_2}}} };
       \node (i1) at (-1, 0) { \NB{\tikz[font= \tiny, xscale=\myxscale, yscale=\myyscale]{\input{sch/reidem3_3}}} };
       \node (i2) at (-1, 1) { \NB{\tikz[font= \tiny, xscale=\myxscale, yscale=\myyscale]{\input{sch/reidem3_4}}} };
       \node (i3) at (-1, -1) { \NB{\tikz[font= \tiny, xscale=\myxscale, yscale=\myyscale]{\input{sch/reidem3_4}}} };
       \node (i4) at (-2, 1) { \NB{\tikz[font= \tiny, xscale=\myxscale, yscale=\myyscale]{\input{sch/reidem3_5}}} };
       \node (i5) at (-2, 0) { \NB{\tikz[font= \tiny, xscale=\myxscale, yscale=\myyscale]{\input{sch/reidem3_6}}} };
       \node (i6) at (-2, -1) { \NB{\tikz[font= \tiny, xscale=\myxscale, yscale=\myyscale]{\input{sch/reidem3_7}}} };
       \node (i7) at (-3, 0) { \NB{\tikz[font= \tiny, xscale=\myxscale, yscale=\myyscale]{\input{sch/reidem3_8}}} };
        \draw[->] (i1) -- (i0);
        \draw[->] (i2) -- (i0);
        \draw[->] (i3) -- (i0);
        \draw[->] (i4) -- (i2);
        \draw[->] (i5) -- (i2);
        \draw[->] (i5) -- (i3);
        \draw[->] (i6) -- (i3);
        \draw[->] (i4) -- ([yshift=0.05cm]i1.west);
        \draw[->] (i6) -- ([yshift=-0.05cm]i1.west);
        \draw[->] (i7) -- (i4);
        \draw[->] (i7) -- (i5);
        \draw[->] (i7) -- (i6);
    }} .
\end{equation}

This complex fits in a short exact sequence of $\Witt$-equivariant complexes that splits after forgetting the $\Witt$ action.
\begin{equation} \label{RIII1eq}
  \providecommand{\myxscale}{0.55}
  \providecommand{\myyscale}{0.45}
    \NB{\tikz[xscale = 2.5, yscale = 2.3,font=\tiny]{
    \node (i0) at (0, 0) { \NB{\tikz[font= \tiny, xscale=\myxscale, yscale=\myyscale]{\input{sch/reidem3_2}}} };
    \node (i1) at (-1, 0) { \NB{\tikz[font= \tiny, xscale=\myxscale, yscale=\myyscale]{\input{sch/reidem3_3}}} };
    \node (i2) at (-1, 1) { \NB{\tikz[font= \tiny, xscale=\myxscale, yscale=\myyscale]{\input{sch/reidem3_4}}} };
    \node (i3) at (-1, -1) { \NB{\tikz[font= \tiny, xscale=\myxscale, yscale=\myyscale]{\input{sch/reidem3_4}}} };
    \node (i4) at (-2, 1) { \NB{\tikz[font= \tiny, xscale=\myxscale, yscale=\myyscale]{\input{sch/reidem3_5}}} };
    \node (i5) at (-2, 0) { \NB{\tikz[font= \tiny, xscale=\myxscale, yscale=\myyscale]{\input{sch/reidem3_6}}} };
    \node (i6) at (-2, -1) { \NB{\tikz[font= \tiny, xscale=\myxscale, yscale=\myyscale]{\input{sch/reidem3_7}}} };
    \node (i7) at (-3, 0) { \NB{\tikz[font= \tiny, xscale=\myxscale, yscale=\myyscale]{\input{sch/reidem3_8}}} };
    \node (i8) at (-2.85, .8) {};
    \node (i9) at (-.2, -.8) {};
    \node (i10) at (-2.85, .8) {};
    \node (i11) at (-.2, .8) {};
    \node (k12) at (-3,-2.3) {\NB{\tikz[font= \tiny, xscale= \myxscale, yscale=\myyscale]{\input{sch/reidem3_9}}}};
    \node (k13) at (-2,-2.3) { \NB{\tikz[font= \tiny, xscale= \myxscale, yscale=\myyscale]{\input{sch/reidem3_9}}} };
    \node (j0) at (0, 3.5) { \NB{\tikz[font= \tiny, xscale= \myxscale, yscale=\myyscale]{\input{sch/reidem3_2}}} };
    \node (j1) at (-1, 3.5) { \NB{\tikz[font= \tiny, xscale= \myxscale, yscale=\myyscale]{\input{sch/reidem3_3}}} };
    \node (j2) at (-1, 4.5) { \NB{\tikz[font= \tiny, xscale= \myxscale, yscale=\myyscale]{\input{sch/reidem3_4}}} };
    \node (j3) at (-1, 2.5) { \NB{\tikz[font= \tiny, xscale= \myxscale, yscale=\myyscale]{\input{sch/reidem3_4}}} };
    \node (j4) at (-2, 4.5) { \NB{\tikz[font= \tiny, xscale= \myxscale, yscale=\myyscale]{\input{sch/reidem3_5}}} };
    \node (j5) at (-2, 3.5) { \NB{\tikz[font= \tiny, xscale= \myxscale, yscale=\myyscale]{\input{sch/reidem3_4}}} };
    \node (j6) at (-2, 2.5) { \NB{\tikz[font= \tiny, xscale= \myxscale, yscale=\myyscale]{\input{sch/reidem3_7}}} };
    \node (j7) at (-3, 3.5) { \NB{\tikz[font= \tiny, xscale= \myxscale, yscale=\myyscale]{\input{sch/reidem3_10}}} };
    \node (j8) at (-2.85, 2.75) {};
    \node (j9) at (-.2, 2.75) {};
    \node (j10) at (-2.85, 4.25) {};
    \node (j11) at (-.2, 4.25) {};
    \draw[->] (i1) -- (i0);
    \draw[->] (i2) -- (i0);
    \draw[->] (i3) -- (i0);
    \draw[->] (i4) -- ([yshift=0.05cm]i1.west);
    \draw[->] (i6) -- ([yshift=-0.05cm]i1.west);
    \draw[->] (i4) -- (i2);
    \draw[->] (i5) -- (i2);
    \draw[->] (i5) -- (i3);
    \draw[->] (i6) -- (i3);
    \draw[->] (i7) -- (i4);
    \draw[->] (i7) -- (i5);
    \draw[->] (i7) -- (i6);
    \draw[->] (k12) --(k13) node[pos=0.5, above] {\Id};
    \draw[->] (i7) -- (k12) node[pos=0.5, right] {$ \mapdcap \circ  \mapunzip $};
    \draw[->] (i6) -- (k13)  node [right, pos = 0.5] {$\begin{pmatrix} 0  & \mapdcap & 0 \end{pmatrix} $};
    \node[draw, densely dotted ,fit=(i8) (i9) (i10) (i11) ,inner sep=10ex,rectangle,rounded corners] {};
    \node[draw, densely dotted ,fit=(j8) (j9) (j10) (j11) ,inner sep=10ex,rectangle,rounded corners] {};
    \node[draw, densely dotted ,fit=(k12) (k13) ,inner sep=1ex,rectangle,rounded corners] {};
    \draw[->] (j1) -- (j0);
    \draw[->] (j2) -- (j0);
    \draw[->] (j3) -- (j0);
    \draw[->] (j4) -- (j1);
    \draw[->] (j4) -- (j2);
    \draw[->] (j5) -- (j2);
    \draw[->] (j5) -- (j3);
    \draw[->] (j6) -- (j3);
    \draw[->] (j6) -- (j1);
    \draw[->] (j7) -- (j4);
    \draw[->] (j7) -- (j5);
    \draw[->] (j7) -- (j6);
    \draw[->] (j7) -- (i7) node[pos=0.54, left] {$\Upsilon$};
    \draw[->] (j6) -- (i4)  node [left, pos = 0.5] {$\begin{pmatrix} \Id & 0 & 0 \\ 0 & \mapdcup & 0\\ 0 & 0 & \Id \end{pmatrix} $};
    \draw[->] (j3) -- (i2)  node [right, pos = 0.5] {$\begin{pmatrix} \Id & 0 & 0 \\ 0 & \Id & 0 \\ 0 & 0 & \Id \end{pmatrix} $};
    \draw[->] (j0) -- (i0)  node [right, pos = 0.5] {$\Id$};
  }}
\end{equation}

The thickest edge of the web at the top left is of thickness 3.
The map \(\Upsilon := \mapsingsdle \circ \mapasso \circ \mapdcup \) is the composition of a digon-cup, associativity, and the singular saddle. Its movie is: 
\[\mymoviefour[scale=0.75]{\NB{\tikz[xscale=.6,yscale=.5]{\input{sch/movie_r3_1}}}}{\NB{\tikz[xscale=.6,yscale=.5]{\input{sch/movie_r3_2}}}}{\NB{\tikz[xscale=.6,yscale=.5]{\input{sch/movie_r3_3}}}}{ \NB{\tikz[xscale=.6,yscale=.5]{\input{sch/movie_r3_4}}}}.\]
Since the bottom complex is acyclic, the two other complexes are isomorphic in the relative homotopy category.

There is a second short exact sequence of $\Witt$-equivariant complexes that splits when one forgets the $\Witt$ action.
\begin{equation} \label{RIII3eq}
  \providecommand{\myxscale}{0.55}
  \providecommand{\myyscale}{0.45}
    \NB{\tikz[xscale = 2.5, yscale = 2.3,font=\tiny]{
    \node (i0) at (0, 0.5) { \NB{\tikz[font= \tiny, xscale=\myxscale, yscale=\myyscale]{\input{sch/reidem3_2}}} };
    \node (i1) at (-1, 0) { \NB{\tikz[font= \tiny, xscale=\myxscale, yscale=\myyscale]{\input{sch/reidem3_3}}} };
    \node (i2) at (-1, 1) { \NB{\tikz[font= \tiny, xscale=\myxscale, yscale=\myyscale]{\input{sch/reidem3_4}}} };
    \node (i4) at (-2, 1) { \NB{\tikz[font= \tiny, xscale=\myxscale, yscale=\myyscale]{\input{sch/reidem3_5}}} };
    \node (i6) at (-2, 0) { \NB{\tikz[font= \tiny, xscale=\myxscale, yscale=\myyscale]{\input{sch/reidem3_7}}} };
    \node (i7) at (-3, 0.5) { \NB{\tikz[font= \tiny, xscale=\myxscale, yscale=\myyscale]{\input{sch/reidem3_10}}} };
    \node (i8) at (-2.8, 0.2) {};
    \node (i9) at (-.2, 0.2) {};
    \node (i10) at (-2.8, .8) {};
    \node (i11) at (-.2, .8) {};
    \draw[->] (i1) -- (i0);
    \draw[->] (i2) -- (i0);
    \draw[->] (i4) -- (i1);
    \draw[->] (i6) -- (i1);
    \draw[->] (i4) -- (i2);
    \draw[->] (i6) -- (i2);
    \draw[->] (i7) -- (i4);
    \draw[->] (i7) -- (i6);
    \node[draw, densely dotted ,fit=(i8) (i9) (i10) (i11) ,inner sep=10ex,rectangle,rounded corners] {};
    \node (j0) at (0, 3.5) { \NB{\tikz[font= \tiny, xscale= \myxscale, yscale=\myyscale]{\input{sch/reidem3_2}}} }; 
    \node (j1) at (-1, 3.5) { \NB{\tikz[font= \tiny, xscale= \myxscale, yscale=\myyscale]{\input{sch/reidem3_3}}} };
    \node (j2) at (-1, 4.5) { \NB{\tikz[font= \tiny, xscale= \myxscale, yscale=\myyscale]{\input{sch/reidem3_4}}} };
    \node (j3) at (-1, 2.5) { \NB{\tikz[font= \tiny, xscale= \myxscale, yscale=\myyscale]{\input{sch/reidem3_4}}} };
    \node (j4) at (-2, 4.5) { \NB{\tikz[font= \tiny, xscale= \myxscale, yscale=\myyscale]{\input{sch/reidem3_5}}} };
    \node (j5) at (-2, 3.5) { \NB{\tikz[font= \tiny, xscale= \myxscale, yscale=\myyscale]{\input{sch/reidem3_4}}} };
    \node (j6) at (-2, 2.5) { \NB{\tikz[font= \tiny, xscale= \myxscale, yscale=\myyscale]{\input{sch/reidem3_7}}} };
    \node (j7) at (-3, 3.5) { \NB{\tikz[font= \tiny, xscale= \myxscale, yscale=\myyscale]{\input{sch/reidem3_10}}} };
    \node (j8) at (-2.8, 2.75) {};
    \node (j9) at (-.2, 2.75) {};
    \node (j10) at (-2.8, 4.25) {};
    \node (j11) at (-.2, 4.25) {};
    \node[draw, densely dotted ,fit=(j8) (j9) (j10) (j11) ,inner sep=10ex,rectangle,rounded corners] {};
    \draw[->] (j1) -- (j0);
    \draw[->] (j2) -- (j0);
    \draw[->] (j3) -- (j0);
    \draw[->] (j4) -- (j1);
    \draw[->] (j6) -- (j1);
    \draw[->] (j4) -- (j2);
    \draw[->] (j5) -- (j3);
    \draw[->] (j5) -- (j2);
    \draw[->] (j6) -- (j3);
    \draw[->] (j7) -- (j4);
    \draw[->] (j7) -- (j6);
    \draw[->] (j7) -- (j5);
    \draw[->] (j7) -- (i7) node[pos=0.6, left] {$\Id$};
    \draw[->] (j6) -- (i4)  node [left, pos = 0.5] {$\begin{pmatrix} \Id & 0 & 0 \\0 & 0 & \Id \end{pmatrix} $};
    \draw[->] (j3) -- (i2)  node [right, pos = 0.5] {$\begin{pmatrix} \Id & 0 & \Id \\0 & \Id & 0 \end{pmatrix} $};
    \draw[->] (j0) -- (i0)  node [right, pos = 0.58] {$\Id$};
    \node (k0) at (-1, 5.8) { \NB{\tikz[font= \tiny, xscale= \myxscale, yscale=\myyscale]{\input{sch/reidem3_4}}} };
    \node (k1) at (-2, 5.8) { \NB{\tikz[font= \tiny, xscale= \myxscale, yscale=\myyscale]{\input{sch/reidem3_4}}} };
    \draw[->] (k0) -- (j2)  node [right, pos = 0.4] {$\begin{pmatrix}  -\Id \\0 \\ \Id \end{pmatrix} $};
    \draw[->] (k1) -- (j4) node [right, pos = 0.5] {$\begin{pmatrix}  0 \\ \Id \\ 0 \end{pmatrix} $};
    \draw[->] (k1) -- (k0) node [above, pos = 0.5] {$\Id$};
    \node[draw, densely dotted ,fit=(k0) (k1) ,inner sep=1ex,rectangle,rounded corners] {};
      }}.
\end{equation}

By the same argument as in (\ref{RIII1eq}), $\NB{\tikz[scale=0.4]{\input{sch/reidem3_1}}}$ is isomorphic to the following complex
\begin{equation} \label{RIII3eq1}
  \providecommand{\myxscale}{0.55}
  \providecommand{\myyscale}{0.45}
    \NB{\tikz[xscale = 2.5, yscale = 2.5]{
    \node (i0) at (0, 0.5) { \NB{\tikz[font= \tiny, xscale=\myxscale, yscale=\myyscale]{\input{sch/reidem3_2}}} };
    \node (i1) at (-1, 0) { \NB{\tikz[font= \tiny, xscale=\myxscale, yscale=\myyscale]{\input{sch/reidem3_3}}} };
    \node (i2) at (-1, 1) { \NB{\tikz[font= \tiny, xscale=\myxscale, yscale=\myyscale]{\input{sch/reidem3_4}}} };
    \node (i4) at (-2, 1) { \NB{\tikz[font= \tiny, xscale=\myxscale, yscale=\myyscale]{\input{sch/reidem3_5}}} };
    \node (i6) at (-2, 0) { \NB{\tikz[font= \tiny, xscale=\myxscale, yscale=\myyscale]{\input{sch/reidem3_7}}} };
    \node (i7) at (-3, 0.5) { \NB{\tikz[font= \tiny, xscale=\myxscale, yscale=\myyscale]{\input{sch/reidem3_10}}} };
    \draw[->] (i1) -- (i0);
    \draw[->] (i2) -- (i0);
    \draw[->] (i4) -- (i1);
    \draw[->] (i6) -- (i1);
    \draw[->] (i4) -- (i2);
    \draw[->] (i6) -- (i2);
    \draw[->] (i7) -- (i4);
    \draw[->] (i7) -- (i6);
      }}.
\end{equation}
By the isomorphism
\begin{equation}
    \NB{\tikz[font= \tiny, scale=0.5]{\input{sch/reidem3_10}}} \simeq \NB{\tikz[font= \tiny, scale=0.5,xscale=-1]{\input{sch/reidem3_11}}},
\end{equation}
the complex (\ref{RIII3eq1}) is left-right symmetric making $\NB{\tikz[scale=0.4]{\input{sch/reidem3_1}}}$ and $\NB{\tikz[scale=0.4]{\input{sch/reidem3_1_1}}}$ isomorphic to the same complex in the relative homotopy category.
\end{proof}

%% file: sch/reidem3_2.tex
\begin{scope}
    \coordinate (bl) at (-0.5, -1);
    \coordinate (br) at ( 0.5, -1);
    \coordinate (tl) at (-0.5,  1);
    \coordinate (tr) at ( 0.5,  1);
    \coordinate (ml) at (-0.5,  -.8);
    \coordinate (Ml) at (-0.5,  .8);
    \coordinate (mr) at (0.5,  -.6);
    \coordinate (Mr) at (0.5,  .6);
    \coordinate (brr) at ( 1.5, -1);
    \coordinate (trr) at ( 1.5, 1);
    \draw[->-] (bl) -- (tl);
    \draw[->-] (br) -- (tr); 
    \draw[->-] (brr) -- (trr);
\end{scope}

%% file: sch/movie_r3_1.tex
\begin{scope}
    \coordinate (au) at (-.5,1.3);
    \coordinate (bu) at (.5,1.3);
    \coordinate (cu) at (1.5,1.3); 
    \coordinate (abu) at (0,0.5);
    \coordinate (abcu) at (.5,0);
    \draw[<-] (au) .. controls +( 0, -0.5) and +(0,0) .. (abu);
    \draw[<-]  (bu) .. controls +( 0, -0.5) and +(0,0) .. (abu);
    \draw[double] (abu) .. controls +( 0, -0.5) and +(0,0) .. (abcu); 
    \draw[<-] (cu) .. controls +( 0, -0.7) and +(0,0) .. (abcu);
    \coordinate (ad) at (-.5,-2.3);
    \coordinate (bd) at (.5,-2.3);
    \coordinate (cd) at (1.5,-2.3); 
    \coordinate (abd) at (0,-1.5);
    \coordinate (abcd) at (.5,-1);
    \draw[>-] (ad) .. controls +( 0,0.5) and +(0,0) .. (abd);
    \draw[>-]  (bd) .. controls +( 0,0.5) and +(0,0) .. (abd);
    \draw[double] (abd) .. controls +( 0,0.5) and +(0,0) .. (abcd); 
    \draw[>-] (cd) .. controls +( 0,0.7) and +(0,0) .. (abcd);
    \draw[line width=.7mm] (abcd) -- (abcu) coordinate[pos = 0.5] (ga);
\end{scope}

%% file: sch/movie_r3_2.tex
\begin{scope}
    \coordinate (au) at (-.5,1.3);
    \coordinate (bu) at (.5,1.3);
    \coordinate (cu) at (1.5,1.3); 
    \coordinate (abu) at (0,0.7);
    \coordinate (abup) at (0.07,.5);
    \coordinate (abupp) at (0.25,.2);
    \coordinate (abcu) at (.5,0);
    \draw[<-] (au) .. controls +( 0, -0.5) and +(0,0) .. (abu);
    \draw[<-]  (bu) .. controls +( 0, -0.5) and +(0,0) .. (abu);
    \draw[double] (abu) .. controls +( 0, -0.1) and +(0,0) .. (abup); 
    \draw[double] (abupp) .. controls +( 0,0) and +(0,-0.1) .. (abcu); 
    \draw[] (abup) .. controls +(0.4,0.08) and +(0,0) .. (abupp);
    \draw[] (abup) .. controls +(-.2,-0.4) and +(0,0) .. (abupp);
    \draw[<-] (cu) .. controls +( 0, -0.7) and +(0,0) .. (abcu);
    \coordinate (ad) at (-.5,-2.3);
    \coordinate (bd) at (.5,-2.3);
    \coordinate (cd) at (1.5,-2.3); 
    \coordinate (abd) at (0,-1.5);
    \coordinate (abcd) at (.5,-1);
    \draw[>-] (ad) .. controls +( 0,0.5) and +(0,0) .. (abd);
    \draw[>-]  (bd) .. controls +( 0,0.5) and +(0,0) .. (abd);
    \draw[double] (abd) .. controls +( 0,0.5) and +(0,0) .. (abcd); 
    \draw[>-] (cd) .. controls +( 0,0.7) and +(0,0) .. (abcd);
    \draw[line width=.7mm] (abcd) -- (abcu) coordinate[pos = 0.5] (ga);
\end{scope}

%% file: sch/movie_r3_3.tex
\begin{scope}
    \coordinate (au) at (-.5,1.3);
    \coordinate (bu) at (.5,1.3);
    \coordinate (cu) at (1.5,1.3); 
    \coordinate (abu) at (0,0.7);
    \coordinate (abup) at (0.07,.5);
    \coordinate (bcu) at (0.95,.3);
    \coordinate (abcu) at (.5,-.1);
    \draw[<-] (au) .. controls +( 0, -0.5) and +(0,0) .. (abu);
    \draw[<-]  (bu) .. controls +( 0, -0.5) and +(0,0) .. (abu);
    \draw[-<-] (abup) .. controls +(0.4,0.08) and +(0,0) .. (bcu);
    \draw[-<-] (abup) .. controls +(-.15,-0.6) and +(0,0) .. (abcu);
    \draw[<-] (cu) .. controls +( 0, -0.7) and +(0,0) .. (bcu);
    \draw[double] (abu) .. controls +( 0, -0.1) and +(0,0) .. (abup); 
    \draw[double] (bcu) .. controls +(-0.1,-0.3) and +(0,0) .. (abcu); 
    \coordinate (ad) at (-.5,-2.3);
    \coordinate (bd) at (.5,-2.3);
    \coordinate (cd) at (1.5,-2.3); 
    \coordinate (abd) at (0,-1.5);
    \coordinate (abcd) at (.5,-1);
    \draw[>-] (ad) .. controls +( 0,0.5) and +(0,0) .. (abd);
    \draw[>-]  (bd) .. controls +( 0,0.5) and +(0,0) .. (abd);
    \draw[double] (abd) .. controls +( 0,0.5) and +(0,0) .. (abcd); 
    \draw[>-] (cd) .. controls +( 0,0.7) and +(0,0) .. (abcd);
    \draw[line width=.7mm] (abcd) -- (abcu) coordinate[pos = 0.5] (ga);
\end{scope}

%% file: sch/movie_r3_4.tex
 \begin{scope}
    \coordinate (au) at (-.5,1.3);
    \coordinate (bu) at (.5,1.3);
    \coordinate (cu) at (1.5,1.3); 
    \coordinate (abu) at (0,0.6);
    \coordinate (abup) at (0,.15);
    \coordinate (bcu) at (0.95,-.1);
    \coordinate (abcu) at (.95,-.7);
    \draw[<-] (au) .. controls +( 0, -0.5) and +(0,0) .. (abu);
    \draw[<-]  (bu) .. controls +( 0, -0.5) and +(0,0) .. (abu);
    \draw[-<-] (abup) -- (bcu);
    \draw[<-] (cu) .. controls +( 0, -0.7) and +(0,0) .. (bcu);
    \draw[double] (abu) .. controls +( 0, -0.1) and +(0,0) .. (abup); 
    \draw[double] (bcu) .. controls +(0,0) and +(0,0) .. (abcu); 
    \coordinate (ad) at (-.5,-2.3);
    \coordinate (bd) at (.5,-2.3);
    \coordinate (cd) at (1.5,-2.3); 
    \coordinate (abd) at (0,-1.5);
    \coordinate (abcd) at (0,-1);
    \draw[>-] (ad) .. controls +( 0,0.5) and +(0,0) .. (abd);
    \draw[>-]  (bd) .. controls +( 0,0.5) and +(0,0) .. (abd);
    \draw[double] (abd) .. controls +( 0,0.5) and +(0,0) .. (abcd); 
    \draw[>-] (cd) .. controls +( 0,0.7) and +(0,0) .. (abcu);
    \draw[->-] (abcd) -- (abcu);
    \draw[-<-] (abup) -- (abcd);
\end{scope}

%% file: SF_unfrinv.tex
\label{unfram}
The homology $\KRW_N(~.~)$ is an invariant of framed link diagrams, however due to the grading shifts and twisting introduced by Reidemeister I moves, it is not naturally an unframed invariant.
The braiding complexes can be modified slightly to correct this issue.
\begin{equation}
  \NB{\tikz[xscale = 0.6]{\input{sch/positive_crossing}}} :=
\NB{\tikz[xscale = 3.5, yscale = 3]
{\node (i0) at (0.2, 0){$q^{N-1} \NB{\tikz[font= \tiny,scale=0.6]{\input{sch/unframed_1}}}$};
\node (i1) at (-1, 0) {$q^{N} \, \NB{\tikz[font= \tiny,scale=0.6]{\input{sch/unframed_2}}}$ };
\draw[->] (i1) -- (i0) coordinate[pos=0.45] (a);
\node[above] at (a) {\NB{\tikz[font=\tiny, scale=.5]{\input{sch/foam_unzip_empty}}}};
  }}
\end{equation}

\begin{equation}
   \NB{\tikz[xscale =- 0.6]{\input{sch/positive_crossing}}}:=
    \NB{ \tikz[xscale = 3.5, yscale = 3]{
    \node (i0) at (-1, 0) {$q^{1-N} \NB{\tikz[font= \tiny, scale=0.6]{\input{sch/unframed_3}}}$};
    \node (i1) at (0.2, 0) {$q^{-N}\NB{\tikz[font= \tiny, scale=0.6]{\input{sch/unframed_4}}}$ };
    \draw[->] (i0) -- (i1) coordinate[pos=0.6] (b);
    \node[above] at (b) {$\NB{\tikz[font=\tiny, scale=.5]{\input{sch/foam_zip_empty}}}$}; }} \,
\end{equation}
For a link \( L \), let \( {\uKRW}_N(L,\scalars)\) be the unframed ($\Witt$-equivariant) Khovanov--Rozansky $\gl_N$-homology, that is the homology of the cube of resolution defined from the adjusted braiding complexes.
\begin{cor}
    The unframed Khovanov--Rozansky $\gl_N$-homology ${\uKRW}_N(L,\scalars)$ is an invariant of oriented links.
\end{cor}


%% file: sch/unframed_1.tex
\begin{scope}
   \coordinate (bl) at (-0.5, -1);
    \coordinate (br) at ( 0.5, -1);
    \coordinate (tl) at (-0.5,  1);
    \coordinate (tr) at ( 0.5,  1);
    \draw[->] (bl) -- (tl)  node[pos=0.2, draw=\crd, fill= \csrd, circle , inner sep=1.5pt]{} node[pos=0.2, right]{$\frac{1}{2}$};
    \draw[->] (br) -- (tr);
    \draw[decorate,decoration={snake}] (-1.08,-1.4) -- (1.08,-1.4);
\end{scope}

%% file: sch/unframed_3.tex
\begin{scope}
   \coordinate (bl) at (-0.5, -1);
    \coordinate (br) at ( 0.5, -1);
    \coordinate (tl) at (-0.5,  1);
    \coordinate (tr) at ( 0.5,  1);
    \draw[->] (bl) -- (tl)  node[pos=0.2, draw=\crd, fill= \csrd, circle , inner sep=1.5pt]{} node[pos=0.2, right]{$-\frac{1}{2}$};
    \draw[->] (br) -- (tr);
    \draw[decorate,decoration={snake}] (-1.08,-1.4) -- (1.08,-1.4);
\end{scope}

%% file: SF_newfunctoriality.tex

The functoriality of \(\mathfrak{gl}_N\) homology implies that link cobordisms induce \(\scalars_N\)-linear maps between the homologies of the source and target links \cite{ETW17}.
However, there is no \textit{a priori} reason for the induced maps to be \(\Witt\)-equivariant.
In fact, the non-triviality of the action on cups, caps and saddles implies that induced maps will rarely be equivariant.
Fortunately, the action on these building blocks is well-behaved.

\begin{dfn}
Let \(\rLinks\) be the category of links and cobordisms with scalar-labelled red dots on the components of links.
Let \(\mathrm{Mod}(\scalars_N)_{g,p}^{\Witt}\) denote the category of bigraded, projective, \((\Witt \# \scalars_N)\)-modules with \(\scalars_N\)-linear morphisms that are not necessarily \(\Witt\)-equivariant or homogeneous.
The homology of a red dotted link is defined by adding solid red dot twists to the webs in the cube of resolutions.
\end{dfn}

\begin{thm}
There is a well-defined functor \(\uKRW_N: \rLinks \to \mathrm{Mod}(\scalars_N)^{\Witt}_{g,p}\).
\end{thm}
\begin{proof}
This follows immediately from \cite[Theorem 4.5]{ETW17} and the isotopy invariance of the \(\Witt\)-action on foams.
\end{proof}

\begin{thm}
Let \(K_1\) and \(K_2\) be knots decorated with red dots labelled \(\lambda_1\) and \(\lambda_2\) respectively and let
\[C: K_1^{\lambda_1} \to K_2^{\lambda_2}\] be a connected cobordism containing \(k\) Reidemeister I moves counted with sign.
If \(\lambda_2 - \lambda_1 = (\chi(C)-k)/2\), then \[\KRW_N(C): \KRW_N(K_1^{\lambda_1}) \to \KRW_N(K_2^{\lambda_2})\] is \(\Witt\)-equivariant.
If \(\lambda_2 - \lambda_1 = \chi(C)/2\), then \[\uKRW_N(C): \uKRW_N(K_1^{\lambda_1}) \to \uKRW_N(K_2^{\lambda_2})\] is \(\Witt\)-equivariant.
\end{thm}
\begin{proof}
Every knot cobordism is a composition of cups, caps, saddles and Reidemeister moves.
The maps associated to Reidemeister moves II and III are equivariant and Reidemeister I moves are equivariant up to a twist with positive and negative requiring the opposite correction.
Examining the action on cups, caps, and saddles implies that the action depends only on the Euler characteristic and comparing with the correction for Reidemeister I moves shows that the required twists for both are precisely the solid red dot.
\end{proof}

\begin{cor}
Let \(L_1\) and \(L_2\) be links decorated with collections of red dots on each component labelled \(\mathbf{\lambda}_1 = (\lambda^i_1)\) and \(\mathbf{\lambda}_2= (\lambda^j_2)\) respectively and let
\[C: L_1^{\mathbf{\lambda}_1} \to L_2^{\mathbf{\lambda}_2}\] be a connected cobordism containing \(k\) Reidemeister I moves counted with sign.
If \(\sum_i \mathbf{\lambda}^i_2 - \sum_j \mathbf{\lambda}^j_1 = (\chi(C) - k)/2\), then \[\KRW_N(C): \KRW_N(L_1^{\mathbf{\lambda}_1}) \to \KRW_N(L_2^{\mathbf{\lambda}_2})\] is \(\Witt\)-equivariant.
If \(\sum_i \mathbf{\lambda}^i_2 - \sum_j \mathbf{\lambda}^j_1 = \chi(C)/2\), then \[\uKRW_N(C): \uKRW_N(L_1^{\mathbf{\lambda}_1}) \to \uKRW_N(L_2^{\mathbf{\lambda}_2})\] is \(\Witt\)-equivariant.
\end{cor}

An analogous statement for cobordisms with multiple components follows easily.

%% file: SF_newcompute.tex
In this section, we show that good coordinates can be given to the state spaces of simple foams by shifting the algebra of decorations by the inverse of a formal square root of a special polynomial.
The action of \(\wgen_m\) is extended using the power and chain rules:
\[
\wgen_m \cdot f^{-1/2} := \frac{1}{2} f^{-3/2} (\wgen_m \cdot f).
\]

In a previous paper by the second author, the action resricted to \(\sld\) was analyzed using this presentation \cite{Roz23}.
The algebra of decorations is decomposed as a \(\Witt\)-module, but a general analysis of the twisted modules is beyond the scope of the paper.

\subsection{Generalized theta webs}

Let \(\mathbf{a} = (a_1, \dots, a_n)\) be an integer composition such that \(\sum_i a_i = N\). For each \(k\) let \(d_k : = a_1+ \cdots +a_k\) be the \(k\)th partial sum.
Define
\[t_{\mathbf{a}} := 
\prod_{k=1}^n \prod_{i=d_{k-1} +1}^{ d_k} \prod_{j=d_k +1}^{ N} (x_i-x_j).\]

\begin{dfn}
Let \(\theta_{\mathbf{a}}\) and \(\Theta_{\mathbf{a}} \in \mathcal{F}_N(\theta_{\mathbf{a}})\) be the generalized theta web and foam associated to this composition, respectively: 
\end{dfn}

\begin{figure}[h!]
	\centering
	\NB{\tikz[scale=0.85,font=\tiny]{\input{sch/theta_gen_web}}} \hfill \NB{\tikz[scale=0.60,font=\tiny]{\input{sch/Theta_gen_foam}}}
	\caption{Theta web and foam}
\end{figure}

\begin{prop}
\label{twisted-hom}
If \(\lambda=\mu=0\) and \(s = 1/2\), then
as \(\Witt \# \scalars_N\)-modules,
\[
\mathcal{F}_N(\theta_{\mathbf{a}}) \cong t_{\mathbf{a}}^{-1/2} \scalars_{\mathbf{a}}.
\]
\end{prop}
\begin{proof}
The state space of \(\theta_{\mathbf{a}}\) is 
a rank one, free \(\scalars_{\mathbf{a}}\)-module generated by \(\Theta_{\mathbf{a}}\), where \(\scalars_{\mathbf{a}}\) acts as the algebra of decorations.
The proposition is proved by comparing the \(\Witt\)-action on \(\Theta_{\mathbf{a}}\) with the action on \(t_{\mathbf{a}}^{-1/2}\).

The foam \(\Theta_{\mathbf{a}}\) can be decomposed starting from the far right edge as a sequence of concentric cups, one for each facet of thickness \(a_k\), and zips, one for each pair of thin facets with thicknesses \(a_k\) and \(d_{k-1}\).
After setting \(\lambda=\mu=0\),
the action of \(\wgen_m\) on the zip with thin facets \(a_k\) and \(d_{k-1}\)
is multiplication by
\begin{align}
\label{crosstermtwist}
&-\overline{s}  \sum_{i=0}^{m} p_i(x_1, \dots, x_{d_{k-1}}) p_{m-i}(x_{d_{k-1}+1}, \dots, x_{d_k}) \\
&= -\overline{s}  \sum_{i=d_{k-1}+1}^{d_k}\sum_{j=1}^{d_{k-1}} h_m(x_i,x_j).
\end{align}
The action on the cup of thickness \(a_k\) is multiplication by
\begin{align}
\label{zip-crosstermtwist}
&\frac{1}{2}  \sum_{i=0}^{m} p_i(x_{d_{k-1}+1}, \dots, x_{d_{k}}) \left(
p_{m-i}(x_{1}, \dots, x_{d_{k-1}})+
p_{m-i}(x_{d_{k}+1}, \dots, x_N)
\right)\\
&= 
\frac{1}{2} \sum_{i = d_{k-1} +1}^{d_k} \left( \sum_{j =1 }^{d_{k-1}} h_m(x_i,x_j)
+ \sum_{j=d_k+1}^{N} h_m(x_i,x_j) \right).
\end{align}
If \(s=1/2\), the first term of the cup contribution cancels with the zip contribution.
The action of \(\wgen_m\) on all zips and cups together is multiplication by
\begin{equation}
\frac{1}{2} \sum_{k=1}^n \sum_{i=d_{k-1} +1}^{ d_k} \sum_{j=d_k +1}^{ N} h_m(x_i,x_j).
\end{equation}

On the other hand, 
\begin{equation}
\wgen_m \cdot (x_i - x_j)
= - (x_i^{m+1} - x_j^{m+1}) = - h_m(x_i,x_j) (x_i-x_j).
\end{equation}
Thus the action on \(t_{\mathbf{a}}^{-1/2}\) is
\begin{align}
\wgen_m \cdot t_{\mathbf{a}}^{-1/2} &=
\wgen_m \cdot \left[ \prod_{k=1}^n \prod_{i=d_{k-1} +1}^{ d_k} \prod_{j=d_k +1}^{ N} (x_i-x_j)^{-1/2} \right]\\
&= \left[\frac{1}{2} \sum_{k=1}^n \sum_{i=d_{k-1} +1}^{ d_k} \sum_{j=d_k +1}^{ N} h_m(x_i,x_j) \right] t_{\mathbf{a}}^{-1/2},
\end{align}
which agrees with the action on \(\Theta_{\mathbf{a}}\).
\end{proof}

\begin{rmk}
The leading term \(t_{\mathbf{a}}^{-1/2}\) is the value of the foam evaluation formula on the non-closed foam \(\Theta_{\mathbf{a}}\).
Even though foam evaluation is well defined only for closed foams, the \(\Witt\)-action was defined so that it would commute with the local contributions from non-closed foams to the formula.
\end{rmk}

\begin{rmk}
The condition on the parameters can be removed by adding red dot twists to \(\theta_{\mathbf{a}}\) in appropriate places. At each of the below vertices of the theta web, the following twists cancel out the contribution of the \(\lambda, \mu,\) and \(s\) parameters:
\[
	\centering
	\NB{\tikz[scale=0.85,font=\tiny]{\input{sch/setparam}}}.
\]
Additional solid red dots labeled \(\lambda\) at these vertices add twists to the polynomial ring of the form
\[\prod_{i=d_{k-1} +1}^{ d_k} \prod_{j=d_k +1}^{ N} (x_i-x_j)^{\lambda}. \]
\end{rmk}



\begin{cor}
Let \(\varnothing\), \(U\), and \(\Theta\) denote the empty web, unknot, and standard theta web with thickness one thin facets. As \(\scalars_N \# \Witt\)-modules:
\begin{itemize}
\item \(\mathcal{F}_N(\varnothing) \cong \mathcal{F}_N(\theta_{N}) \cong \scalars_N\)
\item \(\mathcal{F}_N(U) \cong \mathcal{F}_N(\theta_{1,N-1}) \cong \prod_{i > 1} (x_1 - x_j)^{-1/2} \scalars_{1,N-1}\).
\item \(\mathcal{F}_N(\theta) \cong \mathcal{F}_N(\theta_{1,1,N-2}) \cong \prod_{i > 1} (x_1 - x_i)^{-1/2}\prod_{j > 2} (x_2 - x_j)^{-1/2} \scalars_{1,1,N-2}\).
\item \(\mathcal{F}_N(\theta_{1^N}) \cong \prod_{i > j} (x_i - x_j)^{-1/2}\scalars_{1^N}\).
\end{itemize}
\end{cor}
\begin{proof}
Examining equations (\ref{crosstermtwist}) and (\ref{zip-crosstermtwist}) when \(k=n\) shows that the facet with thickness \(N\) can be removed \(\Witt\)-equivariantly from \(\Theta_{\mathbf{a}}\).
Thus the webs \(\varnothing, U\) and \(\theta\) can be realized as special cases of the generalized theta web described above.
\end{proof}

\subsection{The algebra of decorations as a \(\Witt\)-module}

In the following, we additionally require \(\scalars\) to be a field.

Recall that the action of \(\Witt\) on \(\scalars_{1^N}\) commutes with the symmetric group action permuting variables.
Therefore it is sufficient to decompose the full algebra of polynomials and then apply projections onto the symmetric subrings.
Our strategy is to find generators for \(\scalars_{1^N}\) as a \(\Wittt_{1}\)-module using Cartan-type harmonic polynomials 
then replace them with highest weight vectors to upgrade to a \(\Witt\)-decomposition.

Harmonic polynomials were first employed to find \(\Wittt_{1}\)-generators by \cite{Nishiyama_1996}.
They were used to produce a full \(\Witt\)-decomposition for \(\scalars_{1^2}\) in characteristic zero by
\cite{Tanaka_2001} and in positive characteristic by 
\cite{ChangYao_2021}.
An \(\sld\)-decomposition of \(\scalars_{1^N}\) was previously found using standard theorems about nilpotent derivations \cite{bedratyuk_2024}.
The contribution of this section is to combine these tools and provide more information on the basis of harmonic polynomials using the commuting \(S_N\)-action.

\begin{dfn}[{\cite[Chapter III]{Helgason_2000}}]
\label{inner-prod}
Let \(\text{Diff}(\scalars_{1^N})\) denote the algebra of polynomial differential operators under composition.
The map \(x_i \mapsto \frac{\partial}{\partial x_i}\) extends to an algebra homomorphism \(\partial_{(-)}: \scalars_{1^N} \to \text{Diff}(\scalars_{1^N})\)
and the pairing \(\langle p, q \rangle := (\partial_p q)(0)\) defines the \textit{Fock inner product} on \(\scalars_{1^N}\).
\end{dfn}

\begin{dfn}
For any \(m \geq -1\), let
\[
\wgen^*_m := - \sum_{i=1}^N x_i \frac{\partial^{m+1}}{\partial x_i^{m+1}},
\]
These operators satisfy the relation \([\wgen_n^*, \wgen_m^*] = (m-n)\wgen_{n+m}^*\).
Let \(\Wittt_{1}^*\) denote the Lie algebra generated by \(\wgen_m^*\) for \(m \geq 1\).
\end{dfn}

\begin{lem}[{\cite[Section 4]{Nishiyama_1996}}]
\label{adjoints}
The operators \(\wgen_m\) and \(\wgen_m^*\) are adjoint with respect to the Fock inner product.
The subspace \({H_N := \bigcap_{m \geq 1} \mathrm{ker} \wgen^*_m}\) is finite-dimensional, orthogonal to \(\Wittt_1 (\scalars_{1^N})\),
and a minimal generating subspace for \(\scalars_{1^N}\) as a \(\Wittt_{1}\)-module.
\end{lem}

\begin{rmk}
The elements of \(H_N\) are called \textit{Cartan-harmonic polynomials}.
With respect to the same inner product, ordinary harmonic polynomials are orthogonal to symmetric polynomials and a result analogous to Lemma \ref{adjoints} follows.
\end{rmk}

\begin{prop}
\label{harmonics}
Let \(\scalars[S_N/S_{k,N-k}]\)
and \(\scalars[S_k]\) denote the 
free \(\scalars\)-modules generated by the sets
\(S_N/S_{k,N-k}\) and \(S_k\) respectively.
The Cartan-harmonic polynomials contains a submodule isomorphic to
\[
\bigoplus_{k=0}^N \scalars[S_N/S_{k,N-k}] \otimes \scalars[S_k].
\]
\end{prop}
\begin{proof}
For any non-negative \(k\), let \(\scalars_{1^N}^k \subseteq \scalars_{1^N}\) denote the submodule generated by monomials each involving a collection of exactly \(k\) distinct variables.
The \(\scalars\)-module isomorphism
\(
\scalars_{1^N} \cong \bigoplus_{k=0}^N \scalars_{1^N}^k
\)
is equivariant with respect to \(S_N\) and \(\Wittt_{1}^*\).
For any sequence of \(k\) distinct variables \(x_{i_1}, \dots, x_{i_k}\), the submodules
\[x_{i_1} \cdots x_{i_k}\scalars[x_{i_1}, \dots, x_{i_k}]\subseteq \scalars_{1^N}^k\] are fixed by \(\Wittt_{1}^*\) and span \(\scalars^k_{1^N}\).
The space of leading coefficients is generated as an \(S_N\)-module by \(e_k:= x_{1} \dots x_k\) and is isomorphic to 
\(S_N/S_{k,N-k}\).
This yields a \(\Wittt_{1}^*\)-equivariant isomorphism
\[
\scalars^k_{1^N} \cong \scalars[S_N/S_{k,N-k}] \otimes e_k\scalars_{1^k},
\]
with a trivial action on the left.

By Chevalley's theorem, there is an \(S_k\)-module isomorphism
\[
e_k\scalars_{1^k} \cong \scalars[S_k] \otimes e_k \scalars_k.
\]
The graded dimension of a free \(\Wittt_{1}\)-module on one generator agrees with the graded dimension of \(\scalars_k\) in the first \(k\) degrees.
Thus there are at least \(k!\) \(\Wittt_{1}\)-generators and, by minimality, Cartan-harmonic polynomials.
The action of \(\Wittt_{1}\) commutes with \(S_N\) so the generators must form a copy of \(\scalars[S_k]\).
\end{proof}

\begin{cor}
\[
\dim(H_N)/N! \geq \frac{1}{N!} \sum_{k=0}^N {N \choose k} k! = \sum_{k=0}^N \frac{1}{k!} \xrightarrow{N \to \infty} e.
\]
\end{cor}

\begin{lem}
The only Cartan-harmonic polynomial in \(e_k\scalars_k\) is \(e_k\) itself. Thus \(e_k\scalars_k\) is indecomposible as a \(\Wittt_{0}\)-module.
\end{lem}
\begin{proof}
Let \(f\) be a homogeneous symmetric polynomial in  \(e_k\scalars_k\) such that the highest degree of \(x_1\) is \(m\).
Suppose \(f\) contains a monomial \[v := x_1^{m} x_2^{a_2} \cdots x_k^{a_k},\]
and thus \(\wgen_{m-1}^*f\) contains the monomial
\[
w:= x_1 x_2^{a_2} \cdots x_k^{a_k}.
\]
If there is another monomial \(u \neq v\) in \(f\) such that \(\wgen_{m-1}^*u\) contains the monomial \(w\), \(u\) must be of the form
\(
x_1 x_2^{a_2} \cdots x_i^{a_i+m-1} \cdots x_k^{a_k}.
\)
The exponent \(a_i + m-1\) must be at least \(m\), otherwise it would vanish under \(\wgen_{m-1}^*\)
and the maximality of \(m\) forces \(a_i+m-1\) to be at most \(m\).
Together these imply \(a_i = 1\).
Thus the exponents of \(v\) and \(u\) are the same up to a permutation so the symmetry of \(f\) guarantees that they have the same coefficient.
Finally, this implies the coefficient of \(w\) in \(\wgen_{m-1}^*f\) is non-zero preventing \(f\) from being Cartan-harmonic.
\end{proof}

\begin{cor}
The base ring \(\scalars_N\) is generated by the elementary symmetric polynomials as a \(\Wittt_{0}\)-module.
\end{cor}

\begin{conj}
The submodule in Proposition \ref{harmonics} is the entire space of Cartan-harmonic polymomials.
\end{conj}

\begin{rmk}
To prove the conjecture, it is sufficient to show that \(\scalars_k\)-linear combinations of the Cartan-harmonic copies of \(\scalars[S_k]\) are never Cartan-harmonic. When \(N=2\) all that has to be shown is that there are no Cartan-harmonic anti-symmetric polynomials. In general this is more complicated because higher dimensional representations of \(S_k\) appear with multiplicity in different degrees.
In \cite{Nishiyama_1996}, Cartan-harmonic polynomials were computed explicitly for \(N=3,4\) with dimensions agreeing with the conjecture.
\end{rmk}

\begin{lem}
The polynomial \((x-y)xy\) is the only Cartan-harmonic polynomial in the submodule \((x-y)xy\scalars_2\).
\end{lem}
\begin{proof}
Let \(f\) be a degree \(n\) homogeneous element of \((x-y)xy\scalars_2\).
Any such antisymmetric polynomial can be written as a sum
\[
f = \sum_{i = 1}^{\lceil n/2 \rceil} a_i (x^{n-i}y^{i} - x^{i}y^{n-i}).
\]
then \(\wgen_1^*f\) can be written as the sum
\begin{align}
\wgen_1^*f &= \sum_{i = 1}^{\lceil n/2 \rceil} a_i (n-i)(n-i-1)(x^{n-i-1}y^{i} - x^{i}y^{n-i-1}) \\
&+\sum_{i = 1}^{\lceil n/2 \rceil} a_i i(i-1)(x^{n-i}y^{i-1} - x^{i-1}y^{n-i}).
\end{align}
Let \(k\) be the largest \(k\) such that \(a_k \neq 0\).

If \(n\) is even, then \(n-k-1 \neq k\) so \((x^{n-k-1}y^k - x^k y^{n-k-1}) \neq 0\).
The coefficient of \((x^{n-k-1}y^{k} - x^k y^{n-k-1})\) in \(\wgen_1^* f\) is precisely \(a_k(n-k)(n-k-1)\).
Since \(f\) is divisible by \((x-y)xy\), \(k\) cannot be equal to \(n\) and the coefficient is non-zero.
Thus \(\wgen_1^*f \neq 0\), so \(f\) is not Cartan harmonic.

If \(n\) is odd and \(k\neq(n-1)/2\), then the same argument applies since \(n-k-1 \neq k\).
Now assume that \(k = (n-1)/2\).
In this case \((x^{n-k-1}y^{k} - x^ky^{n-k-1}) = 0\) so consider the next term in the expansion.
The coefficient of 
\((x^{k+1}y^{k-1} - x^{k-1}x^{k+1})\) in \(\wgen_1^*f\) is
\[
a_{k-1}(k+2)(k+1) + a_{k} k(k-1).
\]
The coefficient of
\((x^{k}y^{k-1}-x^{k-1}y^k)\) in \(\wgen_2^* f\)  is
\[
a_{k-1}(k+2)(k+1)k - a_k (k+1)k(k-1).
\]
Unless \(a_{k-1} = 0\) or \(k=1\), 
coefficients cannot both be zero so at least one of \(\wgen_1^*f\) and \(\wgen_2^*f\) is non-zero.
As long as \(f \neq (x-y)xy\), \(k > 1\).
Thus if \(a_{k-1} = 0\) then \(\wgen_1^*f\) is guaranteed to be non-zero.
We conclude that no \(f \in (x-y)xy \scalars_2\) is Cartan-harmonic except for \((x-y)xy\) itself.
\end{proof}

\begin{prop}
As a \(\Witt\)-module, \(\scalars_{1^N}\) is a direct sum of \(\dim (H_N) -1\) indecomposible \(\Witt\)-submodules generated by highest weight vectors.
\end{prop}
\begin{proof}
The operator \(\wgen_{-1}\) is adjoint to multiplication by \(p_1\)
so in analogy with Lemma (\ref{adjoints}) there is a ring isomorphism
\[\scalars_{1^N} \cong (\ker \wgen_{-1})[p_1].\]
Cartan-harmonic generators which are not killed by \(\wgen_{-1}\)
are in \(\Witt\)-extensions with other \(\Wittt_{1}\)-submodules.
To produce split \(\Witt\)-extensions, the Cartan-harmonic \(\Wittt_{1}\)-generators should be replaced with the constant terms in their \(p_1\) polynomial expansions via the quotient map
\[R_{1^N} \twoheadrightarrow R_{1^N}/p_1R_{1^N} \cong \ker \wgen_{-1}.\]

The polynomial \(p_1\) is Cartan-harmonic and thus cannot be replaced by a highest weight generator.
No other Cartan-harmonic polynomial is divisible by \(p_1\) and therefore a highest weight replacement is possible for all other generators.
Suppose \(q \cdot p_1\) were Cartan-harmonic for some homogeneous polynomial \(q\).
The commutation relation in \(\Witt^*\) implies
\[
\wgen_1^* (p_1q) =  2 \wgen_{0}^* q + p_1 \wgen_1^*q.
\]
Since the operator  \(2\wgen_0^* = 2\wgen_0\) multiplies polynomials by their degree,
\(\wgen_1^*(qp_1)\) vanishes 
when \[p_1 \wgen_1^* q = - (\deg q) q.\]
This only happens if \(q\) is \(1\) or divisible by \(p_1\).
Repeating the commutation relation for \(p_1 \wgen_1^*q\) shows that if \(q \neq 1\), it would be infinitely divisible by \(p_1\) which is impossible.
Thus \(\wgen_1^*(qp_1) \neq 0\) for any \(q \neq 1\)
and there are no Cartan-harmonic polynomial multiples of \(p_1\).
In particular, this implies that all other Cartan-harmonic generators have non-zero images under the quotient map onto \(\ker L_{-1}\).

The \(\Wittt_{1}\)-submodule generated by \(p_1\) is in extension with the trivial submodule generated by \(1\). Together these form a \(\Witt\)-direct summand.
The remaining \(\Wittt_{1}\)-generators can be replaced with highest weight generators forming \(\Witt\)-direct summands.
Since each highest weight replacement is accomplished by subtracting terms belonging to other \(\Wittt_1\)-modules,
the corresponding quotient modules are isomorphic as \(\Wittt_{1}\)-modules.
In particular they are indecomposible.
\end{proof}

\begin{cor}
In the case \(N=2\), \(\ker \wgen_{-1} \cong \scalars[(x-y)]\) and the highest weight \(\Witt\)-generators are \(1, (x-y), (x-y)^2,\) and \((x-y)^3\).
\end{cor}

Universal Lee deformation is recovered from \(\scalars_2\)-linear homology by setting \(E_1 = x+y\) to zero.
The above results suggest that the Lee deformation can be identified as the space of highest weight vectors \(\scalars[(x-y)]\).
In fact the functoriality of the Witt representation shows that this identification is natural with respect to link cobordisms.
Thus the following conjecture is proven when restricted to unlinks and theta foams.

\begin{conj}
If \(\lambda = \mu = 0\) and \(s = 1/2\), the kernel of \(\wgen_{-1}\) is naturally isomorphic to the Lee deformation of Khovanov homology.
\end{conj}

Bar-Natan homology, that is the \((x^2-hx)\)-theory also appears in this picture by setting \(E_2 = xy\) to zero.
The ideal generated by \(xy\) is not invariant under \(\wgen_{-1}\), but it is a \(\Wittt_{0}\)-submodule.
Thus the Bar-Natan homology of the module \(\scalars[x,y]\) appears as the first two \(\Wittt_{0}\)-submodules with one dimensional weight spaces.




\subsection{Twisted representations}

Analyzing the twisted representations of the form \(t_{\mathbf{a}}^{-1/2} \scalars_{\mathbf{a}}\) is harder since the inner product and adjoint Witt operators do not have obvious extensions even in the simplest case \(N=2\).
However, the highest weight structure of \(\scalars_{1^2}\) is simple.
It was shown in \cite{Roz23} that the highest weight submodule of \((x-y)^{\lambda}\scalars_{1^2}\) is \((x-y)^{\lambda}\scalars[(x-y)]\) where \(\lambda \in \frac{1}{2}\Z\).

Direct computation shows that
\begin{align}
\wgen_1^2 (x-y)^{\lambda} &= \lambda \wgen_1 \left((x+y) (x-y)^{\lambda}\right) \\
&= \lambda \left( (x^2+y^2) + \lambda(x+y)^2\right)(x-y)^{\lambda}\\
&= \lambda \left( (1+\lambda)(x^2+y^2) + 2\lambda xy)\right)(x-y)^{\lambda}.
\end{align}
On the other hand,
\begin{align}
\wgen_2 (x-y)^{\lambda} = \lambda (x^2+xy+y^2)(x-y)^{\lambda}.
\end{align}
These are linear combinations precisely when 
\(2\lambda^2 = \lambda(1+\lambda)\).
That is when \(\lambda = 0, 1\).

This implies that unless \(\lambda=0,1\), the \(\sld\)-submodule generated by \((x-y)^\lambda\) is in \(\Witt\)-extension with the one generated by \((x-y)^{\lambda+2}\).
In particular, this shows that the homology of the unknot, \((x-y)^{-1/2}\scalars_{1^2}\) is a direct sum of the two \(\Witt\)-submodules generated by \((x-y)^{-1/2}\) and \((x-y)^{1/2}\).

%% file: SF_forktwist.tex
\begin{lem}
\label{forktwist}
There are the following isomorphisms in the relative homotopy category:
\begin{equation} 
\NB{\tikz[font=\tiny, scale=-0.6, yscale=1]{\input{sch/forktwist_1}}} \cong q^2\NB{\tikz[font=\tiny, scale=-0.6]{\input{sch/forktwist_2c}}} \hspace{3mm}\text{and}\hspace{3mm} \NB{\tikz[font=\tiny, scale=-0.6, xscale=-1]{\input{sch/forktwist_1}}} \cong q^{-2}\NB{\tikz[font=\tiny, scale=-0.6]{\input{sch/forktwist_2d}}}.
\end{equation}
\end{lem}

\begin{proof}
    The construction of both isomorphisms uses the same ideas so we only present the proof for the first isomorphism. Thanks to Lemma \ref{lemm_cross2} the existence of the following isomorphism is enough:
\begin{equation} \label{forktwist1}
\NB{\tikz[font=\tiny, scale=-0.6, yscale=1]{\input{sch/forktwist_1b}}} \cong q^{2}\NB{\tikz[font=\tiny, scale=-0.6]{\input{sch/forktwist_2b}}}.
\end{equation}

The first complex is defined to be:
\[
\NB{\tikz[font=\tiny, scale=-0.6]{\input{sch/forktwist_1b}}} =
\NB{
	\tikz[xscale = 3, yscale = 3]{
		\node (AT) at (0,0) {$q\NB{\tikz[font=\tiny, scale=-0.6]{\input{sch/forktwist_3}}}$};
		\node (CT) at (1.5,0) {$\NB{\tikz[font=\tiny, scale=-0.6]{\input{sch/forktwist_2}}}$};
		\draw[->] (AT) -- (CT) node[pos=0.5, above]{$\mapunzip$};
	}
}.
\]

These two complexes fit in a short exact sequence of $\Witt$-equivariant complexes that splits after forgetting the $\Witt$ action:

\[
\NB{\tikz[xscale = 3, yscale = 3]{
	\node (i0) at (0, 0) {$q\NB{\tikz[font= \tiny, scale=-.6]{\input{sch/forktwist_3}}}$};
	\node (i1) at (1.5, 0) { $\NB{\tikz[font= \tiny,scale=-0.6]{\input{sch/forktwist_2}}}$ };
	\node (i2) at (0, 1.1) {$q^2\NB{\tikz[font= \tiny, scale=-0.6]{\input{sch/forktwist_2b}}}$};
	\node (i3) at (1.5, -1.1) {$\NB{\tikz[font= \tiny, scale=-0.6]{\input{sch/forktwist_2}}}$};
	\node (i4) at (0, -1.1) {\NB{\tikz[font= \tiny, scale=-0.6]{\input{sch/forktwist_2}}}};
	\draw[->] (i0) -- (i1)  node[pos=0.5, above] {$\mapunzip$};
	\draw[->] (i0) -- (i2) node[pos=0.5, left] {$\mapdcap$} ;
	\draw[->] (i4) -- (i0)  node[pos=0.5, left] {$\mapdcup$} ;
	\draw[->] (i3) -- (i1)  node[pos=0.5, right] {$\Id$} ;
	\draw[->] (i4) -- (i3)  node[pos=0.5, above] {$\Id$};
}}.
\]

The vertical maps are given by the digon cup and digon cap.  
Since the bottom complex is contractible, Corollary \ref{lemme_RHC} gives the desired isomorphism.

\end{proof}

%% file: sch/forktwist_1.tex
\begin{scope}[xshift = 2.5cm]
	\draw[->-][double] (0.5,2) -- (0.5,1);
	\draw[->] (1, 0.5) arc (0:180:0.5) -- (0,0.5) .. controls +(0,-0.5) and +(0, 0.5)
		.. (1, -0.5) -- (1,-1);
		\fill[white] (0.5, 0) circle (1mm);
	\draw[<-] (0, -1) -- (0, -0.5) .. controls +(0,0.5) and +(0, -0.5)
		.. (1, 0.5);
\end{scope}

%% file: sch/forktwist_2c.tex
\begin{scope}[xshift = 2.5cm]
	\draw[<->] (1,-1) -- (1, 0.5) arc (0:180:0.5) -- (0,0.5) -- ( 0,-1);
	\draw[double] (0.5,2) -- (0.5,1) node[pos=1, draw=\crd, fill= white, circle , inner sep=1.5pt]{} node[pos=1, yshift=-1mm, right]{$1$};
\end{scope}

%% file: sch/forktwist_2d.tex
\begin{scope}[xshift = 2.5cm]
	\draw[<->] (1,-1) -- (1, 0.5) arc (0:180:0.5) -- (0,0.5) -- ( 0,-1);
	\draw[double] (0.5,2) -- (0.5,1) node[pos=1, draw=\crd, fill= white, circle , inner sep=1.5pt]{} node[pos=1, yshift=-1mm, right]{$-1$};
\end{scope}

%% file: sch/forktwist_2b.tex
\begin{scope}[xshift = 2.5cm]
	\draw[<->] (1,-1) -- (1, 0.5) arc (0:180:0.5) -- (0,0.5) -- ( 0,-1);
	\draw[double] (0.5,2) -- (0.5,1) node[pos=1, draw=\crd, fill= white, circle , inner sep=1.5pt]{} node[pos=1, yshift=-1mm, right]{$s$};
\end{scope}

%% file: sch/forktwist_3.tex
\begin{scope}[xshift = 2.5cm]
	\coordinate (bl) at (0,  -1.5);
	\coordinate (br) at ( 1,  -1.5);
	\coordinate (bm) at (   0.5,-0.8);
	\coordinate (tl) at (0,   0.5);
	\coordinate (tr) at ( 1,   0.5);
	\coordinate (tm) at (   0.5, -0.2);

	\draw[<-]  (bl) .. controls +( 0, 0.5) and +(0,0) .. (bm);
	\draw[<-]  (br) .. controls +( 0, 0.5) and +(0,0) .. (bm);
	\draw[->-]  (0.5,1) .. controls +( -0.5, 0) and +(0,0) ..  (tl) .. controls +( 0, -0.5) and +(0,0) .. (tm)
		node[pos= 0.35 , draw=\crd, fill=white,regular polygon, regular polygon sides=3,inner sep=0.9pt] {} node[pos=0.3, right]{$\mu$};
	\draw[->-]  (0.5,1) .. controls +(0.5,0) and +(0,0) .. (tr) .. controls +( 0, -0.5) and +(0,0) .. (tm)
		node[pos= 0.35 , draw=\crd, fill=white,regular polygon, regular polygon sides=3,inner sep=0.9pt] {} node[pos=0.3, left]{$\lambda$};
	\draw [double] (bm) -- (tm)
		node[pos=0, draw=\crd, fill= white, circle , inner sep=1.5pt]{} node[pos=0, right]{$s$};
	\draw[double] (0.5,2) -- (0.5,1)
		node[pos=1, draw=\crd, fill= white, circle , inner sep=1.5pt]{} node[pos=1,yshift=-1mm, right]{$- \overline{s}$};
\end{scope}

%% file: sch/forktwist_2.tex
\begin{scope}[xshift = 2.5cm]
	\draw[<->] (1,-1) -- (1, 0.5) arc (0:180:0.5) -- (0,0.5) -- ( 0,-1);
	\draw[double] (0.5,2) -- (0.5,1) node[pos=1, draw=\crd, fill= white, circle , inner sep=1.5pt]{} node[pos=1,yshift=-1mm, right]{$- \overline{s}$};
\end{scope}

%% file: SF_toruslinks.tex
Using a standard argument we compute the action on \(T_{2,m}\):
\[
\NB{\tikz[font=\tiny, scale=0.6]{\input{sch/toruslink}}}.
\]
The following complex remains after resolving the top crossing only:
\[
\NB{\tikz[font=\tiny, scale=0.6]{\input{sch/toruslink}}}
=
\NB{
	\tikz[xscale = 3, yscale = 3]{
		\node (AT) at (0,0) {$q~~\NB{\tikz[font=\tiny, scale=0.6]{\input{sch/toruslink_2}}}$};
		\node (CT) at (1.5,0) {$~~\NB{\tikz[font=\tiny, scale=0.6]{\input{sch/toruslink_3}}}$};
		\draw[->] (AT) -- (CT) node[pos=0.5, above]{$\mapunzip$};
	}
}.
\]
Applying fork twists \ref{forktwist} to the left term and induction on the right term, the complex below is obtained.
\[\NB{
	\tikz[xscale = 3, yscale = 3]{
		\node (AT) at (0,0) {$q^{2m-1}\, \NB{\tikz[font=\tiny, scale=0.6]{\input{sch/toruslink_5}}}$};
		\node (BT) at (1.5,0) {$q^{2m-3}\, \NB{\tikz[font=\tiny, scale=0.6]{\input{sch/toruslink_5b}}}$};
		\node (CT) at (2.5,0) {\(\cdots\)};
		\node (AB) at (0,-0.75) {\(\cdots\)};
		\node (BB) at (1,-0.75) {$q^{3}\, \NB{\tikz[font=\tiny, scale=0.6]{\input{sch/toruslink_5c}}}$};
		\node (CB) at (2,-0.75) {$\,0$};
		\node (DB) at (3,-0.75) {$q^{1-N}\NB{\tikz[font=\tiny, scale=0.6]{\input{sch/hopflink_5}}}$};
	\draw[->] (AT) -- (BT) node[pos=0.5, above]{$d^{1}$};
	\draw[->] (BT) -- (CT) node[pos=0.5, above]{$d^{2}$};
	\draw[->] (AB) -- (BB) node[pos=0.5, above]{$d^{m-2}$};
	\draw[->] (BB) -- (CB);
	\draw[->] (CB) -- (DB);
	}
}.\]
The differentials alternate between the following two maps:
\begin{align}
d^{m-2i} &=
\NB{\tikz[font=\tiny, scale=1]{\input{sch/identity_X}}}
-
\NB{\tikz[font=\tiny, scale=1]{\input{sch/identity_Y}}}\\
d^{m-2i-1} &=
\NB{\tikz[font=\tiny, scale=1]{\input{sch/identity_X}}}
-
\NB{\tikz[font=\tiny, scale=1]{\input{sch/identity_X_1}}} = 0
\end{align}

The results of section \ref{sec_EX} make the homology computation easy:
\begin{thm}
\label{homaspoly}
As a graded \(\Witt \# \scalars_N\)-module 
\begin{align*}
\KRW_N^m(T_{2,m})     &\cong \prod_{1<i\leq N}(x_1-x_i)^{-1}\scalars_{1,N-1}, \\
\KRW_N^{m-j}(T_{2,m}) &\cong 0 \text{ for all odd } j, \\
\KRW_N^{m-j}(T_{2,m}) &\cong (x_1-x_2)^{j-2}\scalars_{1,1,N-2}/(x-y)^{j-1}\scalars_{1,1,N-1} \text{ for all even } 0 < j < m,\\
\KRW_N^{0}(T_{2,m})   &\cong \prod_{1<i\leq N}(x_1-x_i)^{m-2}\scalars_{1,N-1} \text{ if } m \text{ is even}.
\end{align*}
\end{thm}

\subsection{The case \(N=2\)}

The first goal is to describe all \(\Witt \# \scalars_2\)-equivariant maps \(\uKRW(\varnothing) \to \uKRW(T_{2,m}^{\lambda})\).
Unfortunately, the techniques of Section \ref{sec_EX} do not carry over well to the twisted case because the inner product is not well defined when twists are negative or fractional.
Even when twists are positive, the adjoints are not always well defined.
Thus \textit{ad hoc} methods and direct computation are necessary.

The twisted homology can be described easily when \(N=2\).
We present the unframed version to avoid issues with counting Reidemeister moves.
\begin{cor}
As graded \(\Witt \# \scalars_N\)-modules
\begin{align*}
\uKRW_2^0(T_{2,m}^{\lambda})     &\cong (x-y)^{\frac{\lambda+m}{2}-1}\scalars[x,y], \\
\uKRW_2^{-j}(T_{2,m}^{\lambda}) &\cong 0 \text{ for all odd } j, \\
\uKRW_2^{-j}(T_{2,m}^{\lambda}) &\cong (x-y)^{\frac{\lambda+m}{2}+j-2}\scalars[x,y]/(x-y)^{\frac{\lambda+m}{2}+j-1}\scalars[x,y] \\&\quad\quad \text{ for all even } 0 < j < m,\\
\uKRW_2^{-m}(T_{2,m}^{\lambda})   &\cong (x-y)^{\frac{\lambda+m}{2}+m-2}\scalars[x,y] \text{ if } m \text{ is even}.
\end{align*}
\end{cor}

Any such \(\scalars_N\)-linear map is determined by the image of \(1\) and \(\Witt\)-equivariance is guaranteed if the image is trivial module.
When \(N=2\) the decomposition above agrees with the \(\sld\) decomposition in \cite{Roz23} of the framed homology.
It was shown that the only highest weight vectors are powers of \((x-y)\) and in particular the only trivial submodule is generated by \((x-y)^0\).
Therefore the only maps which can be induced by cobordisms \(C\) are scalar multiples of the inclusion of symmetric polynomials into \((x-y)^{(\lambda+m)/2-1}\scalars[x,y]\).
These only exist when \(\lambda+m\) is even and \(\chi(C)/2 = \lambda \leq 2 - m\).
When \(m\) is odd, so that \(T_{2,m}\) is a knot, \(\chi(C) = 1-2g\) and \(m-3/2 \leq m - 1 \leq g \) recovering the genus bound from the \(s\)-invariant.



%% file: sch/toruslink_2.tex
\begin{scope}[xshift = 2.5cm]

	\draw[-] (0.5, -2.5) ..controls +(0,0.3) and +(0,-0.3) .. (-0.5,-1.5);
	\fill[white] (0,-2) circle (2mm);
	\draw[-] (-0.5, -2.5) ..controls +(0,0.3) and +(0,-0.3) .. (0.5,-1.5);

	\node at (0,-1) {$\cdots$};

	\draw[-] (0.5, 0.5) arc (0:180:0.5) -- (-0.5,0.5) .. controls +(0,-0.5) and +(0, 0.5) .. (0.5, -0.5);
	\fill[white] (0, 0) circle (1mm);
	\draw[-] (-0.5, -0.5) .. controls +(0,0.5) and +(0, -0.5) .. (0.5, 0.5);

	\draw[-]  (0.5, 2) .. controls +( 0, -0.5) and +(0,0) .. (0,1.5)
		node[pos= 0.35 , draw=\crd, fill=white,regular polygon, regular polygon sides=3,inner sep=0.9pt] {} node[pos=0.3, right]{$\mu$};
	\draw[-]  (-0.5, 2) .. controls +( 0, -0.5) and +(0,0) .. (0,1.5)
		node[pos= 0.35 , draw=\crd, fill=white,regular polygon, regular polygon sides=3,inner sep=0.9pt] {} node[pos=0.3, left]{$\lambda$};
	\draw[-][double] (0,1.5) -- (0,1)
		node[pos=1, draw=\crd, fill= white, circle , inner sep=1.5pt]{} node[pos=1,yshift=1mm, right]{$s$};

	\draw[->-] (-0.5, 2) arc (-180:0:-0.5)  -- (-1.5,-2.5) arc (-180:0:0.5)  ;
	\draw[->-] (0.5, 2) arc (180:0:0.5)  -- (1.5,-2.5) arc (0:-180:0.5)  ;
\end{scope}

%% file: sch/toruslink_5.tex
\begin{scope}[xshift = 2.5cm]

	\draw[-] (0.5, 0.5) arc (0:180:0.5) -- (-0.5,0.5);


	\draw[-]  (0.5, 2) .. controls +( 0, -0.5) and +(0,0) .. (0,1.5)
		node[pos= 0.35 , draw=\crd, fill=white,regular polygon, regular polygon sides=3,inner sep=0.9pt] {} node[pos=0.3, right]{$\mu$};
	\draw[-]  (-0.5, 2) .. controls +( 0, -0.5) and +(0,0) .. (0,1.5)
		node[pos= 0.35 , draw=\crd, fill=white,regular polygon, regular polygon sides=3,inner sep=0.9pt] {} node[pos=0.3, left]{$\lambda$};
	\draw[-][double] (0,1.5) -- (0,1)
		node[pos=1, draw=\crd, fill= white, circle , inner sep=1.5pt]{} node[pos=1,yshift=1mm, right]{$s+m-1$};

	\draw[->-] (-0.5, 2) arc (-180:0:-0.5)  -- (-1.5,0.5) arc (-180:0:0.5)  ;
	\draw[->-] (0.5, 2) arc (180:0:0.5)  -- (1.5,0.5) arc (0:-180:0.5)  ;
\end{scope}

%% file: sch/toruslink_5b.tex
\begin{scope}[xshift = 2.5cm]

	\draw[-] (0.5, 0.5) arc (0:180:0.5) -- (-0.5,0.5);


	\draw[-]  (0.5, 2) .. controls +( 0, -0.5) and +(0,0) .. (0,1.5)
		node[pos= 0.35 , draw=\crd, fill=white,regular polygon, regular polygon sides=3,inner sep=0.9pt] {} node[pos=0.3, right]{$\mu$};
	\draw[-]  (-0.5, 2) .. controls +( 0, -0.5) and +(0,0) .. (0,1.5)
		node[pos= 0.35 , draw=\crd, fill=white,regular polygon, regular polygon sides=3,inner sep=0.9pt] {} node[pos=0.3, left]{$\lambda$};
	\draw[-][double] (0,1.5) -- (0,1)
		node[pos=1, draw=\crd, fill= white, circle , inner sep=1.5pt]{} node[pos=1,yshift=1mm, right]{$s+m-2$};

	\draw[->-] (-0.5, 2) arc (-180:0:-0.5)  -- (-1.5,0.5) arc (-180:0:0.5)  ;
	\draw[->-] (0.5, 2) arc (180:0:0.5)  -- (1.5,0.5) arc (0:-180:0.5)  ;
\end{scope}

%% file: sch/toruslink_5c.tex
\begin{scope}[xshift = 2.5cm]

	\draw[-] (0.5, 0.5) arc (0:180:0.5) -- (-0.5,0.5);

	\draw[-]  (0.5, 2) .. controls +( 0, -0.5) and +(0,0) .. (0,1.5)
		node[pos= 0.35 , draw=\crd, fill=white,regular polygon, regular polygon sides=3,inner sep=0.9pt] {} node[pos=0.3, right]{$\mu$};
	\draw[-]  (-0.5, 2) .. controls +( 0, -0.5) and +(0,0) .. (0,1.5)
		node[pos= 0.35 , draw=\crd, fill=white,regular polygon, regular polygon sides=3,inner sep=0.9pt] {} node[pos=0.3, left]{$\lambda$};
	\draw[-][double] (0,1.5) -- (0,1)
		node[pos=1, draw=\crd, fill= white, circle , inner sep=1.5pt]{} node[pos=1,yshift=1mm, right]{$s+1$};

	\draw[->-] (-0.5, 2) arc (-180:0:-0.5)  -- (-1.5,0.5) arc (-180:0:0.5)  ;
	\draw[->-] (0.5, 2) arc (180:0:0.5)  -- (1.5,0.5) arc (0:-180:0.5)  ;
\end{scope}

%% file: sch/hopflink_5.tex
\begin{scope}[xshift = 2.5cm]
	\draw (-0.5, 2) arc (-180:0:-0.5)  -- (-1.5,0.5) arc (-180:0:0.5)  -- (-0.5,2)
	      node[draw=\crd, fill=\csrd, circle, inner sep=1.5pt, pos=0]{} node[pos=0, right]{$-\frac{1}{2}$};
\end{scope}

%% file: sch/identity_X.tex
\begin{scope}
  \begin{scope}
    \coordinate (L1) at (0.2,0.4);
    \coordinate (L2) at (0,0);
    \coordinate (R1) at (2.2,0.4);
    \coordinate (R2) at (2,0);
    \coordinate (ML) at (0.6, 0.2);
    \coordinate (MR) at (1.6, 0.2);
    \draw[double] (ML) -- (MR) node[above, midway] {};
    \draw (MR) .. controls +(0, 0) and +(-0.3,0) .. (R1);
    \draw (MR) .. controls +(0, 0) and +(-0.3,0) .. (R2);
    \draw (L1) .. controls +( 0.3, 0) and +(0,0) .. (ML);
    \draw (L2) .. controls +( 0.3, 0) and +(0,0) .. (ML);
    \draw[dotted] (L2) ..controls +(-.8,-.8) and +(.6,-.6) .. (R2);
    \draw[dotted] (L1) ..controls +(-.4,.4) and +(.6,.6) .. (R1);
  \end{scope}  
 \begin{scope}[yshift = -1cm]
    \coordinate (L1B) at (0.2,0.4);
    \coordinate (L2B) at (0,0);
    \coordinate (R1B) at (2.2,0.4);
    \coordinate (R2B) at (2,0);
    \coordinate (MLB) at (0.6, 0.2);
    \coordinate (MRB) at (1.6, 0.2);
    \coordinate (R3B) at (2.2,.5);
    \draw[double] (MLB) -- (MRB) node[above, midway] {};
    \draw (MRB) .. controls +(0, 0) and +(-0.3,0) .. (R1B) ;
    \draw (MRB) .. controls +(0, 0) and +(-0.3,0) .. (R2B) node[above, pos =0.8] {$\deco_1$};
    \draw (L1B) .. controls +( 0.3, 0) and +(0,0) .. (MLB);
    \draw (L2B) .. controls +( 0.3, 0) and +(0,0) .. (MLB);
    \draw[dotted] (L2B) ..controls +(-.8,-.8) and +(.6,-.6) .. (R2B);
    \draw[dotted] (R1B) ..controls +(.1,.1) and +(.1,.1) .. (R3B);
 \end{scope}  
  \draw (R1) -- (R1B);
  \draw (R2) -- (R2B);
  \draw (L1) -- (L1B);
  \draw (L2) -- (L2B);
  \draw (ML) -- (MLB);
  \draw (MR) -- (MRB);
\end{scope}

%% file: sch/identity_Y.tex
\begin{scope}
  \begin{scope}
    \coordinate (L1) at (0.2,0.4);
    \coordinate (L2) at (0,0);
    \coordinate (R1) at (2.2,0.4);
    \coordinate (R2) at (2,0);
    \coordinate (ML) at (0.6, 0.2);
    \coordinate (MR) at (1.6, 0.2);
    
    \draw[double] (ML) -- (MR) node[above, midway] {};
    \draw (MR) .. controls +(0, 0) and +(-0.3,0) .. (R1);
    \draw (MR) .. controls +(0, 0) and +(-0.3,0) .. (R2);
    \draw (L1) .. controls +( 0.3, 0) and +(0,0) .. (ML) node[below, pos =0.3] {$\deco_1$};
    \draw (L2) .. controls +( 0.3, 0) and +(0,0) .. (ML);
    \draw[dotted] (L2) ..controls +(-.8,-.8) and +(.6,-.6) .. (R2);
    \draw[dotted] (L1) ..controls +(-.4,.4) and +(.6,.6) .. (R1);
  \end{scope}  
 \begin{scope}[yshift = -1cm]
    \coordinate (L1B) at (0.2,0.4);
    \coordinate (L2B) at (0,0);
    \coordinate (R1B) at (2.2,0.4);
    \coordinate (R2B) at (2,0);
    \coordinate (MLB) at (0.6, 0.2);
    \coordinate (MRB) at (1.6, 0.2);
    \coordinate (R3B) at (2.2,.5);
    \draw[double] (MLB) -- (MRB) node[above, midway] {};
    \draw (MRB) .. controls +(0, 0) and +(-0.3,0) .. (R1B) ;
    \draw (MRB) .. controls +(0, 0) and +(-0.3,0) .. (R2B);
    \draw (L1B) .. controls +( 0.3, 0) and +(0,0) .. (MLB);
    \draw (L2B) .. controls +( 0.3, 0) and +(0,0) .. (MLB);
    \draw[dotted] (L2B) ..controls +(-.8,-.8) and +(.6,-.6) .. (R2B);
    \draw[dotted] (R1B) ..controls +(.1,.1) and +(.1,.1) .. (R3B);
 \end{scope}  
  \draw (R1) -- (R1B);
  \draw (R2) -- (R2B);
  \draw (L1) -- (L1B);
  \draw (L2) -- (L2B);
  \draw (ML) -- (MLB);
  \draw (MR) -- (MRB);
\end{scope}

%% file: sch/identity_X_1.tex
\begin{scope}
  \begin{scope}
    \coordinate (L1) at (0.2,0.4);
    \coordinate (L2) at (0,0);
    \coordinate (R1) at (2.2,0.4);
    \coordinate (R2) at (2,0);
    \coordinate (ML) at (0.6, 0.2);
    \coordinate (MR) at (1.6, 0.2);
    \draw[double] (ML) -- (MR) node[above, midway] {};
    \draw (MR) .. controls +(0, 0) and +(-0.3,0) .. (R1);
    \draw (MR) .. controls +(0, 0) and +(-0.3,0) .. (R2);
    \draw (L1) .. controls +( 0.3, 0) and +(0,0) .. (ML);
    \draw (L2) .. controls +( 0.3, 0) and +(0,0) .. (ML);
    \draw[dotted] (L2) ..controls +(-.8,-.8) and +(.6,-.6) .. (R2);
    \draw[dotted] (L1) ..controls +(-.4,.4) and +(.6,.6) .. (R1);
  \end{scope}  
 \begin{scope}[yshift = -1cm]
    \coordinate (L1B) at (0.2,0.4);
    \coordinate (L2B) at (0,0);
    \coordinate (R1B) at (2.2,0.4);
    \coordinate (R2B) at (2,0);
    \coordinate (MLB) at (0.6, 0.2);
    \coordinate (MRB) at (1.6, 0.2);
    \coordinate (R3B) at (2.2,.5);
    \draw[double] (MLB) -- (MRB) node[above, midway] {};
    \draw (MRB) .. controls +(0, 0) and +(-0.3,0) .. (R1B) ;
    \draw (MRB) .. controls +(0, 0) and +(-0.3,0) .. (R2B);
    \draw (L1B) .. controls +( 0.3, 0) and +(0,0) .. (MLB);
    \draw (L2B) .. controls +( 0.3, 0) and +(0,0) .. (MLB) node[above, pos =0.2] {$\deco_1$};
    \draw[dotted] (L2B) ..controls +(-.8,-.8) and +(.6,-.6) .. (R2B);
    \draw[dotted] (R1B) ..controls +(.1,.1) and +(.1,.1) .. (R3B);
 \end{scope}  
  \draw (R1) -- (R1B);
  \draw (R2) -- (R2B);
  \draw (L1) -- (L1B);
  \draw (L2) -- (L2B);
  \draw (ML) -- (MLB);
  \draw (MR) -- (MRB);
\end{scope}

%% file: SF_future.tex
\subsection{Knot-theoretic interpretation}
In the introduction and Section \ref{sec_torus} we described how to recover information about the Euler characteristics of link cobordisms from the structure of the \(\Witt\)-representation when \(N=2\).
The key feature was that the space of highest weight vectors is isomorphic to \(\scalars[(x-y)]\), possibly with a shift.
When \(N > 2\), the space of highest weight vectors is larger and its significance is unclear.
For example, when \(N=3\), the highest weight vectors in the homology of the thickness-1 unknot are generated by \((x-y)+(x-z)\) and \((x-y)(x-z)\).
The second generator corresponds to a handle attachment and its square root, which controls the red dot twists, corresponds to a Reidemeister I move.
For degree reasons, the first generator cannot correspond to a cobordism and its significance is unclear.

\subsection{Further representation theoretic analysis}
Several aspects of the \(\Witt\) structure were not analyzed in this paper.
First, the techniques for analyzing the twisted homology in \(N=2\) do not extend well and thus we do not have good structural results for the general case.
In particular, the inner product of \ref{inner-prod} does not extend in a simple way to the twisted modules.

Second, the parameters \(\lambda, \mu,\) and \(s\) were not fully analysed.
The placement of dots in the braiding complexes guarantees that in all computed examples, the homology does not depend on the parameters \(\lambda, \mu\) and \(s\).
This is not the only choice of twists for the differential and it is not clear how another choice would affect the representations.
One hint that deformations of our set up would be useful comes from considering a specialization to ordinary, non-equivariant Khovanov homology.
The ideal generated by symmetric polynomials is not preserved by the \(\Witt\)-action in characteristic zero and therefore there is no induced action on specializations.
Direct inspection of the module structure shows that even such a specialization could be defined, it would be trivial precisely because the \(\Witt\)-action commutes with the \(S_N\)-action.
However, the triangular deformations can break this symmetry and potentially lead to more fruitful specializations.

Third, the tensor products of the representations corresponding to disjoint unions were not analyzed.
Since such tensor products are taken over symmetric or partially symmetric polynomial rings, there is a natural connection to singular Soergel bimodules.
On the Soergel side, actions of \(\sld\) were essential to the work of Elias and Williamson \cite{Elias_2014}.
Elias and Williamson's \(\sld\)-action is given by multiplication operators and our operators are given by derivations making them clearly distinct.

\subsection{Connections to geometry}
The algebras of decorations for a generalized theta webs can be realized geometrically
as the \(U(N)\)-equivariant cohomology of a flag variety in \(\C^N\).
The commuting action of \(S_N\) is the action of the Weyl group and
the \(\Wittt_{1}\)-action is an integral lift of certain cohomology operations.
More specifically, the rational algebra generated by \((\wgen_k)_{k\geq 1}\) under composition is isomorphic to the Landweber-Novikov algebra of cohomology operations on complex cobordism.
The mod \(p\) reduction of a special subalgebra is isomorphic to the subalgebra of the mod \(p\) Steenrod algebra generated by power operations \cite{Wood_1997}.
Moreover, the action of \(\Witt\) on our polynomial rings agrees with the action of the Landweber-Novikov algebra on the cohomology of the corresponding flag variety \cite{Lenart_1998}.
Via the Thom isomorphism, the irreducible submodule \(x_1 \cdots x_k \scalars[x_1, \dots, x_k]^{S_k}\) is the ordinary cohomology of \(MU(k)\) and its irreducibility corresponds to classical results about cohomology operations.
The other twisted modules might also be realized via the Thom isomorphism.
Actions of the Steenrod algebra have previously been studied on categorified quantum groups and Soergel bimodules \cite{Beliakova_2018,Kitchloo_2013}.
Unfortunately, the \(L_{-1}\) action does not have a clear interpretation in this context since it has negative degree and is therefore unstable.